\documentclass[fleqn,11pt]{article}
\usepackage{amssymb,amsmath,amsfonts}
\usepackage[margin=1in]{geometry}
\usepackage{dsfont}
\usepackage{color}
\usepackage{graphicx}
\usepackage{latexsym}
\usepackage{amsmath}
\usepackage{amssymb}
\usepackage{graphics}
\usepackage[dvips]{epsfig}
\usepackage{mathrsfs}
\usepackage{leqno}

\numberwithin{equation}{section}
\newtheorem{thm}{Theorem}[section]
\newtheorem{coro}[thm]{Corollary}
\newtheorem{lem}[thm]{Lemma}
\newtheorem{prop}[thm]{Proposition}
\newtheorem{rmk}[thm]{Remark}
\newtheorem{defn}[thm]{Definition}
\newcommand{\ea}{\epsilon}
\newcommand{\ta}{\theta}

\newcommand{\da}{\delta}

\newcommand{\al}{\alpha}
\newcommand{\za}{\zeta}

\newcommand{\pl}{\partial}

\newcommand{\oa}{\omega}

\newcommand{\ba}{\beta}

\newcommand{\na}{\nabla}
\newcommand{\iy}{\infty}

\newcommand{\Da}{\Delta}
\newcommand{\Oa}{\Omega}

\newcommand{\R}{\mathbb R}

\newcommand{\idm}{{\mathscr{D}}}

\newcommand{\lt}{\left}
\newcommand{\rt}{\right}

\newcommand{\be}{\begin{equation}}
\newcommand{\bs}{\begin{split}}
\newcommand{\es}{\end{split}}
\newcommand{\ee}{\end{equation}}
\newcommand{\bee}{\begin{equation*}}
\newcommand{\eee}{\end{equation*}}

\newcommand{\ef}{\eqref}

\begin{document}
\begin{center}
\large{ \bf On the Free Surface Motion of Highly Subsonic Heat-conducting Inviscid  Flows}
\end{center}
\centerline{ Tao Luo, Huihui Zeng}
\begin{abstract}
 For a free surface problem of a highly subsonic heat-conducting inviscid flow,   motivated by a geometric approach developed by  Christodoulou and Lindblad  in the study of the free surface problem of incompressible inviscid flows,  the {\it a priori} estimates of Sobolev norms  in 2-D and 3-D are proved under the Taylor sign condition by identifying a suitable higher order energy functional.
The estimates for some geometric quantities such as the second fundamental form and the injectivity radius of the normal exponential map of the free surface are also given.
The novelty in our analysis includes dealing with the strong coupling of large variation of temperature, heat-conduction,  compressibility of fluids and the evolution of free surface, loss of symmetries of equations, and   loss of derivatives  in closing the argument which is a key feature compared with Christodoulou and Lindblad's work.
The motivation of this paper is to  contribute to the program of understanding the role played by the heat-conductivity to free surface motions of inviscid compressible flows and the behavior of such motions  when the Mach number is small.
\end{abstract}

\section{Introduction}
Fluids free surface problems have been receiving much attentions due to their physical importance and challenge in the mathematical analysis.
For incompressible inviscid flows,  the local-in-time well-posedness in Soblev spaces was obtained first in \cite{wu1, wu2} for the irrotational case,  and then in \cite{ABZ, AM, CL00,Coutand,Ebin,HaoLuo,L2,LN,SZ,SWZ,zhang} for some extensions; the global or almost global existence was achieved recently in \cite{GMS, wangxc, wu4,wu3}; and the singularity formation was proved in \cite{CCFG,wu5}. One may refer to the survey  \cite{L2} for more references.
For compressible inviscid  flows, the local-in-time well-posedness of smooth solutions was established for liquids  in  \cite{22,35}; while for gases  with  physical vacuum singularity (cf. \cite{Liu,liuyang}), the related results  can be found in  \cite{10, 10',GL, 16, 16',LXZ} for the local-in-time theories, and in \cite{hadjang1, hadjang2, LZ,HHZeng} for the global-in-time ones.
Most of the above results are either for incompressible or isentropic fluids without taking the effect of heat-conductivity into account.
In many physical situations, the heat-conductivity is an important driving  force to motions of fluids free surfaces, for example,  for a gaseous star.
As noted in \cite{lebovitz3}, the heat-conductivity plays an important role to  driving the evolution of a star in the phase of secular evolution, while the viscosity plays much less role.
In general, it is important and necessary to understand the role played by heat-conductivities to the dynamics of  fluids free surfaces.
However,  as far as we know, there have been no results on  the free surface  problem of  heat-conductive inviscid flows, though some results are available for viscous and heat-conductive flows, for example, in \cite{SV}, where the viscosity plays an essential role to the regularity of solutions.
If the effect of heat-conductivities is taken into account, the analysis will become  difficult due to the strong coupling  among the large variations of temperature, heat-conduction, entropy,  velocity fields and evolutions of free surfaces, while no much difficulties will be created  compared with isentropic flows if  heat-conductivities  are ignored,  because entropy is  transported along particle paths.
This is analogue to the low Mach number limit problem of compressible fluids.
As   a  step towards an understanding of the role played by  heat-conductivity to the well-posedness of free surface problems for inviscid fluids, we consider in this paper the  problem of a highly subsonic flow  which is used to approximate  general heat-conductive compressible inviscid flows when the Mach number is small  (cf. \cite{Alazard}).

We consider the following problem:  for $n=2$ and $n=3$
\begin{subequations}\label{l2}\begin{align}
& (\pl_t  + v^k \pl_k) v_j + \mathcal{T} \pl_j p =0, \ \ j=1,\cdots,n  \ \ & {\rm in} \ \  \mathscr{D}, \label{l2-a}\\
&  {\rm div }v= \kappa \Delta \mathcal{T} , \ \  (\pl_t  + v^k \pl_k) \mathcal{T} =  \kappa \mathcal{T} \Delta \mathcal{T}  , \ \ & {\rm in} \ \  \mathscr{D}, \label{l2-b}
 \end{align} \end{subequations}
describing the motion of a highly subsonic heat-conducting flow,
where the velocity field $v=(v_1, \cdots, v_n)$,  the temperature $\mathcal{T}$,  the  pressure $p$  and the domain
$\idm\subset [0,  T]\times \R^n$ are the unknowns to be determined, $\kappa>0$ is the  (scaled) heat-conductive coefficient.
Here $v^k=\da^{ki}v_i=v_k$, and we have used  the Einstein summation convention.
Given a  connected bounded domain $\mathscr{D}_0\subset \mathbb{R}^n$ and initial data $(v_0, \mathcal{T}_0)$ satisfying $ {\rm div }v_0= \kappa \Delta \mathcal{T}_0$, we want to find a set $\mathscr{D} \subset [0,T]\times \mathbb{R}^n$,  a vector field $v$ and  scalar functions $p$ and $\mathcal{T}$ solving \ef{l2} and satisfying the initial conditions:
\begin{align}\label{initial}
\mathscr{D}_0=\{x: \ (0,x) \in \mathscr{D}\} \ \ {\rm and} \ \
(v, \mathcal{T})=(v_0, \mathcal{T}_0) \ \ {\rm on} \ \ \{0\}\times \mathscr{D}_0.
\end{align}
Let $\mathscr{D}_t=\{x\in \mathbb{R}^n: \ (t,x) \in \mathscr{D} \}$. We require the boundary conditions on  the free surface $\pl\mathscr{D}_t$,
\begin{align}\label{s110}
&p=0 , \ \     \mathcal{T}= \mathcal{T}_b   \ \   {\rm and}  \ \  v_\mathcal{N}=\varpi \ \
 {\rm on} \ \  \pl \mathscr{D}_t
 \end{align}
for each $t$, where $\mathcal{T}_b $ is a  positive constant, $\mathcal{N} $  is the exterior unit normal to $\pl \mathscr{D}_t$, $v_\mathcal{N}=\mathcal{N}^iv_i$, and $\varpi$ is the normal velocity of $\pl \mathscr{D}_t$. For the derivation and physical background of system \ef{l2}, one may refer to \cite{Alazard}.

Another motivation  of this article is to serve as a step to understand the behavior of  free surface motions of inviscid heat-conducting compressible flows for small Mach numbers.  When the effect of heat-conductivity is ignored, it is well-known that the incompressible Euler equations can be derived from the compressible Euler equations when the Mach number is small.    System \ef{l2} can be derived in the same spirit,  by taking the effect of the heat conductivity into account however. 
The low Mach number limit was rigorously justified in \cite{Alazard} for the initial value  problem in entire $\mathbb{R}^n$-space or  the periodic problem of heat-conducting flows (see also \cite{DJO,JJLX} for the related results).   For fluids free surface problems, the only available result on the low Mach number limit is quite recent due to Lindblad and Luo\cite{lindluo}  for  {\it isentropic } flows  where the the low Mach number limit equations are incompressible since the effect of heat-conductivity  is ignored for an isentropic flow. However, as observed in \cite{Alazard},  the low Mach number limit equations are not incompressible anymore for  heat-conducting flows, and  limit problem becomes  more complicated and subtle. Toward to this direction of low Mach number limit problems with free surfaces for heat-conducting flows, it is important to gain a good understanding of  solutions to limiting flows since they may be used as the leading order approximation.   This is one of motivations for us to study problem   \ef{l2}-\ef{s110}.

We will prove {\em a priori} estimates for  problem \ef{l2}-\ef{s110} in Sobolev spaces when the initial data satisfies
\be\label{is11}\min_{\pl \mathscr{D}_0 }\lt(-\pl_\mathcal{ N} p\rt) >0, \ee
which implies, as we will prove, that for some $T>0$ and $0\le t\le T$,
\be\label{s11}-\pl_\mathcal{ N} p \ge \ea_b >0   \ \  {\rm on} \ \  \pl \mathscr{D}_t ,\ee
where $\pl_\mathcal{ N} =\mathcal{ N}^j \pl_j$, and $\ea_b=2^{-1}\min_{\pl \mathscr{D}_0 }\lt(-\pl_\mathcal{ N} p\rt)$.
\ef{s11} is a natural stability condition called the {\em physical condition} or  the {\em  Taylor sign condition} for an incompressible inviscid fluid in literatures (cf. \cite{bealhou,CL00,Coutand,Ebin,L1,L2,LN,SZ,wu1, wu2,zhang}),  excluding  the possibility of the Rayleigh-Taylor type instability (cf. \cite{Ebin}).
Since system \ef{l2}  keeps unchanged when a constant is added to $p$, the condition $p=0$ on  $\pl \mathscr{D}_t $ is equivalent to
that of  $p$ being a  constant  on $\pl \mathscr{D}_t $. Therefore, the boundary conditions $p=0 $ and $    \mathcal{T}= \mathcal{T}_b $ on $\pl \mathscr{D}_t $ is to match the exterior media with the constant pressure and temperature. The boundary condition $p=0$ on $\pl \mathscr{D}_t $  is commonly used for incompressible flows without surface tensions in literatures (cf. \cite{CCFG,CL00,Coutand,Ebin,GMS,GW,HaoLuo,Lindblad2,LN,SZ,wu1,wu3,zhang} and references therein).

Note  that system \ef{l2} is reduced to the usual incompressible Euler equation  if  the heat-conductivity coefficient $\kappa=0$ or   temperature $\mathcal{T}=$constant, for which a geometric approach was introduced  in \cite{CL00} to study  free surface problems  without assuming that the flow is irrotational.
We  adopt this approach to study problem \ef{l2}-\ef{s110} for a highly subsonic  heat-conductive flow.
Several new and essential analytic difficulties occur  in extending the analysis in \cite{CL00} to  heat-conductive flows, including  loss of symmetries of equations,  loss of derivatives  in closing   arguments,   the strong coupling of the large variation of temperature, heat-conduction,   compressibility of fluids due to the non-zero divergence of the velocity filed,  and the evolution of free surfaces.   These issues will be addressed  in Section \ref{sec1.3}.
We  construct a higher order energy functionals from which the Sobolev norms of $H^s(\idm_t)$ ($s=0, 1,\cdots,n+2$) of solutions can be derived.
This energy functional involves space-time derivatives of the divergence of  velocity fields,   ${\rm div} v$, for which the estimates are quite involved.
This is a key difference between  the constructions of the higher order energy functional compared with that in \cite{CL00}.   Besides the {\em a priori} estimates of Sobolev norms for problem  \ef{l2}-\ef{s110},   estimates for some geometric quantities of free surfaces, for example, the  $L^{\infty}$-bound for the second fundamental form and lower bound for the injective radius of the normal exponential map are also given.
The bounds for those geometric quantities are not only needed to bound   Sobolev norms of solutions,   but also vital  to the understanding of the evolution of the geometry of free surfaces, for example, to the study
of  the formation of singularities, such as the curvature blow-up or the self-intersection. It should be noted that the singularities  such as the splash singularity or splat singularity in \cite{CCFG,CS}   and wave crests in  \cite{wu5}   all  occur on free surfaces.

Note that, besides in \cite{CL00}, the Riemannian geometry tools (parallel transports, vector fields and covariant differentiations) were intensively used in \cite{Lindblad2,LN,lindluo,SZ,SZ1,SZ2} to study fluids free surface problems, of which one of advantages is to make full use of some intrinsic properties of the studied problems independent of choice of coordinates.  The geometric approach used in \cite{CL00} was also adopted to study free surface problems of incompressible MHD flows  in \cite{HaoLuo} and incompressible Neo-Hookean elastodynamics  in \cite{HaoWang}.

\subsection{Main results}\label{sec1.2}
We will prove that the temporal derivative of the constructed higher order energy functional  is controlled by itself.  This higher order functional consists of a boundary part and an interior part. In order to define the boundary integral,    we   project the equations   to the tangent space of the boundary as in \cite{CL00}.   The orthogonal projection $\Pi$ to the tangent space of the boundary of a $(0, r)$ tensor $\alpha$ is defined to be the projection of each component along the normal as follows:

\begin{defn}\label{defn1.1} The orthogonal projection $\Pi$ to the tangent space of the boundary of a $(0,r)$ tensor $\al$ is defined to be the projection of each component along the normal:
$$
(\Pi \al)_{i_1\cdots i_r} =  \Pi_{i_1}^{j_1}\cdots \Pi_{i_r}^{j_r} \al_{j_1\cdots j_r}, \ \ {\rm where} \ \  \Pi_{i}^{j}=\da_i^j-\mathcal{N}_i\mathcal{N}^j .
$$
The tangential derivative of the boundary is defined by $\bar{\pl}_i=\Pi_i^j\pl_j$, and the second fundamental form of the boundary is defined by $\ta_{ij}=\bar{\pl}_i \mathcal{ N} _j$.
\end{defn}

As  in \cite{CL00}, we also need  a positive definite quadratic form $Q(\alpha, \beta)$ for tensors $\alpha$ and $\beta$ of the same order which is the inner product of the tangential components  when restricted to the boundary, and $Q(\al,\al)$ increases to the norm $|\al|^2$ in the interior. For this purpose, we  extend the normals on the boundary to the interior as follows:

\begin{defn}\label{defn1.2}
Let $\iota_0$ be the injectivity radius of the normal exponential map of $\pl \mathscr{D}_t$, i.e., the largest number such that the map
$$
\pl \mathscr{D}_t \times (-\iota_0, \iota_0)   \  \to \ \{x\in \mathbb{R}^n: \  {\rm dist} (x,\pl \mathscr{D}_t)< \iota_0 \} :  \ \
(\bar x, t)   \   \mapsto \ x=\bar x + \iota \mathcal{N}(\bar x)
$$
is an injection.
\end{defn}

As in \cite{CL00},  we also need the definition of  $\iota_1$ as follows:
\begin{defn}\label{defn2.3}
Let $0<\ea_1\le 1/2$ be a fixed number, and let $\iota_1=\iota_1(\ea_1)$ be the largest number such that
$$
|\mathcal{N}(\bar x_1)-\mathcal{N}(\bar x_2)| \le \ea_1 \ \ {\rm whenever} \ \
|\bar x_1 - \bar x_2|\le \iota_1, \ \     \bar x_1, \bar x_2 \in \pl\mathscr{D}_t.
$$
\end{defn}
As shown in \cite{CL00}, $\iota_1$ is equivalent to $\iota_0$ in conjunction with a bound of the second fundamental form $\ta$. We do not need it for the statement of the main theorem, but we will need it when we illustrate the main idea of the proof of our theorem, so we give its definition here, together with $\iota_0$.

\begin{defn}\label{defn1.3}
Let $d_0$ be a fixed number such that $\iota_0/16 \le d_0 \le \iota_0/2$, and $\eta$ be a smooth cutoff function on $[0, \iy)$ satisfying  $0\le \eta(s)\le 1$, $\eta(s)=1$ when $s\le  d_0/4$, $\eta(s)=0$ when $ s\ge d_0/2$, and  $|\eta'(s)|\le 8/d_0$.
Define
$$
\varrho^{ij}(t,x)=\da^{ij}-\eta^2(d(t,x)) \mathcal{N}^i(t,x) \mathcal{N}^j(t,x) \ \   {\rm in} \ \ \mathscr{D}_t,
$$
where
$$  \mathcal{N}^j(t,x)=\da^{ij} \mathcal{N}_i(t,x),         \ \ \mathcal{N}_i (t,x) =\pl_i d(t,x)  \ \  {\rm and}
\ \  d(t, x)={\rm dist}(x,\pl \mathscr{D}_t).   $$
In particular, $\varrho$ gives the the induced metric on the tangential space to the boundary:
$$ \varrho^{ij}=\da^{ij}-\mathcal{N}^i\mathcal{N}^j , \ \ \varrho_{ij}=\da_{ij}-\mathcal{N}_i\mathcal{N}_j \ \  {\rm on} \ \ \pl \mathscr{D}_t.$$
\end{defn}
With this setting, the above mentioned quadratic form $Q(\alpha, \beta)$ for $(0, r)$ tensors is defined by
\be\label{quadratic}
Q(\al, \ba)=  \varrho^{i_1j_1}\cdots  \varrho^{i_rj_r} \al_{i_1\cdots i_r}  \ba_{j_1\cdots j_r} .\ee

We are concerned with the problem for the fixed $\kappa$ in this paper,  so we set $\kappa=1$ from now on for the simplicity of the presentation. The  energy functionals of each order are then defined by
\begin{subequations}\label{norm'}\begin{align}
  E_0 (t) = &  \int_{\mathscr{D}_t} \mathcal{T}^{-1}  |v|^2 dx ,  \label{norm'-a}  \\
  E_r (t) = &   \int_{\mathscr{D}_t}  \mathcal{T}^{-1}  \da^{mn}   Q(\pl^r v_m , \pl^r v_n)  dx   +  \int_{\mathscr{D}_t} |\pl^{r-1} {\rm curl} v |^2    dx
 +   \int_{\mathscr{D}_t}  |\pl D_t^{r-1}{\rm div} v|^2   dx \notag\\
& +   \int_{\pl\mathscr{D}_t }   Q(\pl^r  p,  \pl^r  p) (-\pl_\mathcal{ N} p)^{-1} d S , \ \ \ \ r\ge 1 , \label{norm'-b}
\end{align}\end{subequations}
where $D_t=\pl_t + v^k \pl_k$. The higher order energy functional is defined by $\sum_{r=0}^{n+2}E_r(t)$.

In order to state the main result  of the paper, we set
\begin{subequations}\label{17initial}\begin{align}
&{\rm Vol} \mathscr{D}_0=\int_{\mathscr{D}_0}dx ,  \ \
 K_0 = \max_{x\in\pl \mathscr{D}_0 } \{ |\ta(0, x)| + | \iota_0^{-1}(0, x)| \} ,   \ \  \underline{\ea}_0=\min_{x\in\pl \mathscr{D}_0 }( -\pl_\mathcal{N}  p )(0,x),  \label{s118}\\
& \underline{\mathcal{T}}_0= \min_{x\in\mathscr{D}_0 }\mathcal{T}(0,x), \ \  \overline{\mathcal{T}}_0 =\max_{x\in\mathscr{D}_0 }\mathcal{T}(0,x)  ,
\label{s112}\\
&   M_0= \max_{x\in\mathscr{D}_0 } \{ |\pl p(0,x)| + |\pl v(0,x)| + |\pl  \mathcal{T}(0,x) |    \} .\label{s113}
\end{align}\end{subequations}
The initial pressure $p_0(x)=p(0, x)$ is determined by the following Dirichelet problem:
$${\rm div}(\mathcal{T}_0\pl p_0)=-\pl_i v_0^j \pl_j v_0^i-D_t {\rm div} v\big|_{t=0} \ \  {\rm in}\ \ \mathscr{D}_0,\ \  p_0=0\ \ {\rm on} \  \ \pl\mathscr{D}_0,$$
where $D_t {\rm div} v |_{t=0}$ can be given in terms of initial values $v_0$ and $\mathcal{T}_0$ via the equations
${\rm div }v=  \Delta \mathcal{T}$ and  $ D_t \mathcal{T} =   \mathcal{T} \Delta \mathcal{T}$.
With these notations, the main theorem of the present work is stated as follows:
\begin{thm}\label{mainthm} Let $n=2,3$.    Suppose that
 $$ 0<{\rm Vol} \mathscr{D}_0, \ \underline{\ea}_0,\  \underline{\mathcal{T}}_0, \ \overline{\mathcal{T}}_0 <\infty, \ K_0, \ M_0<\infty. $$
Then there are continuous functions
$\mathscr{T}_n$   such that if
$$
T\le \mathscr{T}_n\lt({\rm Vol}\mathscr{D}_0 , K_0,  \underline{\ea}_0^{-1},  \underline{\mathcal{T}}_0^{-1}, \overline{\mathcal{T}}_0 ,M_0,  E_0(0),\cdots,E_{n+2}(0) \rt),
$$
then any smooth solution of the free surface problem \ef{l2}-\ef{s110} for $0\le t\le T$ satisfies
\begin{subequations}\begin{align}
& \sum_{s=0}^{n+2}E_s(t)\le   2 \sum_{s=0}^{n+2} E_s(0), \ \  \ \  0\le t\le T,\\
& 2^{-1} {\rm Vol}\mathscr{D}_0     \le  {\rm Vol} \mathscr{D}_t     \le  2 {\rm Vol}\mathscr{D}_0, \ \ \ \  0\le t\le T,\\
& \underline{\mathcal{T}}_0 \le \mathcal{T} \le \overline{\mathcal{T}}_0  \ \ {\rm in} \ \  \mathscr{D}_t,   \ \ \ \    0\le t\le T,  \\
& |\ta|+ |  \iota_0^{-1} |\le C K_0  \ \ {\rm on} \ \  \pl\mathscr{D}_t,  \ \  \ \   0\le t\le T,
\\
& -\pl_\mathcal{N}  p  \ge 2^{-1}\underline{\ea}_0    \ \ {\rm on} \ \  \pl\mathscr{D}_t,   \ \ \ \    0\le t\le T,
 \end{align}\end{subequations}
for a certain constant $C$, where ${\rm Vol} \mathscr{D}_t=\int_{\mathscr{D}_t}dx$.
\end{thm}

\begin{rmk} The  bound for $\|\pl (v, p)(0)  \|_{L^\iy(   \mathscr{D}_0)}$ was not needed  in \cite{CL00} to show their result for the problem of  incompressible Euler equations with a free surface, because it could be controlled by initial values of their higher order energy functional,  ${\rm Vol} \mathscr{D}_0$  and  $\iota_1(0, x)$ via
Sobolev lemmas and elliptic estimates. We need this bound, which cannot be controled as  in \cite{CL00},  to prove Theorem \ref{mainthm}, which reflects the subtlety of our problem due to the complicated coupling of variation of
temperature, heat-conduction, and compressibility of the fluid in our analysis.
\end{rmk}

We  give some remarks  on the  choice of the higher order energy functional, and  explain briefly the reason why we need $n+2$ derivatives in this functional,  while only $n+1$ derivatives were needed in \cite{CL00} when $n=2,3$.
Let
\be\label{era}\begin{split}
   E_r^a (t) = &   \int_{\mathscr{D}_t}  \mathcal{T}^{-1}  \da^{mn}   Q(\pl^r v_m , \pl^r v_n)  dx   +  \int_{\mathscr{D}_t} |\pl^{r-1} {\rm curl} v |^2    dx    \\
   &+   \int_{\pl\mathscr{D}_t }   Q(\pl^r  p,  \pl^r  p) (-\pl_\mathcal{ N} p)^{-1} d S ,  \  \ \ \  r \ge 1.
\end{split}\ee
Note that $ E_r^a$ ($r\ge 1$)  correspond to the energy functionals employed in \cite{CL00} for the study of an incompressible flow when $\mathcal{T}$ is constant.
In order to control the $L^2$-norm of $\pl^r v$ for compressible flows, one may attempt to use the following:
\be\label{comp-norm}\begin{split}
  \widetilde{E}_r (t) =  {E}_r^a (t) + \int_{\mathscr{D}_t} |\pl^{r-1} {\rm div} v |^2    dx  , \ \ \ \ r\ge 1 .
\end{split}\ee
However, \ef{comp-norm} does not work for the study the problem \ef{l2}-\ef{s110}. In fact, ${\rm div } v$ satisfies the following parabolic type equation:
 \be\label{equfordiv}
 D_t  {\rm div } v  - \mathcal{T} \Delta {\rm div } v = \textrm{other terms},
\ee
which requires a control of one more spatial derivative of ${\rm div} v$, besides
$\pl^r v$.  Here and thereafter, ``other terms" means something that does not affect the terms we single out to discuss.
So, one may try to include the $L^2$-norm of $\pl^{r} {\rm div} v $ into the $r$-th order energy functional:
\be\label{comp-norm-new}\begin{split}
  \overline{E}_r (t) =  {E}_r^a (t) + \int_{\mathscr{D}_t} |\pl^{r} {\rm div} v |^2    dx  , \ \ \ \ r\ge 1 .
\end{split}\ee
Due to \ef{equfordiv} and the boundary condition that ${\rm div } v=0$ on $\pl\mathscr{D}_t$, it is more convenient to  include the $L^2$-norm of $\pl D_t^{r-1}{\rm div} v$, instead of $\pl^{r} {\rm div} v $, into the $r$-th order energy functional \ef{norm'-b}. This is one of reasons why we choose such an energy functional.  Indeed, one can see from the proof that it is not sufficient to study the problem \ef{l2}-\ef{s110} even adopting \ef{comp-norm-new}.

The choice of the higher order energy functional $\sum_{r=0}^{n+2}E_r$ enables us to  prove that the temporal derivative of it  can be controlled by itself under the following  {\it a priori} assumptions:
\begin{subequations}\label{apriori17}\begin{align}
& \underline{V}  \le  {\rm Vol}\mathscr{D}_t (t)     \le   \overline{V}  \  \  & {\rm on} \ \   &[0,T] ,  \label{apriori17a} \\
&|\ta|+ 1/\iota_0 \le K, \ \  \ \  -\pl_\mathcal{N} p \ge \ea_b >0  \ \  & {\rm on} \ \  & \pl \mathscr{D}_t , \label{apriori17b} \\
&  \sum_{i=1}^{n-1}   (|\pl_\mathcal{N}  D_t^i p| +|\pl_\mathcal{N}  D_t^i {\rm div} v|) +|\pl^2 p|  \le L   \ \ &   {\rm on} \ \ & \pl \mathscr{D}_t , \label{apriori17c} \\
&   |\pl  p | + |\pl v|  + |\pl  \mathcal{T} | + |\pl  {\rm div} v|   \le M      \ \  &  {\rm in} \ \ &     {\mathscr{D}_t}, \label{apriori17d} \\
&  |D_t  p| +  | D_t  {\rm div} v |  + |\pl^2 \mathcal{T}|  \le \widetilde{M}    \ \  &  {\rm in} \ \ &    \mathscr{D}_t . \label{apriori17e}
\end{align}\end{subequations}
It should be noted that the {\em a priori} assumptions adopted in \cite{CL00} for incompressible flows are the following:
\begin{align*}
&|\ta|+ 1/\iota_0 \le K, \ \  \ \  -\pl_\mathcal{N} p \ge \ea_b >0  \ \  & {\rm on} \ \  & \pl \mathscr{D}_t ,  \\
&     |\pl_\mathcal{N}  D_t p|  +|\pl^2 p|  \le L   \ \ &   {\rm on} \ \ & \pl \mathscr{D}_t ,  \\
&   |\pl  p | + |\pl v|     \le M      \ \  &  {\rm in} \ \ &     {\mathscr{D}_t}. \end{align*}
In closing the argument, the  {\em a priori} assumptions, for example, on the $L^{\infty}$-bounds for $\pl (v,p)$ in $\mathscr{D}_t$ and $\ta$ on $\pl\mathscr{D}_t$, need to be verified  both in \cite{CL00} and this article.
In fact,  these $L^{\infty}$-bounds could be controlled in \cite{CL00} by their higher order energy functional, ${\rm Vol} \mathscr{D}_0$, $\min_{\pl \mathscr{D}_t}(-\pl_\mathcal{N} p)$ and $\max_{\pl \mathscr{D}_t}(\iota_1^{-1})$ via Sobolev lemmas, elliptic estimates and projection formulae.  But this is not the case for the problem studied in this paper,  that is, we do not have such simple and neat control of  $\pl(v,p)$   and $\ta$.  This is a key feature for the problem studied here.
Instead of using the method adopted in \cite{CL00}, we employ the evolution equations for $\pl(v,p)$ and $\ta$, which  causes the loss of derivatives.  For example, in order to control  $\|\theta\|_{L^{\infty}(\pl\mathscr{D}_t)}$,  we will need the control of  $\| \pl^2 v\|_{L^{\infty}(\pl\mathscr{D}_t)}$, while a projection formula was used in \cite{CL00} to control $\|\theta\|_{L^{\infty}(\pl\mathscr{D}_t)}$ for which there is no need to control   $\| \pl^2 v\|_{L^{\infty}(\pl\mathscr{D}_t)}$. This loss of derivative in the control of the $L^{\infty}$-bound for $\theta$ forces us to use $n+2$ derivatives in the higher order functional, while only $n+1$ derivatives were needed  in \cite{CL00} for $n=2, 3$.  We will address these issues in more details in the next subsection.

\subsection{Main issues and novelty in analysis}\label{sec1.3}
We first highlight  the main issues in extending the analysis in \cite{CL00} to  problem  \ef{l2}-\ef{s110} and then present the main strategy of the proof of Theorem \ref{mainthm}. The big obstacle in the analysis lies in the strong coupling of large variation of temperature, heat-conduction and compressibility of fluids due to the non-zero divergence of the velocity filed, which creates  essential and new challenges in the analysis. It should be noted that the analysis in this work for $\mathcal{T}=$constant or ${\rm div  } v=0$ reduces to that in \cite{CL00}.   Indeed,  the sharp estimates in \cite{CL00,lindluo} use  all the
symmetries of the incompressible or isentropic Euler equations, which are missing for   \ef{l2} we  consider here. The loss of symmetries of the equations  we study is reflected by the following facts:  for the problems of incompressible  or isentropic Euler equations studied in \cite{CL00,lindluo}, the zero-th order energy functional  is conserved in time, and the temporal derivative of the $r$-th $(r\ge 1)$ order energy functional  can be controlled  by  lower order functionals  under
some suitable   a priori  assumptions. However, in our case, the temporal derivatives of the zero-th  and the first order energy functionals $E_0$  and   $E_1$ depend on the higher order ones. The fact that $E_0$  is not conserved indicates some kind of loss of symmetries of the equations studied in this paper.

Another difficulty in our analysis is to deal with the problem  of loss of derivatives when we work on   evolution equations for some quantities
in the {\em a priori} assumptions to obtain the bounds for them to close the argument.  The first one is on the second fundamental form $\ta$ for  free surfaces.    The projection formula,
\be\label{projection}\Pi (\pl^2 p)=\theta\pl_{\mathcal{N}} p \ \ {\rm on} \ \  \pl \mathscr{D}_t,\ee  was used  to estimate the $L^{\infty}$-bound for $\theta$ in \cite{CL00}.  The reason that this can work in \cite{CL00} is because one may obtain the $L^{\infty}$-bound for $\pl^2 p$ on $\pl\mathscr{D}_t$ independent of the $L^{\infty}$-bound for $\theta$, which, together with  the lower bound for $-\pl_{\mathcal{N}} p$ due to the Taylor sign condition, gives the $L^{\infty}$-bound for $\ta$.
Indeed,  it was proved in \cite{CL00} that
\be\label{9/11-1}
\|\pl^2 p\|_{L^{\infty}(\pl \mathscr{D}_t)}\le C(K_1) \sum_{r=2}^{n+1}  \|\pl^r p\|_{L^2( \pl \mathscr{D}_t)} \le C(K_1, \mathcal{E}_0, \cdots, \mathcal{E}_{n+1}, {\rm Vol} \mathscr{D}_t),  \ \  n=2,3,
\ee
where $\mathcal{E}_0=E_0$ and $\mathcal{E}_r=E_r^a $ $(r\ge 1)$ with $\mathcal{T}=1$,  $K_1$ is the upper bound for $1/\iota_1$ on $\pl\mathscr{D}_t$ with $\iota_1$ given in Definition \ref{defn2.3}.  In the same spirit,   the $L^{\infty}$-bound for $\theta$ was obtained  in \cite{lindluo} for isentropic Euler equations by replacing the pressure $p$ in \ef{projection} by the enthalpy $h$.   However,  we can only obtain,  for problem  \ef{l2}-\ef{s110} that
 \be\label{9/11-2}\begin{split}
\|\pl^2 p\|_{L^{\infty}( \pl \mathscr{D}_t)} \le   &  C(K_1, {\rm Vol} \mathscr{D}_t)  \| \mathcal{T}^{-1}\|_{L^\iy(   \mathscr{D}_t)} \|\pl p\|_{L^\iy(   \mathscr{D}_t)} \|\ta \|_{L^\iy( \pl \mathscr{D}_t)} \| \pl^{n} \mathcal{T}  \|_{L^2( \pl  \mathscr{D}_t)} \\
  & +  \textrm{other terms}, \ \  n=2,3,
 \end{split}\ee
from which it is clear that the projection formula used in \cite{CL00} to give the $L^{\infty}$-bound for $\ta$ cannot work directly for our problem. Indeed, \ef{9/11-2} follows from Sobolev lemmas and the following estimates:
\be\label{9/11-3}\begin{split}
   \|\pl^{n+1} p\|_{L^2( \pl \mathscr{D}_t)}
\le C  \|\Pi \pl^{n+1} p\|_{L^2( \pl \mathscr{D}_t)}  +   C(K_1, {\rm Vol} \mathscr{D}_t) \sum_{r=0}^{n} \|\pl^{r} \Delta p\|_{L^2(   \mathscr{D}_t)},
 \end{split}\ee
\be\label{9/11-4}\begin{split}
   \|\pl^{n} \Delta p\|_{L^2(   \mathscr{D}_t)}
   \le \| \mathcal{T}^{-1}\|_{L^\iy(   \mathscr{D}_t)} \|\pl p\|_{L^\iy(   \mathscr{D}_t)} \| \pl^{n+1} \mathcal{T}  \|_{L^2(   \mathscr{D}_t)} + \textrm{other terms},
 \end{split}\ee
\be\label{9/11-5}\begin{split}
& \| \pl^{n+1} \mathcal{T}  \|_{L^2(   \mathscr{D}_t)}   \le  C \|\Pi \pl^{n+1} \mathcal{T} \|_{L^2( \pl \mathscr{D}_t)}   +  \textrm{other terms}  \\
& \le C   \|\pl_\mathcal{N} \mathcal{T} \|_{L^\iy( \pl \mathscr{D}_t)}  \| \bar\pl^{n-1} \ta  \|_{L^2( \pl  \mathscr{D}_t)} + C \|\ta \|_{L^\iy( \pl \mathscr{D}_t)} \| \pl^{n} \mathcal{T}  \|_{L^2( \pl  \mathscr{D}_t)}      +  \textrm{other terms}.
 \end{split}\ee
Here \ef{9/11-3}, \ef{9/11-4} and \ef{9/11-5} follow from elliptic estimates, the equation $\mathcal{T}\Delta p = -(\pl \mathcal{T})\cdot \pl p$+other terms, and the
 projection formula,  respectively.

Instead of using the projection formula, we need to use the evolution equations for $\ta$.  By doing so, we are led to the following estimate:
$$|D_t \ta| \le |\pl^2 v| + C |\ta||\pl v|,$$
from which it is clear that we need to  get the  $L^{\infty}$-bounds for both  $\pl v$ and $\pl^2 v$ on $ \pl  \mathscr{D}_t$, while only the $L^{\infty}$-bound for  $ \pl v$ was sufficient   in \cite{CL00}. Thus,  the $L^{\infty}$-bound for one more derivative of the velocity field than  that in \cite{CL00} is needed. This causes the loss of one more derivative than that in \cite{CL00}.
Hence,  we need to estimate $n+2$ derivatives in the energy functionals to close the argument, while  only $n+1$ derivatives were
needed in \cite{CL00} for $n=2,3$.
It should be noted that  only $\pl  v$ enters equations \ef{l2}, but not $\pl^2 v$,  and thus one may think that the  estimate of $\pl v$ may be sufficient to close the argument as done in \cite{CL00}.
But the above argument suggests that this is not the case for the problem \ef{l2}-\ef{s110} which reflects the subtlety of this problem .
It is extremely involved to bound $\pl^2 v$  before one obtains the $L^{\infty}$-bound for $\ta$ in our case,   due to the strong coupling of variation of temperature, heat-conduction,  compressibility of the fluid and the evolution of the free surface.

In fact, even for the  $L^\iy$-bound for $ \pl v$ in  $\mathscr{D}_t$, we will have to use the evolution equation of $ \pl v$, while  it was obtained by the Sobolev lemma in \cite{CL00}:
\be\label{plv}
\|\pl v \|^2_{L^{\infty}( \mathscr{D}_t)}\le C(K_1)\sum_{r=1}^{3} \|\pl^r v\|^2_{L^{2}( \mathscr{D}_t)} \le   C(K_1 ) \sum_{r=1}^{3} \mathcal{E}_r \le C(K_1 ) \sum_{r=1}^{n+1} \mathcal{E}_r ,  \ \  n=2,3.
\ee
For the problem considered in this paper, we do not have such a simple and neat estimate due to the complicated coupling as mentioned above. Indeed, if we
try to use the Sobolev lemma as in \cite{CL00},  we can only get a bound depending on the $L^{\infty}$-bound for $\ta$ that cannot be controlled by $n+1$ derivatives,  as shown in the  following:
\be\label{911}\begin{split}
\|\pl v \|_{L^{\infty}( \mathscr{D}_t)}^2 \le &C(K_1, {\rm Vol} \mathscr{D}_t)  \| \mathcal{T}^{-1}\|_{L^\iy(   \mathscr{D}_t)}^2 \| \pl \mathcal{T} \|_{L^\iy(   \mathscr{D}_t)}  \| v \|_{L^\iy(   \mathscr{D}_t)}  \\
&  \times \|\ta \|_{L^\iy( \pl \mathscr{D}_t)} \| \pl^{2} \mathcal{T}  \|_{L^2( \pl  \mathscr{D}_t)}   \| \pl^2  v   \|_{L^2(   \mathscr{D}_t)}    +  \textrm{other terms},
 \end{split}\ee
which follows from Sobolev lemmas and the following estimates:
\be\label{911-1}\begin{split}
\|\pl^3 v \|_{L^{\infty}(   \mathscr{D}_t)}^2 \le    C (E^a_3 +  \|\pl^2 {\rm div} v \|^2_{L^2(   \mathscr{D}_t)} )
\le  C  E^a_3   +  C(K_1, {\rm Vol} \mathscr{D}_t)\|\Delta {\rm div} v \|^2_{L^2(   \mathscr{D}_t)},
 \end{split}\ee
\be\label{911-2}\begin{split}
\|\Delta {\rm div} v \|_{L^2(   \mathscr{D}_t)}^2
\le \| \mathcal{T}^{-1}\|_{L^\iy(   \mathscr{D}_t)}^2 \| \pl \mathcal{T} \|_{L^\iy(   \mathscr{D}_t)}  \| v \|_{L^\iy(   \mathscr{D}_t)}  \sum_{k=1}^3   \| \pl^k \mathcal{T}  \|_{L^2(   \mathscr{D}_t)} \sum_{j=0}^2  \| \pl^j  v   \|_{L^2(   \mathscr{D}_t)} ,
 \end{split}\ee
\be\label{911-3}\begin{split}
 \| \pl^3 \mathcal{T}  \|_{L^2(   \mathscr{D}_t)}  \le  C \|\ta \|_{L^\iy( \pl \mathscr{D}_t)} \| \pl^{2} \mathcal{T}  \|_{L^2( \pl  \mathscr{D}_t)}      +  \textrm{other terms}.
 \end{split}\ee
Here \ef{911-1}, \ef{911-2}, and \ef{911-3} follow from the divergence-curl decomposition,
the equation $\mathcal{T}\Delta  {\rm div} v = \pl^2 \mathcal{T} \cdot \pl v $+other terms,
and \ef{9/11-5}, respectively.
The evolution equation  $D_t \pl v =- \mathcal{T} \pl^2 p $+other terms  and the Sobolev lemma lead to
\be\label{dtplv}\begin{split}
\| D_t \pl v  \|_{L^\iy(   \mathscr{D}_t)}  \le &  \| \mathcal{T} \|_{L^\iy(   \mathscr{D}_t)} \| \pl^2 p \|_{L^\iy(   \mathscr{D}_t)} +\textrm{other terms}\\
\le &  C(K_1)\| \mathcal{T} \|_{L^\iy(   \mathscr{D}_t)} \sum_{r=2}^4\| \pl^r p \|_{L^2(   \mathscr{D}_t)}  +\textrm{other terms}, \ \ n=2,3,
 \end{split}\ee
from  which it is clear again  that $n+2$ derivatives are needed to obtain  the  $L^\iy$-bound for $  \pl v $ in the case of $n=2$.

\subsection{The strategy of the proof}

Next, we present  the strategy of the proof. We want to prove that the temporal derivative of the higher order energy functional $\sum_{r=0}^{n+2} E_r$ can be bounded by itself under
the {\em a priori } assumptions \ef{apriori17}.
For $E^a_r$ given in \ef{era}, we can prove that
\begin{subequations}\label{dtea}\begin{align}
  \frac{d}{dt}E_1^a(t)\le   &  C_1( \cdot )
       \|\pl (v, \ p, \ {\rm div} v )\|_{L^2(\mathscr{D}_t)}^2  ,
 \\
 \frac{d}{dt} E_r^a(t) \le & C_r(\cdot) \lt( \lt\|\pl^r \lt(v, \ p, \ \mathcal{T}, \ {\rm div} v \rt)\rt\|_{L^2(\mathscr{D}_t)}^2 +  \lt\| \lt( \pl^{r-1} v, \  \pl^r p, \  \Pi \pl^r D_t p  \rt)\rt\|_{L^2(\pl\mathscr{D}_t)}^2 \rt) \notag\\
 & + \textrm{other terms},  \ \  r=2,\cdots,n+2,
\end{align}\end{subequations}
where and thereafter $C_r(\cdot)$ stands for a constant depending continuously on the bounds in the a priori assumptions \ef{apriori17}. We need to control the quantities  on the right-hand side of \ef{dtea} by the energy functionals. We will mainly discuss the estimates for the pressure $p$ which  appear also in \cite{CL00} but require  additional works in our problem due to the involvement  of ${\rm div} v $  in  these estimates.

It follows from the definition of the energy functional  $E_r^a$ that for $r\ge 2$,
$$
\|\Pi\pl^r p\|_{L^2(\pl\mathscr{D}_t)}^2  \le  \| \pl_\mathcal{N} p \|_{L^\iy(\pl\mathscr{D}_t)} E_r^a.
$$
Since $-\mathcal{T} \Delta p= D_t {\rm div} v +$other terms, one can use {\em elliptic estimates} to control all components of $\pl^r p$ from the tangential components $\Pi\pl^r p$ in the energy:
\be\label{9/7-1}\begin{split}
    \|\pl^r p\|_{L^2(\pl\mathscr{D}_t)}^2  &\le  C\|\Pi\pl^r p\|_{L^2(\pl\mathscr{D}_t)}^2 + C(K,{\rm Vol}\mathscr{D}_t)  \sum_{s=0}^{   r-1}\|\pl^s \Delta p\|^2_{L^2(\mathscr{D}_t)}  \\
 &\le C_r(\cdot)\lt(E_r^a +\|\pl^{r-1}D_t {\rm div} v\|_{L^2(\mathscr{D}_t)} ^2 \rt)+ \textrm{other terms}.
 \end{split}\ee
Under the physical condition $-\pl_\mathcal{N} p \ge \ea_b >0$, we can use the higher order version of the projection formula to get
\be\label{9/7-3}\begin{split}
    \|\bar \pl^{r-2}\theta\|_{L^2(\pl\mathscr{D}_t)}^2
    &\le   C_r(\cdot)\lt( \|\Pi \pl^r p \|_{L^2(\pl\mathscr{D}_t)}^2+ \|\pl^{r-1}p \|_{L^2(\pl\mathscr{D}_t)}^2\rt)+  \textrm{other terms}\\
    &\le C_r(\cdot)\lt(E_r^a+E_{r-1}^a + \|\pl^{r-2}D_t {\rm div} v\|_{L^2(\mathscr{D}_t)}^2\rt)+ \textrm{other terms}.
 \end{split}\ee
Once we have the bound for the second fundamental form, we can get estimates for solutions of any Dirichlet problem of elliptic equations. So, we can get estimates for $\mathcal{T}$,  $ {\rm div} v $  and $D_t p$, which satisfy  elliptic equations:
$\Delta \mathcal{T} = {\rm div} v$, $ \mathcal{T} \Delta {\rm div } v =  D_t  {\rm div } v +$other terms, and $ \mathcal{T} \Delta D_t p= - D_t^2 {\rm div} v +$other terms. Since the equation for  $ D_t p$ involves the highest order temporal derivative of ${\rm div} v $, we  show  here how to control $D_tp$.
\be\label{9/7-2}\begin{split}
  &  \|\Pi\pl^{r} D_t p\|_{L^2(\pl\mathscr{D}_t)}^2 \\
    \le & C\| \pl_\mathcal{N}  D_t p\|_{L^\iy(\pl\mathscr{D}_t)}^2\| \bar \pl^{r-2}\theta \|_{L^2(\pl\mathscr{D}_t)}^2 + C_r(\cdot)  \|\pl^{r-1}  D_t p  \|_{L^2( \pl \mathscr{D}_t)}^2+  \textrm{other terms}\\
   \le &  C_r(\cdot)\lt( \| \bar \pl^{r-2}\theta \|_{L^2(\pl\mathscr{D}_t)}^2+ \|\pl^{r-2} \Delta D_t p  \|_{L^2( \mathscr{D}_t)}^2\rt)+  \textrm{other terms}\\
     \le &C_r(\cdot)\lt(E_r^a+E_{r-1}^a + \|\pl^{r-2}D_t {\rm div} v\|_{L^2(\mathscr{D}_t)}^2 +\|\pl^{r-2}D_t^2 {\rm div} v\|_{L^2(\mathscr{D}_t)}^2  \rt)+ \textrm{other terms}.
 \end{split}\ee

Since the terms involving ${\rm div} v$ on the right-hand side of \ef{9/7-1} and \ef{9/7-2} cannot be controlled by $\sum_{s\le r}E^a_s$, we introduce
$$
 E_r^d(t)=    \int_{\mathscr{D}_t}  |\pl D_t^{r-1}{\rm div} v|^2   dx, \ \  r=1,\cdots, n+2,
$$
so that
$$
E_r(t)= E_r^a (t) +E_r^d(t),  \ \ r=1,\cdots, n+2.
$$
It can be proven that the terms on the right-hand side of \ef{dtea} can be controlled by $C_r(\cdot)\sum_{s\le r} E_s$ when $r\ge 2$. For example, the idea of the estimates for $\pl^r {\rm div} v$ can be illustrated as follows:
\bee\label{}\begin{split}
     & \|\pl^r {\rm div} v\|_{L^2( \mathscr{D}_t)}^2  \\
     \le & C \|\Pi \pl^r {\rm div} v\|_{L^2(\pl\mathscr{D}_t)}^2+  C(K,{\rm Vol}\mathscr{D}_t)\sum_{s=0}^{r-2}\|\pl^{s}\Delta {\rm div} v\|_{L^2(\mathscr{D}_t)}^2    \\
    \le &  C \| \pl_\mathcal{N} {\rm div} v\|_{L^2(\pl\mathscr{D}_t)}^2 \|\bar \pl^{r-2} \theta\|_{L^2(\pl\mathscr{D}_t)}^2 +C_r(\cdot) \lt(\|\pl^{r-1} {\rm div} v\|_{L^2(\pl\mathscr{D}_t)}^2+\|\pl^{r-2}\Delta {\rm div} v\|_{L^2(\mathscr{D}_t)}^2\rt) \\
    & +  \textrm{other terms} \\
    \le & C_r(\cdot) \lt( \|\bar \pl^{r-2} \theta\|_{L^2(\pl\mathscr{D}_t)}^2 +\|\pl^{r-2}\Delta {\rm div} v\|_{L^2(\mathscr{D}_t)}^2\rt) +  \textrm{other terms} \\
      \le & C_r(\cdot) \lt(E_r^a  + \|\pl^{r-2}D_t {\rm div} v\|_{L^2(\mathscr{D}_t)}^2  +\|\pl^{r-2}\Delta {\rm div} v\|_{L^2(\mathscr{D}_t)}^2\rt) +  \textrm{other terms},
 \end{split}\eee
 where \ef{9/7-3} has been used to derive the last inequality. Using the  equation $ \mathcal{T} \Delta {\rm div } v =  D_t  {\rm div } v +$other terms,  we may obtain
\bee\label{}\begin{split}
       \|\pl^r {\rm div} v\|_{L^2( \mathscr{D}_t)}^2
      \le & C_r(\cdot) \lt(E_r^a  + \|\pl^{r-2}D_t {\rm div} v\|_{L^2(\mathscr{D}_t)}^2  \rt) +  \textrm{other terms} \\
      \le & C_r(\cdot) \lt(E_r^a  +  E_{r}^d  +  E_{r-1}^d   \rt) +  \textrm{other terms}.
 \end{split}\eee
Here the equation   $ \mathcal{T} \Delta D_t {\rm div } v =  D_t^2  {\rm div } v +$other terms has been used to obtain the last inequality for $r=4,5$.

We need to estimate the temporal derivative of $E_r^d$. One may get
\bee\label{}\begin{split}
  \frac{d}{dt}  E_r^d(t)
 \le &C_r(\cdot )  \lt(    \sum_{j=1}^{r-2}   \lt\| \pl^2 D_t^{j}(p, \ {\rm div} v) \rt\|_{L^2(\mathscr{D}_t)}^2   +  \sum_{j=1}^{r-3}   \lt\| \pl^3 D_t^{j} p  \rt\|_{L^2(\mathscr{D}_t)}^2  \rt. \\
& \lt. + \lt\| \pl^r(v, \ p , \ {\rm div} v )    \rt\|_{L^2(\mathscr{D}_t)}^2 +    E_r^d \rt)
 + \textrm{other terms} ,  \ \  r=1,\cdots, n+2.
\end{split}\eee
The task is then to control $\pl^2 D_t^{r-2}p$ $(r\ge 3)$ and  $\pl^3 D_t^{r-3} p$ $(r\ge 4)$ by    energy functionals.
It follows from the equation $ -\mathcal{T} \Delta D_t^{r-2} p=   D_t^{r-1} {\rm div} v +$other terms ($r\ge 3$) that
\bee\label{}\begin{split}
  \|\pl^{2}D_t^{r-2} p \|_{L^2(\mathscr{D}_t)}^2
  \le & C(K,{\rm Vol}\mathscr{D}_t) \|\Delta D_t^{r-2} p \|_{L^2(\mathscr{D}_t)}^2  \\
  \le & C_r(\cdot) \|D_t^{r-1} {\rm div} v\|_{L^2(\mathscr{D}_t)}^2 + \textrm{other terms} \\
  \le & C_r(\cdot) \|\pl D_t^{r-1} {\rm div} v\|_{L^2(\mathscr{D}_t)}^2 + \textrm{other terms}\\
  =& C_r(\cdot) E_r^d+ \textrm{other terms}.
\end{split}\eee
It follows from the equation $ -\mathcal{T} \Delta D_t^{r-3} p=   D_t^{r-2} {\rm div} v +$other terms ($r\ge 4$) that
\bee\label{}\begin{split}
&  \|\pl^{3}D_t^{r-3} p \|_{L^2(\mathscr{D}_t)}^2 \\
  \le & C \| \Pi\pl^{3}D_t^{r-3} p \|_ {L^2(\pl\mathscr{D}_t)}^2 + C(K,{\rm Vol}\mathscr{D}_t) \sum_{s=0}^1 \|\pl^s \Delta D_t^{r-3} p \|_{L^2(\mathscr{D}_t)}^2   \\
  \le & C \| \pl_\mathcal{N} D_t^{r-3} p \|_ {L^\iy(\pl\mathscr{D}_t)}^2 \| \bar\pl \ta \|_ {L^2(\pl\mathscr{D}_t)}^2   + C_r(\cdot) \|\pl D_t^{r-2} {\rm div} v\|_{L^2(\mathscr{D}_t)}^2 + \textrm{other terms} \\
  \le & C_r(\cdot) \lt(E_{r-1}^a + E_{r-1}^d\rt) + \textrm{other terms} .
\end{split}\eee
It should be pointed out that we need the bound of $\| \pl_\mathcal{N} D_t^{r-3} p \|_ {L^\iy(\pl\mathscr{D}_t)}$ to control $\|\pl^{3}D_t^{r-3} p \|_{L^2(\mathscr{D}_t)}$ by  energy functionals. This is why the a priori assumptions \ef{apriori17} we made for the case of $n=3$ include  the bound of $\| \pl_\mathcal{N} D_t^{2} p \|_ {L^\iy(\pl\mathscr{D}_t)}$, which was not needed in \cite{CL00}.
Similarly, it can be seen from \ef{9/7-2} that
 the a priori assumption on the bound of $\| \pl_\mathcal{N} D_t^{2} {\rm div} v \|_ {L^\iy(\pl\mathscr{D}_t)}$ is also needed when $n=3$. The verification of the a priori assumptions on these bounds is difficult, even on that for
$\| \pl_\mathcal{N} D_t  p \|_ {L^\iy(\pl\mathscr{D}_t)}$, which will be discussed later.

We may conclude, under the a priori assumptions \ef{apriori17}, that there are continuous functions $C_r$ ($0\le r\le n+2$) such that
\bee
\frac{d}{dt} E_0(t)\le
     \|p\|_{L^2(\mathscr{D}_t)}^2 +   \|{\rm div} v\|_{L^2(\mathscr{D}_t)}^2  \le C_0\lt( \overline{V},    M, \underline{\mathcal{T}}_0^{-1}, \overline{\mathcal{T}}_0 \rt) \lt(E_1(t)+E_2(t)\rt) ,
\eee
\bee
\frac{d}{dt}  E_1(t) \le C_1 \lt( \overline{V},K,   M, \underline{\mathcal{T}}_0^{-1}, \overline{\mathcal{T}}_0  \rt)\lt(E_1(t)+E_2(t)\rt),
\eee
\bee
 \frac{d}{dt}  E_2(t) \le C_2 \lt( \overline{V},K,  \ea_b^{-1},L, M, \widetilde{M}, \underline{\mathcal{T}}_0^{-1}, \overline{\mathcal{T}}_0  \rt)\lt(E_1(t)+E_2(t)\rt),
 \eee
\bee
 \frac{d}{dt}  E_r(t) \le C_r \lt( \overline{V},K,  \ea_b^{-1},L, M, \widetilde{M}, \underline{\mathcal{T}}_0^{-1}, \overline{\mathcal{T}}_0, \sum_{s=1}^{r-1} E_s(t)  \rt)\sum_{s=1}^r E_s(t), \ \   3 \le r\le n+2.
\eee

In order to close the arguments, we need to get the estimates for the a priori bounds in terms of the energy functionals $E_r$ ($0\le r\le n+2$),  for which the clear and detailed dependence of $C_r(\cdot)$ on the quantities in the a priori assumptions is crucial. The lower and upper bounds for ${\rm Vol} \mathscr{D}_t $, the $L^{\infty}$-bound  for $\ta$ and the lower bound for $\pl_{\mathcal{N}} p$ on $\pl \mathscr{D}_t$, and the $L^{\infty}$-bound for $\pl(p, v, \mathcal{T})$ in $  \mathscr{D}_t$
can be obtained by looking at the  evolution  of these quantities.  The estimate  for the lower bound for  $\iota_0$
follows from the same idea as in \cite{CL00}. The estimates for other quantities in the a priori assumptions follow from Sobolev lemmas, the projection formula, and elliptic estimates. Here we point out some main differences compared with \cite{CL00}.
The estimate on  $\|\pl_\mathcal{N} D_tp \|_{L^\iy(\pl\mathscr{D}_t)}$ given by \cite{CL00} cannot work for our problem. Indeed, the bound for $\|\pl_\mathcal{N} D_tp \|_{L^\iy(\pl\mathscr{D}_t)}$ was obtained in \cite{CL00} by use of
the following fact: If $q=0$ on $\pl\mathscr{D}_t$, then
\be\label{9/9-1}\begin{split}
\|\pl_\mathcal{N} q \|_{L^\iy(\pl\mathscr{D}_t)} \le  C \|\pl^{n-1}\Delta q\|_{L^2(\mathscr{D}_t)} + C \lt(K,{\rm Vol}\mathscr{D}_t,\|\ta \|_{L^2(\pl\mathscr{D}_t)}, \cdots, \rt. \\
\lt. \|\bar\pl^{n-2}\ta \|_{L^2(\pl\mathscr{D}_t)} \rt)\sum_{s=0}^{n-2}\|\pl^{s}\Delta q\|_{L^2(\mathscr{D}_t)},  \ \ n=2,3.
\end{split}\ee
This can be found in Proposition 5.10 of \cite{CL00}.  If we apply \ef{9/9-1} to our problem, we get
 \be\label{9/9-2}\begin{split}
 \|\pl_\mathcal{N} D_tp \|_{L^\iy(\pl\mathscr{D}_t)}  \le & C   \|\mathcal{T}^{-1}\|_{L^\iy (\mathscr{D}_t)}\|\pl \mathcal{T} \|_{L^\iy (\mathscr{D}_t)}  \|\bar\pl \ta \|_{L^2(\pl\mathscr{D}_t)} \|\pl_\mathcal{N} D_tp \|_{L^\iy(\pl\mathscr{D}_t)}\\
  &  + \textrm{other terms}, \ \ n=3,
\end{split}\ee
which cannot give the bound of $\|\pl_\mathcal{N} D_tp \|_{L^\iy(\pl\mathscr{D}_t)}$ for our problem. In fact, \ef{9/9-2} follows  from the equation $\mathcal{T} \Delta D_t p = - (\pl \mathcal{T}) \cdot \pl  D_t p $+other terms, and
\bee\label{}\begin{split}
  \|\pl^2\Delta D_tp \|_{L^2(\mathscr{D}_t)} \le \|\mathcal{T}^{-1}\|_{L^\iy (\mathscr{D}_t)}\|\pl \mathcal{T} \|_{L^\iy (\mathscr{D}_t)} \|\pl^3   D_tp \|_{L^2(\mathscr{D}_t)} + \textrm{other terms} \\
  \le \|\mathcal{T}^{-1}\|_{L^\iy (\mathscr{D}_t)}\|\pl \mathcal{T} \|_{L^\iy (\mathscr{D}_t)} \|\pl_\mathcal{N} D_tp \|_{L^\iy(\pl\mathscr{D}_t)} \|\bar\pl \ta \|_{L^2(\pl\mathscr{D}_t)} + \textrm{other terms}.
\end{split}\eee
To overcome the difficulty appearing in \ef{9/9-2}, we refine \ef{9/9-1} to show: If $q=q_b$ on $\pl\mathscr{D}_t$ with $q_b$ being a constant, then for any $\delta>0$,
\be\label{9/9-3}\begin{split}
\|\pl_\mathcal{N} q \|_{L^\iy(\pl\mathscr{D}_t)} \le  \delta \|\pl^{n-1}\Delta q\|_{L^2(\mathscr{D}_t)} + C \lt(\delta^{-1}, K,{\rm Vol}\mathscr{D}_t,\|\ta \|_{L^2(\pl\mathscr{D}_t)}, \cdots, \rt. \\
\lt. \|\bar\pl^{n-2}\ta \|_{L^2(\pl\mathscr{D}_t)} \rt)\sum_{s=0}^{n-2}\|\pl^{s}\Delta q\|_{L^2(\mathscr{D}_t)},  \ \ n=2,3.
\end{split}\ee
Clearly, the bound for $\|\pl_\mathcal{N} D_tp \|_{L^\iy(\pl\mathscr{D}_t)}$ can be obtained by use of \ef{9/9-3}.

 \section{Preliminaries}\label{sec2}
In this section, we introduce Lagrangian transformation, the metric and covariant differentiation associated with it,  the
induced metric on the boundary, the geometry and regularity of the boundary, Sobolev lemmas, interpolation inequalities and estimates for the boundary. Those materials are basically from \cite{CL00}. We list them here for the convenience of readers and
the easier reference.

\subsection{Lagrangian coordinates, the metric, and covariant differentiation in the interior}
Let  $x= x(t,y) $ be the change of variables given by
\be\label{s2.1}
\pl_t x(t,y) =   v \left(t,  x(t,y)\right)  \ \  {\rm and} \ \ x(t,y) =x_0(y), \ \  y \in \Omega.
\ee
Initially, when $t=0$, we can start with  either the Euclidean coordinates in $\Omega= \mathscr{D}_0$ or some other coordinates $x_0: \Omega\to  \mathscr{D}_0$ where $x_0$ is a diffeomorphism in which the domain $\Omega$ becomes simple. For each $t$ we will then have a change of coordinates $x: \Omega \to \mathscr{D}_t $, taking $y \to x(t,y)$.
The Euclidean metric
$\da_{ij}$ in $ \mathscr{D}_t$ then induces a metric
\be\label{s2.2}g_{ab}(t,y)=\da_{ij}\frac{\pl x^i}{\pl y^a}\frac{\pl x^j}{\pl y^b}
\ee
on $\Omega$ for each fixed $t$.

The {\it covariant differentiation} of a $(0,r)$ tensor $w(t,y)$, is the $(0,r+1)$ tensor given by
$$\na_a w_{a_1\cdots a_r}=\frac{\pl w_{a_1\cdots a_r}}{\pl y^a} - \Gamma_{aa_1}^e w_{ea_2\cdots a_r}-\cdots- \Gamma_{aa_r}^e w_{ a_1\cdots a_{r-1} e},$$
where $\Gamma^{c}_{ab}$ are the Christoffel symbols given by
$$\Gamma^{c}_{ab}=\frac{g^{cd}}{2} \lt( \frac{\pl g_{bd}}{\pl y^a} + \frac{\pl g_{ad}}{\pl y^b} - \frac{\pl g_{ab }}{\pl y^d} \rt)= \frac{\pl y^c}{\pl x^i} \frac{\pl^2 x^i}{ \pl y^a \pl y^b} $$
with $g^{ab}$ being the inverse of $g_{ab}$. If $\oa(t,x)$ is the $(0,r)$ tensor expressed in the $x$-coordinates, then the same tensor $w(t,y)$ expressed in the $y$-coordinates is given by
$$w_{a_1\cdots a_r}(t,y)=\frac{\pl x^{i_1}}{\pl y^{a_1}}\cdots \frac{\pl x^{i_r}}{\pl y^{a_r}} \oa_{i_1\cdots i_r}(t,x), \ \  x=x(t,y),$$
and by the transformation properties for tensors,
$$\na_a w_{a_1\cdots a_r}(t,y) =\frac{\pl x^{i}}{\pl y^{a}} \frac{\pl x^{i_1}}{\pl y^{a_1}}\cdots \frac{\pl x^{i_r}}{\pl y^{a_r}} \frac{\pl\oa_{i_1\cdots i_r}(t,x) }{\pl x^i} . $$
So that the norms of tensors are invariant under change of coordinates:
\begin{align}
g^{a_1b_1}\cdots  g^{a_rb_r} w_{a_1\cdots a_r}w_{b_1\cdots b_r}
=\da^{i_1j_1}\cdots  \da^{i_rj_r} \oa_{i_1\cdots i_r}\oa_{j_1\cdots j_r}. \label{sss2.3}
\end{align}
Since the curvature vanishes in the $x$-coordinates, it must do so in the $y$-coordinates, and hence
$$[\na_a, \na_b]=0.$$
Set
$${w_{a\cdots}}\ ^b \ _{\cdots c}=g^{bd}w_{a\cdots d\cdots c}.$$

The material derivative is defined as
$$
D_t= \pl_t |_{x=const.} + v^k \pl_k = \pl_t |_{y=const.} ,  \ \   \pl_k =\frac{\pl}{\pl x^k} = \frac{\pl y^a}{\pl x^k}\frac{\pl}{\pl y^a}.
$$
Let $\al$ be a $(0,s)$ tensor and $\ba$ be a $(0,r)$ tensor. Then $\al\widetilde{\otimes}\ba$ is used to denote some partial symmetrization of the tensor product $\al\otimes \ba$, i.e., a sum over some subset of the permutations of the indices divided by the number of permutations in that subset. Moreover $\al \ \widetilde{\cdot}\ \ba$ is used to denote a partial symmetrization of the dot product $\al{\cdot}\ba$, which in turn is defined to be a contraction of the last index of $\al$ with the first index of $\ba$: $(\al\cdot\ba)_{i_1\cdots i_{r+s-2}}=g^{ij}\al_{i_1\cdots i_{s-1}i}\ba_{ji_s\cdots i_{r+s-2}}$.

The following lemmas are for  temporal derivatives of the change of coordinates and commutators between temporal derivative and spatial derivatives, which are Lemmas 2.1 and 2.4 in \cite{CL00}, and  will be used to calculate the higher order equations in Lagrangian coordinates.

 \begin{lem} Let $x=x(t,y)$ be the change of variables given by \ef{s2.1},  and let $g_{ab}$ be the metric given by \ef{s2.2}. Let $v_i=\da_{ij}v^j=v^i$ and $d\mu_g=\sqrt{{\rm det}g} dy$. Set
$$u_a(t,y)= \frac{\pl x^j}{\pl y^a} v_j(t,x), \ \ h_{ab}= (\na_a u_b + \na_b u_a)/2 , \ \ h^{ab}=g^{ac} g^{bd} h_{cd}, \ \ {\rm div} u=g^{ab}\na_a u_b.$$
 Then,
\begin{align}
& D_t \frac{\pl x^i}{\pl y^a} = \frac{\pl x^k}{\pl y^a}\frac{\pl v_i}{\pl x^k}, \ \
D_t \frac{\pl y^a}{ \pl x^i}=- \frac{\pl y^a}{\pl x^k} \frac{\pl v_k}{\pl x^i}, \label{CL-2.12}\\
& D_t g_{ab}=  2 h_{ab},  \ \ D_t g^{ab}=-2 h^{ab}   ,  \ \  D_t d\mu_g =  g^{ab} h_{ab} d\mu_g=({\rm div} u) d\mu_g.
 \label{CL-2.13}
 \end{align}
\end{lem}

\begin{lem} Let $w_{a_1\cdots a_r}$ be a (0,r) tensor, $q$ be a function, and
  $\Delta=g^{cd}\na_c\na_d$. Then,
 \begin{align}
&[D_t, \na_a ] w_{a_1\cdots a_r} = - (\na_{a_1}\na_a u^e)w_{ea_2\cdots a_r}-\cdots-(\na_{a_r}\na_a u^e)w_{a_1\cdots a_{r-1} e},
\label{CL-2.20} \\
&[D_t, \Delta] q= - 2 h^{ab} \na_a\na_b q -(\Delta u^e) \na_e q,  \label{CL-2.22}
\end{align}
Furthermore,
\begin{equation}\label{CL-2.23}\lt.\begin{split}
[D_t, \na^r ] q=- \sum_{s=1}^{r-1} \lt(\begin{split} &\ \ r \\ &s+1 \end{split}\rt) (\na^{s+1} u)  \cdot \na^{r-s} q.
\end{split}\rt.\end{equation}
\end{lem}

\subsection{The geometry and regularity of the boundary}
As in \cite{CL00}, we extend the normal to the boundary to the interior by a geodesic extension, which enables us to define a pseudo-Riemann metric in the whole domain whose restriction on the boundary  is then the induced metric on the tangential space to the boundary. Using this induced metric, we can define  the orthogonal projection of a tensor to the boundary, the covariant differentiation on the boundary, and the second fundamental form of the boundary as follows:

\begin{defn}\label{gama} Let $d(t,y)={\rm dist}_g (y, \pl\Omega)$ be the geodesic distance to the boundary, which is the same as the Euclidean distance in the $x$-variables, and $\eta$ be the smooth cut-off function given by Definition \ref{defn1.3}. Set $N_a(t,y)=\na_a d(t,y)$ and $N^a(t,y)=g^{ab}(t,y)N_b(t,y)$. Define
$$
 \zeta^{ab}(t,y)=g^{ab}(t,y)-    \widetilde{N}^a(t,y) \widetilde{N}^b(t,y), \ \ {\rm where} \ \   \widetilde{N}^a(t,y)=\eta (d(t,y) )  {N}^a(t,y).
$$
In particular, $\zeta$ gives the induced metric on the tangent space to the boundary:
$$\zeta_{ab}=g_{ab}-N_aN_b \ \ {\rm and} \ \ \zeta^{ab}=g^{ab} -N^a N^b \ \  {\rm on} \ \ \pl\Oa. $$
The orthogonal projection of a $(0,r)$ tensor $w (t,y)$  to the boundary is given by
$$(\Pi w)_{a_1\cdots a_r} = \zeta_{a_1}^{c_1}\cdots \zeta_{a_r}^{c_r} w_{c_1\cdots c_r},
\ \ {\rm where} \ \ \zeta_a^c = \da_a^c -N_aN^c.$$
The covariant differentiation on the boundary $\overline{\na}$ is given by $\overline{\na}_a = \zeta_a^c \na_c$. The second fundamental form of the boundary is given by $\ta_{ab}(t,y)=\overline{\na}_a N_b$.
\end{defn}
It follows from Definitions \ref{defn1.1}, \ref{defn1.3} and \ref{gama} that
$$
 N_a(t,y)=  \frac{\pl x^j}{\pl y^a} \mathcal{N}_j (t,x)  \ \ {\rm and} \ \  \ta_{ab}(t,y)= \frac{\pl x^i}{\pl y^a}\frac{\pl x^j}{\pl y^b} \ta_{ij}(t,x).
$$

\begin{lem}  Let $N$ be the unit normal to $\pl\Omega$ and $d\mu_\za =\sqrt{{\rm det}g/ (\sum N_a^2) } dS $ with $dS$ being the Euclidean surface measure. On $[0,T]\times\pl\Omega$ we have
\begin{align}
&D_t N_a=h_{NN} N_a, \ \  D_t N^c=-2 h^{c}_d N^d + h_{NN} N^c, \ \ {\rm where} \ \ h_{NN}=N^aN^b h_{ab}, \label{CL-3.21} \\
&D_t \za^{ab}=-2 \za^{ac}\za^{bd}h_{cd},  \label{CL-3.22}\\
& D_t d\mu_\za=\lt(g^{ab}h_{ab}-h_{NN}\rt)d\mu_\za =\lt({\rm div} u - h_{NN}\rt)d\mu_\za . \label{CL-3.23}
\end{align}
\end{lem}
This is Lemma 3.9 in \cite {CL00}, where the proof can be found.

\begin{defn}\label{defn2.4} For the multi-indices $I=(i_1, \cdots, i_r)$ and $J=(j_1, \cdots, j_r)$, set
$g^{IJ}=g^{i_1j_1} \cdots g^{i_r j_r}$ and   $\zeta^{IJ}=\zeta^{i_1j_1} \cdots \zeta^{i_r j_r} $ .
If $\alpha$ and $\beta$ are $(0,r)$ tensors, let
\begin{align*}
&<\alpha, \beta>=g^{IJ} \alpha_I \beta_J  \ \ {\rm and} \  \  |\alpha|^2=<\alpha, \alpha>= g^{IJ} \alpha_I \alpha_J .
\end{align*}
Then for the projection $(\Pi \beta)_I=\zeta_I^J \beta_J$,
$$
<\Pi\alpha, \Pi\beta>=\zeta^{IJ} \alpha_I \beta_J \  \ {\rm and} \ \ |\Pi\alpha|^2=\zeta^{IJ} \alpha_I \alpha_J  \ \  {\rm on} \ \ \pl\Omega .
$$
Let
$$
 \| \alpha  \| = \| \alpha  \| _{L^2(\Omega)}   = \lt( \int_\Omega |\alpha|^2 d\mu_g \rt)^{1/2}  \ \ {\rm and} \ \    |\| \alpha \| |    = \|  \alpha \| _{L^2(\pl\Omega)}  = \lt( \int_{\pl\Omega} |\alpha|^2  d\mu_\zeta \rt)^{1/2}.
$$
Moreover, we define the following notations:
$$\| \cdot  \| _{L^p}=\| \cdot  \| _{L^p(\Omega)} \ \ {\rm and} \ \
|\| \cdot  \||_{L^p}=|\| \cdot  \||_{L^p(\pl\Omega)} \ \  {\rm for} \ \ 1\le p \le \iy.$$
\end{defn}

\begin{lem}   With the notations in Definitions \ref{defn1.2}, \ref{gama} and \ref{defn2.4}, we have
\begin{align}
\| \na \za\|_{L^\iy} \le 512\lt(|\| \theta \||_{L^\iy}+ 1/\iota_0\rt) \ \
{\rm and} \ \  \|D_t \zeta\|_{L^\iy} \le 128 \|h\|_{L^\iy}. \label{CL-3.28}
\end{align}
\end{lem}
This is Lemma 3.11 in \cite{CL00}, where the proof can be found.

The following Lemma which is Lemma 3.6 in \cite{CL00} shows that $\iota_1$ given in Definition \ef{defn2.3} is equivalent to $\iota_0$ in conjunction with a bound of the second fundamental form.
\begin{lem}\label{lemkk1}
Suppose that $|\ta|\le K$, and let $\iota_0$ and $\iota_1$ be as in  Definitions \ref{defn1.2} and \ref{defn2.3}. Then
$$
\iota_0\ge \min\{\iota_1/2, \ 1/K\} \ \ {\rm and} \ \  \iota_1 \ge \min\{2\iota_0, \ \ea_1/K\}.
$$
\end{lem}
The advantage to using $\iota_1$, instead of $\iota_0$, is that it is easier to control the evolution off.

\subsection{Sobolev lemmas, interpolation inequalities and estimates for the boundary}
\begin{lem} {\rm(}{\rm Lemmas} {\rm A.1-A.4} {\rm in} {\rm \cite{CL00}}{\rm )} Let $\al$ be a $(0,r)$ tensor and $\iota_1\ge 1/K_1$. Assume   $k$, $m$  are positive integers, and $p\ge 1$.  We have\\
 (i) if $2\le p\le s\le q\le \iy$ and $m/s=k/p+(m-k)/q$,
\begin{align}
  &     |\| \overline{\na}^k \al \||_{L^s}^{m} \le C(k,m,n,s) |\|  \al \||_{L^q}^{m-k}
|\| \overline{\na}^m \al \||_{L^p}^{k}\label{CL-A.4}, \\
& \lt(\sum_{i=0}^k \| \na^i \al \|_{L^s} \rt)^{m} \le C(k,m,n,s)  \|  \al \|_{L^q}^{m-k}   \lt(\sum_{i=0}^m K_1^{m-i} \| \na^i \al \|_{L^p} \rt)^{k}; \label{CL-A.12}
\end{align}
(ii) for any $\da>0$,
\begin{align}
&|\| \al \||_{L^{{(n-1)p}/{(n-1-kp)}}} \le C(k,n,p)  \sum_{i=0}^k K_1^{k-i}  |\|\na^i \al \||_{L^{p}}, \ \    1\le p < (n-1)/k, \label{CL-A.7} \\
& |\| \al \||_{L^{\iy}} \le \da |\| \na^k \al \||_{L^{p}} +  C(\da^{-1}, K_1, k, n,p) \sum_{i=0}^{k-1} |\| \na^i \al \||_{L^{p}},   \ \     p>(n-1)/k, \label{CL-A.8} \\
&\| \al \|_{L^{{np}/{(n-kp)}}} \le C(k,n,p)  \sum_{i=0}^k K_1^{k-i} \|\na^i \al \|_{L^{p}},  \ \   1\le p < n/k, \label{CL-A.14} \\
&\| \al \|_{L^{\iy}} \le C(k,n,p)  \sum_{i=0}^k K_1^{k-i} \|\na^i \al \|_{L^{p}},  \ \     p > n/k. \label{CL-A.15}
\end{align}
\end{lem}

\begin{lem}\label{2ndfun} Let $q=q_b$ on $\pl \Omega$ with $q_b$ being a constant, then for $r=2,3,4$,
\begin{align}
 |\|\Pi \na^r q \| | \le & 2 |\|\na_N q\||_{L^\iy} |\|\overline{\na}^{r-2}\ta\|| + C \sum_{k=1}^{r-1} |\| \ta\||_{L^\iy}^k |\| \na^{r-k}  q\||
 \notag  \\
& +C \sum_{k=1}^{r-3} |\|    \ta  \| |_{L^\iy} |\|   \na_N q  \| |_{L^\iy} |\|  \overline{\na}^{k} \ta  \| | ,
 \label{hb1}
\end{align}
where $C$ is a positive number. If, in addition, $|\na_N q|\ge \ea$ and $|\na_N q| \ge 2 \ea |\|\na_N q\||_{L^\iy}$ on $\pl \Omega$ for a certain positive constant $\ea$, then there exists a positive number $C$ such that
\begin{align}
 |\| \overline{\na}^{r-2}\ta\| | \le & \ea^{-2}    |\|\Pi \na^r q\|| + C \ea^{-3} \sum_{k=1}^{r-1} |\| \ta\||_{L^\iy}^k |\| \na^{r-k}  q\||  \notag \\
& +C  \ea^{-2} \sum_{k=1}^{r-3} |\|    \ta  \| |_{L^\iy} |\|   \na_N q  \| |_{L^\iy} |\|  \overline{\na}^{k} \ta  \| | ,     \ \ r=2,3,4.
\label{hb2}
\end{align}
\end{lem}

{\em Proof}.  Simple calculations give that on $\pl\Omega$,
\begin{align}
\Pi \na^2 q= &  ( \na_N q) \ta ,  \ \
\Pi \na^3 q=   ( \na_N q) \overline{\na} \ta + 3 ( \overline{\na} \na_N q) \widetilde{\otimes}  \ta  +  2(\ta \ \widetilde{\cdot}\ \ta)  \widetilde{\otimes}  \na q,  \label{hb3} \\
 \Pi \na^4 q=&( \na_N q) \overline{\na}^2 \ta +  4 (\overline{\na} \ta) \widetilde\otimes \overline{\na}\na_N q +6 \ta \widetilde\otimes \overline{\na}^2\na_N q  + 3 (\ta \widetilde\otimes \ta ) \na_N^2 q  \notag\\
& + 7 ( (\overline{\na} \ta)   \ \tilde\cdot \  \ta   )\widetilde\otimes \na q
 - \ta \widetilde\otimes ( \ta \ \tilde\cdot \  \ta) \na_N q +  3 (( \ta   \ \tilde\cdot \  \ta   )  \ \tilde\cdot \  \ta )\widetilde\otimes N  \widetilde\otimes \na q  \notag\\
& + 8 ( \ta \ \tilde\cdot \  \ta) \widetilde\otimes ( \overline{\na} \na_N q ) \widetilde\otimes N  + 3 \ta \widetilde\otimes ( \ta \ \tilde\cdot \    \overline{\na} \na_N q ) \widetilde\otimes N , \notag
\end{align}
where we have used the facts that $\overline{\na} q =0$, $\ta\cdot N=0$, $ (\overline{\na} \ta)\cdot N=- \ta \cdot  \ta$ and $(\overline{\na}^2 \ta)\cdot N = - 3  (\overline{\na} \ta)  \ \tilde\cdot \    \ta $. Clearly, \ef{hb1} and \ef{hb2} hold for $r=2,3$, because of
$ \overline{\na} \na_N q= N^e \overline{\na} \na_e q$ and
\begin{align}
|\|\Pi \na^3 q - ( \na_N q) \overline{\na} \ta\| |
\le  3 |\| \ta\||_{L^\iy} |\| \na^2 q\||  +  2 |\| \ta|\|_{L^\iy}^2 |\| \na q |\|\label{hb4}.
\end{align}
For $r=4$, we first derive from  \ef{CL-A.4}, H$\ddot{\rm o}$lder's inequality and Young's inequality that  for any positive constant $\da$,
\begin{align*}
|\| |\overline{\na} \ta|    |\overline{\na}\na_N q |    \| |
\le &  |\|  \overline{\na} \ta  \| |_{L^4}  |\| \overline{\na}\na_N q  \| |_{L^4} \\
 \le & C |\|   \ta  \| |_{L^\iy}^{1/2}  |\|  \overline{\na}^2 \ta  \| |^{1/2}
 |\|  \na_N q  \| |_{L^\iy}^{1/2}  |\| \overline{\na}^2\na_N q  \| |^{1/2} \\
 \le& ( \da /4) |\|\na_N q\||_{L^\iy} |\|\overline{\na}^2 \ta\||
 + C \da^{-1} |\|   \ta  \| |_{L^\iy} |\| \overline{\na}^2\na_N q  \| |;
\end{align*}
which, together with $\na_N^2 q = N^e \na_N \na_e q$ and
$|\overline{\na}^2 \na_N q|\le |\na^3 q|+ 3|\ta||\na^2 q|$, gives that  for any positive constant $\da$,
\begin{align}
|\| \Pi \na^4 q - & ( \na_N q) \overline{\na}^2 \ta \| |
\le  \da |\|\na_N q\||_{L^\iy} |\|\overline{\na}^2 \ta\|| \notag \\
& + 7 |\|    \ta  \| |_{L^\iy} |\|   \na_N q  \| |_{L^\iy} |\|  \overline{\na} \ta  \| |
  + C\da^{-1} \sum_{k=1}^{3} |\| \ta\||_{L^\iy}^k |\| \na^{4-k}  q\||. \label{hala}
\end{align}
Here $C$ is a positive number. Clearly, choose $\da=1$ in \ef{hala} to prove
\ef{hb1}. Note that
\begin{align*}
2\ea |\|\na_N q\||_{L^\iy} |\|\overline{\na}^2 \ta\|| \le   |\| (\na_N q )\overline{\na}^2 \ta\|| \le |\| \Pi \na^4 q \| | + |\| \Pi \na^4 q -   ( \na_N q) \overline{\na}^2 \ta \| |,
\end{align*}
we then prove \ef{hb2} by choosing $\da=\ea$ in \ef{hala}.
\hfill$\Box$

\begin{lem}\label{172ndfun} Let $q=q_b$ on $\pl \Omega$ with $q_b$ being a constant, then
\begin{align}
 |\|\Pi \na^5 q \| | \le & 2 |\|\na_N q\||_{L^\iy} |\|\overline{\na}^{3}\ta\||    + C \sum_{k=1}^{4} |\| \ta\||_{L^\iy}^k |\| \na^{5-k}  q\|| \notag\\
& +C(K_1, |\| \ta\||_{L^\iy})(|\|\overline{\na}^{2}\ta\||+ |\|\overline{\na} \ta\||)\sum_{k=1}^{4}  |\| \na^{k}  q\||
 \label{17hb1}.
\end{align}
If, in addition, $|\na_N q|\ge \ea$ and $|\na_N q| \ge 2 \ea |\|\na_N q\||_{L^\iy}$ on $\pl \Omega$ for a certain positive constant $\ea$,
\begin{align}
 |\| \overline{\na}^{3}\ta\| | \le & \ea^{-2}    |\|\Pi \na^5 q\|| +  \ea^{-3} C \sum_{k=1}^{4} |\| \ta\||_{L^\iy}^k |\| \na^{5-k}  q\|| \notag\\
& +\ea^{-3}C(K_1, |\| \ta\||_{L^\iy})(|\|\overline{\na}^{2}\ta\||+ |\|\overline{\na} \ta\||)\sum_{k=1}^{4}  |\| \na^{k}  q\||.
\label{17hb2}
\end{align}
\end{lem}

{\em Proof}. This lemma can be shown in a similar way to proving Lemma \ref{2ndfun} by noticing the following fact:
\be\label{17331}\begin{split}
& |\| |\overline{\na} \ta|    |\overline{\na}^2 \na_N q |    \| |
+ |\| |\overline{\na}^2 \ta|    |\overline{\na}  \na_N q |    \| | \\
\le &  |\|  \overline{\na} \ta  \| |_{L^6}  |\| \overline{\na}^2 \na_N q  \| |_{L^3} + |\|  \overline{\na}^2 \ta  \| |_{L^3}  |\| \overline{\na}  \na_N q  \| |_{L^6} \\
 \le & C |\|   \ta  \| |_{L^\iy}^{2/3}  |\|  \overline{\na}^3 \ta  \| |^{1/3}
 |\|  \na_N q  \| |_{L^\iy}^{1/3}  |\| \overline{\na}^3\na_N q  \| |^{2/3}
 \\
 & +C |\|   \ta  \| |_{L^\iy}^{1/3}  |\|  \overline{\na}^3 \ta  \| |^{2/3}
 |\|  \na_N q  \| |_{L^\iy}^{2/3}  |\| \overline{\na}^3\na_N q  \| |^{1/3} \\
 \le&  \da |\|\na_N q\||_{L^\iy} |\|\overline{\na}^3 \ta\||
 + C \da^{-1} |\|   \ta \||_{L^\iy} |\| \overline{\na}^3\na_N q  \| |
\end{split}\ee
for any positive constant $\da$, and
\be\label{17331'}\begin{split}
 |\|    | \na^2 q | |\overline{\na} \ta|    \| |
\le    |\|  {\na}^2   q  \| |_{L^\iy}   \|| \overline{\na} \ta  \|| \le C(K_1) |\|  \overline{\na} \ta  \|  \sum_{k=2}^{4}  |\| \na^{k}  q\||.
\end{split}\ee
Here \ef{17331} follows from  \ef{CL-A.4}, H$\ddot{\rm o}$lder's inequality and Young's inequality, and  \ef{17331'} follows from  \ef{CL-A.8}.
\hfill$\Box$

\begin{rmk}  $\ea$ appearing on the right-hand side of \ef{hb2} and \ef{17hb2} can be chosen as
\begin{align}
\ea= |\| (\na_N q)^{-1} \| |_{L^\iy}^{-1} \min\lt\{1,  \   2^{-1} |\|  \na_N q  \| |_{L^\iy}^{-1}\rt\}. \label{rt1}
\end{align}
In particular, it follows from \ef{hb3} and \ef{hb4} that
\begin{align}
& |\|  \ta\| |_{L^s} \le  |\| (\na_N q)^{-1} \| |_{L^\iy}    |\|  \Pi \na^2 q\| |_{L^s} ,  \ \ 2\le s \le \iy, \label{4.25.1}\\
& |\| \overline{\na} \ta\| | \le |\| (\na_N q)^{-1} \| |_{L^\iy} \lt(   |\|  \Pi \na^3 q\| |  + 3 \sum_{k=1}^{2} |\| \ta\||_{L^\iy}^k |\| \na^{3-k}  q\|| \rt). \label{4.25.2}
\end{align}
\end{rmk}

\begin{lem}\label{lem5.5} {\rm(}{\rm Lemmas} {\rm 5.5-5.6} {\rm in} {\rm \cite{CL00}}{\rm )}
Let $w$ be a $(0,1)$ tensor and define  a scalar ${\rm div} w=g^{ab}\na_a w_b$ and a $(0, 2)$ tensor ${\rm curl} w_{ab}=\na_a w_b -\na_b w_a$. If $|\ta|+1/\iota_0\le K$, then for any nonnegative integer $r$,
\begin{align}
&|\na^{r+1} w|^2 \le C\lt(g^{ij} \za^{kl} \za^{IJ} (\na_k \na^r_I w_i) \na_l \na^r_J w_j + |\na^ r {\rm div} w|^2 + |\na^r {\rm curl} w|^2 \rt), \label{CL-5.16} \\
&\|\na^{r+1} w\|^2 \le C\int_{\Omega} \widetilde{N}^i\widetilde{N}^j g^{kl} \za^{IJ}
(\na_k\na^{r }_I w_i)\na_l\na^{r }_J w_j d\mu_g \notag\\
&\qquad \qquad  \quad + C\lt(\|\na^{r}{\rm div} w\|^2 + \|\na^{r}{\rm curl} w\|^2 + K^2  \|\na^{r} w\|^2  \rt) ,\label{CL-5.17}\\
&|\|\na^r w\||^2 \le C \lt(\|\na^{r+1}w\|  + K \|\na^r w\| \rt)\|\na^r w\|,\label{CL-5.19}\\
&|\|\na^r w\||^2 \le C|\|\Pi \na^r w\||^2 + C\lt(\|\na^{r}{\rm div} w\|  + \|\na^{r}{\rm curl} w\|   +  K \|\na^r w\| \rt)\|\na^r w\|,\label{CL-5.20}\\
&\|\na^{r+1} w\|^2 \le C|\|  \na^{r+1} w\||     |\|  \na^{r} w\|| + C\lt(\|\na^{r}{\rm div} w\|^2 + \|\na^{r}{\rm curl} w\|^2   \rt) ,\label{CL-5.21}\\
&\|\na^{r+1} w\|^2 \le C|\| \Pi \na^{r+1} w\|||\| \Pi (N^i   \na^{r} w_i)\|| \notag\\
&\qquad \qquad \quad + C\lt(\|\na^{r}{\rm div} w\|^2 + \|\na^{r}{\rm curl} w\|^2 + K^2  \|\na^{r} w\|^2  \rt) ,\label{CL-5.22}\\
&\|\na^{r+1} w\|^2 \le C|\| \Pi (N^i \na^{r+1} w_i)\|||\| \Pi \na^{r} w\|| \notag\\
&\qquad \qquad  \quad + C\lt(\|\na^{r}{\rm div} w\|^2 + \|\na^{r}{\rm curl} w\|^2 + K^2  \|\na^{r} w\|^2  \rt).\label{CL-5.23}
\end{align}
\end{lem}
Indeed, the proof of \ef{CL-5.16}-\ef{CL-5.17} can also be found in \cite{Lindblad2}.
The proof of \ef{CL-5.19}-\ef{CL-5.23} are based on the divergence theorem, and  \ef{CL-5.22}-\ef{CL-5.23} are based additionally on \ef{CL-5.17}.

\begin{lem}\label{lemA.5} {\rm(}{\rm Lemma} {\rm A.5} {\rm in} {\rm \cite{CL00}}{\rm )} Suppose that $q=0$ on $\pl\Omega$. Then
\begin{align}
\|q\| \le C ({\rm Vol}\Omega)^{1/n} \| \na q \| \ \  {\rm and} \ \
\|\na q\| \le C ({\rm Vol}\Omega)^{1/n} \| \Delta q \| \label{CL-A.17}.
\end{align}
\end{lem}

As a consequence of Lemmas \ref{lem5.5} and \ref{lemA.5}, we have

\begin{coro}\label{prop5.8}
Let $q=q_b$ on $\pl \Omega$ with $q_b$ being a constant. If $|\ta|+1/\iota_0\le K$,  we have  for any $r\ge 2$ and $\da>0$,
\begin{align}
&\|q-q_b\|  \le C({\rm Vol}\Omega)^{1/n} \|\na q\| , \ \  \|\na q\| + \|\na^2  q\| \le C(K, {\rm Vol}\Omega) \|\Delta q\| ,  \label{310-1} \\
&\|\na^r q\| + |\|\na^r q\|| \le C|\| \Pi \na^r q\|| + C(K, {\rm Vol}\Omega) \sum_{s=0}^{r-1} \|\na^s \Delta q \|, \label{CL-5.28}\\
& \|\na^r q\| + |\|\na^{r-1} q\||  \le \da |\| \Pi \na^r q\|| + C(\da^{-1}, K, {\rm Vol}\Omega) \sum_{s=0}^{r-2} \|\na^s \Delta q \|. \label{CL-5.29}
\end{align}
\end{coro}
Clearly, \ef{310-1} is a consequence of \ef{CL-A.17} and \ef{CL-5.23}. The proof of \ef{CL-5.28} and  \ef{CL-5.29} can be found in Proposition 5.8, \cite{CL00}. (Indeed, \ef{CL-5.28} follows from \ef{CL-5.19}-\ef{CL-5.21} and \ef{CL-A.17}, and \ef{CL-5.29} follows from \ef{CL-5.19}, \ef{CL-5.21},  \ef{CL-5.22} and \ef{CL-A.17}.)

\begin{lem}\label{prop5.10} Let $q=q_b$ on $\pl \Omega$ with $q_b$ being a constant. If $|\ta|+1/\iota_0\le K$ and  $\iota_1\ge 1/K_1$, then for any $\da>0$,
\be\label{CL-5.34}
 |\|\na  q\||_{L^\iy} \le   \lt\{ \begin{split}
  & \da \|\na  \Da q \| +  C(\da^{-1}, K,K_1,|\| \ta\||, {\rm Vol}\Oa )  \|  \Da q \| ,    &n=2;\\
  &\da \|\na^2  \Da q \| +  C  (\da^{-1},K,K_1,|\|\overline{\na}\ta\||, {\rm Vol}\Oa )  ( \|\na  \Da q \| +  \|  \Da q \|),    & n=3.
\end{split}\rt.\ee
\end{lem}
{\em Proof}. When $n=3$, it follows from \ef{hb4} and \ef{CL-A.8} that for any $\da>0$,
\begin{align}
|\|\Pi \na^3 q\|| \le \da |\|\na^3 q\|| + C  (\da^{-1},K,K_1,|\|\overline{\na}\ta\|| ) (|\|\na^2 q\|| + |\|\na  q\||). \label{t-4}
\end{align}
In view of \ef{CL-5.29}, \ef{CL-5.19} and \ef{310-1}, we see that for any $\da_1>0$,
\begin{align}
&|\|\na^2 q\||\le \da_1 |\|\Pi \na^3 q\||   + C(\da_1^{-1},K, {\rm Vol}\Oa ) (\|\na \Da q\| +  \| \Da q\| ), \label{t-6}\\
&|\|\na  q\|| \le C(K) (\|\na^2 q\| +  \| \na q\| ) \le C( K, {\rm Vol}\Oa ) \| \Da q\| . \label{t-7}
\end{align}
Substitute \ef{t-6} and \ef{t-7} into \ef{t-4} and choose suitable small $\da_1$ to obtain for any $\da>0$,
\begin{align*}
2^{-1}|\|\Pi \na^3 q\|| \le \da |\|\na^3 q\|| + C  (\da^{-1},K,K_1,|\|\overline{\na}\ta\||, {\rm Vol}\Oa )  ( \|\na  \Da q \| +  \|  \Da q \|). \label{}
\end{align*}
This, together with \ef{CL-5.28}, gives
\begin{align}
|\|\na^3 q\||   \le & C (K ,{\rm Vol}\Oa )   \|\na^2  \Da q \|
  +  C (K,K_1,|\|\overline{\na}\ta\||, {\rm Vol}\Oa )  ( \|\na  \Da q \| +  \|  \Da q \|) \label{CL-5.34'} .
\end{align}
So, \ef{CL-5.34} follows from \ef{CL-A.8}, \ef{CL-5.34'}
\ef{t-6} and \ef{t-7} in the case of $n=3$.  Similarly, \ef{CL-5.34} can be shown when $n=2$. \hfill$\Box$

\section{Higher Order Equations}
Let $u(t,y)$ be the same tensor of the velocity $v(t,x)$ expressed in the $y$-coordinates, i.e.,
$$u_a(t,y)= \frac{\pl x^j}{\pl y^a} v_j(t,x).
$$
Then, system \ef{l2} can be rewritten as
\begin{subequations}\begin{align}
&D_t u_a +   \mathcal{T} \nabla_a  p = (\nabla_a u_c) u^c , \label{e6}\\
& {\rm div } u =    \Delta \mathcal{T}, \ \
   D_t   \mathcal{T}  =    \mathcal{T} \Delta  \mathcal{T}. \label{e7}
\end{align}\end{subequations}
 It follows from $\ef{e6}$  and \ef{CL-2.20} that
\begin{align}
& D_t \nabla_b u_a+     \mathcal{T}  \nabla_b \nabla_a p  = -  (\nabla_b \mathcal{T} ) \nabla_a p  + (\nabla_a u_e ) \nabla_b u^e ,
\label{e25}\end{align}
which implies
 \begin{align}
& D_t {\rm curl} u_{ab}  =      ( \nabla_b  \mathcal{T} ) \nabla_a p -( \nabla_a  \mathcal{T} ) \nabla_b p,   \label{e27}  \\
 & D_t{\rm div} u  + \mathcal{T} \Delta p    =    -  (\nabla  \mathcal{T})\cdot \nabla  p - (\nabla_e u )\cdot\nabla u^e. \label{lappi}
\end{align}
Moreover, we can derive  from  \ef{e7} and \ef{CL-2.22} that

\begin{align}
 D_t  {\rm div } u & -   \mathcal{T} \Delta {\rm div } u \notag\\
= & -   (\Delta u^e )   \nabla_e  \mathcal{T} - 2  g^{ab} (\nabla_e \nabla_a\mathcal{T} )  \nabla_b  u^e +    ({\rm div } u) ^2 + 2 g^{ab}  (  \nabla_a \mathcal{T} )  \nabla_b  {\rm div } u. \label{e9}
\end{align}

\begin{lem}\label{lem6.1} Let $q$ be any given function. Then for any integer $r \ge 2$,
\begin{align}
&\lt|  D_t \nabla^r u +    \mathcal{T}  \nabla^{r +1}    p  \rt| + \lt| D_t \nabla^{r-1} {\rm curl} u \rt| +|\mathcal{T} \na^{r-1} \Delta p + \na^{r-1} D_t{\rm div}  u |  \notag \\
& \qquad\qquad \le
 C \sum_{s=0}^{r-1} \lt( \lt|\nabla^{1+s} u  \rt|   \lt| \nabla^{r-s} u \rt| + \lt|\nabla^{1+s} \mathcal{T} \rt|   \lt| \nabla^{r-s} p \rt|\rt)\label{CL-6.17}, \\
& \lt|  D_t \nabla^r q + (\nabla^r u)\cdot \nabla q - \nabla^r D_t q   \rt| \le  C \sum_{s=1}^{r-2}  \lt|  (\nabla^{1+s} u) \cdot    \nabla^{r-s} q   \rt|, \label{CL-6.18}
\\
&\lt| \Pi( D_t \nabla^r q + (\nabla^r u)\cdot \nabla q - \nabla^r D_t q )  \rt| \le  C \sum_{s=1}^{r-2}  \lt| \Pi( (\nabla^{1+s} u) \cdot    \nabla^{r-s} q  ) \rt|. \label{CL-6.18'}
\end{align}

\end{lem}
{\em Proof}. This lemma can be proved in a similar way to deriving Lemma 6.1 in \cite{CL00}, so we sketch the proof and omit the details. First, we can apply the following fact
\begin{equation*}\label{}\lt.\begin{split}
[D_t, \pl_i]=-(\pl_i v^k) \pl_k \ \  {\rm and} \ \ [D_t, \pl^r]=- \sum_{s=0}^{r-1} \lt(\begin{split} &\ \ r \\ &s+1 \end{split}\rt) (\pl^{1+s} v) \cdot \pl^{r-s}
\end{split}\rt.\end{equation*}
to \ef{l2} and change coordinates to obtain
\begin{equation}\label{l4}\lt.\begin{split}
& D_t \nabla^r u_a+  \nabla^r \lt( \mathcal{T} \nabla_a p \rt)= (\nabla_a u_c- \nabla_c u_a) \nabla^r u^c -\sum_{s=1}^{r-2} \lt(\begin{split} &\ \ r \\ &s+1 \end{split}\rt) \lt(\nabla^{1+s} u  \rt) \cdot \nabla^{r-s} u_a.
\end{split}\rt.\end{equation}
The estimate for ${\rm curl} u$ can be shown similarly.
For any $r\ge 0$, take $\na^r$ of \ef{lappi}  to get
\begin{align*}
\mathcal{T} \na^r \Delta p + \na^r D_t{\rm div}  u = - \lt(\na^r (\mathcal{T}   \Delta p) - \mathcal{T} \na^r  \Delta p \rt)    - \na^r \lt( (\nabla  \mathcal{T})\cdot \nabla  p + (\nabla_e u )\cdot\nabla u^e \rt).
\end{align*}
This proves \ef{CL-6.17}. Clearly, \ef{CL-6.18} and \ef{CL-6.18'} follow directly from \ef{CL-2.23}.
\hfill $\Box$

\begin{lem}\label{lem3.7} For  $r \ge 0$ and $s=0,1,2$, we have
\begin{align}
& |\na^r D_t\mathcal{T}| + |D_t \na^r \mathcal{T} | \le  C |\mathcal{T} \na^r {\rm div} u| +  C \sum_{s=0}^{r-1}|\na^{s+1}u||\na^{r-s}\mathcal{T}|
 \label{r-3.28}, \\
&|\na^r D_t^2 \mathcal{T}| \le C \sum_{s_1=0}^r |\na^{s_1} \mathcal{T} |
\lt(|\na^{r-s_1} D_t {\rm div} u|  + \sum_{s_2=0}^{r-s_1} |\na^{r-s_1-s_2} {\rm div} u||\na^{s_2}{\rm div} u| \rt)    ,    \label{r-3.29}  \\
& |\na^s D_t^3 \mathcal{T}| \le C\sum_{i=0}^{s} |\na^i \mathcal{T}|H_{s-i},  \label{rr-3.30}
\end{align}
where
\begin{align*}
&H_0=|D_t^2 {\rm div} u| +  (|D_t  {\rm div} u|+  |\na u|^2) |\na u|  ,\\
&H_1=|\na D_t^2 {\rm div}u| +|\na D_t {\rm div} u||\na u| +  (|D_t  {\rm div} u|+  |\na u|^2) |\na {\rm div} u|,\\
&H_2=|\na^2 D_t^2 {\rm div} u| + |\na^2  D_t {\rm div} u||\na u| +
|\na   D_t {\rm div} u||\na {\rm div} u|   \notag\\
& \qquad  + |\na^2 {\rm div} u|  (|D_t  {\rm div} u|+  |\na u|^2)  + |\na {\rm div} u|^2 |\na u|.
\end{align*}
\end{lem}
{\em Proof}. It follows from  \ef{e7} that
\begin{align}
&D_t\mathcal{T}= \mathcal{T} {\rm div} u, \ \    D_t^2 \mathcal{T} =  \mathcal{T} \lt( ({\rm div} u)^2 +  D_t  {\rm div} u \rt), \label{r-3.27} \\
&D_t^3 \mathcal{T} =  \mathcal{T} \lt( ({\rm div} u)^3 + 3 ({\rm div} u) D_t {\rm div} u +  D_t^2  {\rm div} u \rt), \notag
\end{align}
which, together with \ef{CL-6.18}, proves  \ef{r-3.28}-\ef{rr-3.30}.   \hfill $\Box$

\begin{lem}\label{lem3.8} For any integer $r\ge 0$, we have
\begin{align}
& |\na^r D_t \na u| \le |\mathcal{T}||\na^{r+2} p| +  C \sum_{s=0}^r \lt(|\na^{s+1}\mathcal{T}||\na^{r+1-s} p| +|\na^{s+1} u| |\na^{r+1-s}u| \rt), \label{bfordtu} \\
&|D_t^2 \nabla u |
  \le  C \sum_{r=1}^2 |\na^{2-r} \mathcal{T}|\lt(|\na^{r} D_t p|
+ \sum_{s=1}^r |\na^s p| |\na^{r+1-s} u| \rt)  + C   |\na u|^3  ,
\label{bfordt2u} \\
& |D_t^2 \na^2 u| \le  C \sum_{r=1}^3 |\na^{3-r} \mathcal{T}|\lt(|\na^{r} D_t p|
+ \sum_{s=1}^r |\na^s p| |\na^{r+1-s} u| \rt)  + C |\na^2 u||\na u|^2 , \label{q-3.29} \\
&|D_t^3 \nabla u |
\le C\{ |\mathcal{T}| (|\na^2 D_t^2 p| + |\na^2 D_t p| |\na u|) +    \mathcal{T}^2 |\na^3 p| |\na p|   + |\na  \mathcal{T}| | \na D_t^2 p| +    \mathfrak{L}_1 \},
\label{bfordt3u}\end{align}
where
\begin{align}
  & \mathfrak{L}_1 =   |\na D_t p|(|\mathcal{T}\na^2 u|+|\na\mathcal{T}||\na u|)
 +( |\mathcal{T}\na^2 p|  + |\na  \mathcal{T}| |\na p| +|\na u|^2   )(|D_t{\rm div}u| \notag\\
 & + |\mathcal{T}\na^2 p|    + |\na  \mathcal{T}| |\na p| +|\na u|^2 )+  |\mathcal{T}\na p| (| \na D_t  {\rm div} u |   + |\na^2  \mathcal{T}| |\na p|   +|\na^2 u| |\na u|   ) . \label{lot}
\end{align}
\end{lem}
{\em Proof}.  Clearly, \ef{bfordtu} follows from \ef{e25}. It follows from \ef{e25} and \ef{CL-2.23} that
\begin{align*}
 D_t^2 \nabla_b u_a= & -  \mathcal{T} \lt( \na_b \nabla_a D_t p - (\na_a\na_b u^e) \na_e p  \rt)
- (D_t \mathcal{T} )\na_b \nabla_a p \notag\\
 & - (\nabla_b D_t\mathcal{T}) \nabla_a p - (\nabla_b \mathcal{T})  \nabla_a D_t p + D_t \lt( (\nabla_a u_e ) \nabla_b u^e \rt)
\label{}\end{align*}
and
\begin{align*}
& D_t^3 \nabla_b u_a=     -  \mathcal{T}\lt( \na_b \nabla_a D_t^2 p - (\na_a\na_b u^e) \na_e D_t p  \rt) + \mathcal{T}D_t \lt(  (\na_a\na_b u^e) \na_e p  \rt) \\
&\qquad - 2(D_t \mathcal{T}) \lt( \na_b \nabla_a D_t p - (\na_a\na_b u^e) \na_e p  \rt)
- (D_t^2 \mathcal{T} )\na_b \nabla_a p - (\nabla_b D_t^2\mathcal{T}) \nabla_a p \\
&\qquad  -2(\nabla_b D_t\mathcal{T}) \nabla_a D_t p - (\nabla_b \mathcal{T})  \nabla_a D_t^2 p + D_t^2 \lt( (\nabla_a u_e ) \nabla_b u^e \rt),
\label{}\end{align*}
which, together with \ef{CL-2.13}, \ef{r-3.28}, \ef{r-3.29} and \ef{CL-6.17}, imply \ef{bfordt2u} and \ef{bfordt3u}.
Choose $r=2$ in \ef{l4} to get
\begin{align}
D_t \nabla^2 u_a+     \mathcal{T} \nabla^2 \nabla_a p  = -  ( \nabla^2  \mathcal{T} ) \nabla_a p -2   ( \nabla   \mathcal{T} ) \widetilde{\otimes} \na \nabla_a p   +  ( \nabla^2 u^e ) {\rm curl} u_{ae}. \label{q-3.31}
\end{align}
Taking $D_t$ of  \ef{q-3.31} and noticing \ef{CL-2.13}, \ef{e27}, \ef{CL-6.17}, \ef{CL-6.18} and \ef{r-3.28}, we prove \ef{q-3.29}.
\hfill $\Box$

\begin{lem} We have
\begin{align}
& |  D_t^2  \na^2 \mathcal{T} | \le C \{ |\na^2   \mathcal{T}| ( |D_t {\rm div} u|  +  |\na u|^2 + |\na  \mathcal{T}||\na p| )     + |\na  \mathcal{T}|  (|\na D_t {\rm div} u|     + |\na  \mathcal{T}||\na^2 p| \notag\\
&\quad+ |\na^2 u||\na u|) +   \mathcal{T}(|\na^2 D_t {\rm div} u| + |\na  \mathcal{T}||\na^3 p| +|\na u| |\na^3  u|+   |\na {\rm div} u||\na^2 u| ) \}  \label{4/1-1} , \\
& |  D_t^3  \na^2 \mathcal{T} | \le   C ( | \mathcal{T}  \na^2 D_t^2 {\rm div} u|
 +\mathfrak{L}_2 ),  \label{4/1-2}
\end{align}
where
\begin{align}
& \mathfrak{L}_2 =   \mathcal{T} \{  |\na^2  D_t {\rm div} u||\na u|
+
(|\na   D_t {\rm div} u|+|\na    {\rm div} u||\na u|)|\na^2 u|  + |\na^2 {\rm div} u|  (|D_t  {\rm div} u| \notag\\
&  +  |\na u|^2)  + |\na {\rm div} u|  | \mathcal{T} \na^3 p|  \}   +  |\na \mathcal{T}|  \{|\na D_t^2 {\rm div}u| +|\na D_t {\rm div} u||\na u|  +   |\na^2 u|(|\na u|^2 \notag \\
& +  |D_t  {\rm div} u| )  \} +
|\na^2 \mathcal{T}| \{ |D_t^2 {\rm div} u| +  (|D_t  {\rm div} u|+  |\na u|^2) |\na u|  + |\mathcal{T} \na {\rm div} u||\na p| \} \notag \\
 &+ |\na \mathcal{T}| \sum_{r=1}^3 |\na^{3-r} \mathcal{T}|\lt(|\na^{r} D_t p|
+ \sum_{s=1}^r |\na^s p| |\na^{r+1-s} u| \rt) .
    \label{lot3}
\end{align}
\end{lem}
{\em Proof}.   By repeat use of \ef{CL-2.20}, we have for any given function $q$,
\begin{align*}
& D_t^2 \na^2 q = - D_t ( (\na^2 u)\cdot \na q ) -
(\na^2 u)\cdot \na D_t q + \na^2 D_t^2 q, \\
&D_t^3 \na^2 q = \na^2 D_t^3 q - (\na^2 u ) \cdot \na D_t^2 q - D_t \lt((\na^2 u ) \cdot \na D_t q  \rt)- D_t^2 \lt((\na^2 u ) \cdot \na  q\rt),
 \end{align*}
which, together with \ef{CL-2.13}, implies
\begin{align}
& |D_t^2 \na^2 q | \le  |\na^2 D_t^2 q|+ C (|\na^2 u||\na D_t q| +
|D_t \na^2 u||\na q| + |\na^2 u||\na u| |\na q| ),\label{3/31-1} \\
& |D_t^3 \na^2 q | \le  |\na^2 D_t^3 q|+ C \{|\na^2 u||\na D_t^2 q| + |D_t \na^2 u|
(|\na D_t q|+|\na u||\na q|) \notag\\
& \qquad+ |D_t^2 \na^2 u| |\na q| + (|D_t \na u| +|\na u|^2 )|\na^2 u| |\na q| + |\na^2 u| |\na D_t q| |\na u|  \} . \notag
 \end{align}
With this fact, we can use \ef{CL-6.17}, \ef{r-3.28}-\ef{rr-3.30}  and \ef{q-3.29} to show \ef{4/1-1} and \ef{4/1-2}.

\begin{lem} For $r\ge 0$, we have
\begin{align}
&|\na^r D_t \Delta u| \le  | \na^{r+1} D_t  {\rm div} u | + C \sum_{s=0}^{r+1}
\lt(|\na^{1+s}u||\na^{r+2-s} u| +|\na^{1+s} \mathcal{T} | |\na^{r+2-s} p|  \rt)  \label{dtdau},\\
&|D_t^2 \Delta u| \le  | \na D_t^2  {\rm div} u | +  C \sum_{r=1}^3 |\na^{3-r} \mathcal{T}|  \sum_{s=1}^r |\na^s p| |\na^{r+1-s} u|   +  C \sum_{r=1}^2 |\na^{3-r} \mathcal{T}| |\na^{r} D_t p|\notag\\
& \qquad \qquad
+ C |\na^2 u||\na u|^2,
 \label{dt2dau}\\
&|D_t^3 \Delta u| \le  | \na D_t^3  {\rm div} u | + C (  |\na u|| \mathcal{T}\na^3 D_t p| + |\na \mathcal{T}||\na^2 D_t^2 p|+ |\na^2 \mathcal{T}||\na D_t^2 p| +    \mathfrak{L}_3  ),  \label{dt3dau}
\end{align}
where
\begin{align}
  \mathfrak{L}_3 = & ( | \mathcal{T}\na^2 p| +  |\na \mathcal{T}| |\na  p| +  |\na u|^2)
( |\na \mathcal{T}||\na^2 p| + |\na^2 \mathcal{T}| |\na p|
+|\na^2 u||\na u|)\notag \\
&+|\na^2 p|   ( |\mathcal{T} \na D_t {\rm div} u|  + |\na \mathcal{T}||D_t {\rm div} u|   ) +|\na p| \{|\mathcal{T}\na^2 D_t {\rm div} u|\notag \\
& +|\na \mathcal{T}|(|\mathcal{T} \na^3 p| + |\na D_t{\rm div} u| ) + |\na^2 \mathcal{T}||D_t{\rm div} u|\} + |\na D_t p| |\mathcal{T} \na^2 {\rm div} u|
\notag\\
&
+ |\na^2 u|  \sum_{r=1}^2 |\na^{2-r} \mathcal{T}|\lt(|\na^{r} D_t p|
+ \sum_{s=1}^r |\na^s p| |\na^{r+1-s} u| \rt)    \notag\\
&+|\na  u| \lt\{ \sum_{r=1}^3 |\na^{3-r} \mathcal{T}|\lt(|\na^{r} D_t p|
+ \sum_{s=1}^r |\na^s p| |\na^{r+1-s} u| \rt) -   | \mathcal{T}\na^3 D_t p|\rt\}  .
  \label{lot2}
\end{align}
\end{lem}
{\em Proof}.
Take $D_t^r$ $(r=1,2,3)$ of \ef{e27} and use \ef{r-3.28} and \ef{r-3.29} to get
\begin{align}
& |D_t {\rm curl} u| \le C |\na \mathcal{T}||\na p|, \label{3/31-2}\\
& |D_t^2 {\rm curl} u| \le C |\na \mathcal{T}| (|\na D_t p| + |\na p|  |\na u| )+ C \mathcal{T} |\na p| |\na {\rm div} u|. \label{3/31-3}
\end{align}
 It follows from \ef{e27} and \ef{CL-2.20} that
\begin{align}
D_t \na_c {\rm curl} u_{ab} = &(\na_c\na_b \mathcal{T})\na_a p + (\na_b \mathcal{T}) \na_c\na_a p - (\na_c \na_a \mathcal{T}) \na_b p -(\na_a \mathcal{T})\na_c \na_b p \notag\\
& -(\na_a\na_c u^e) {\rm curl} u_{eb} -(\na_b\na_c u^e) {\rm curl} u_{ae}, \label{3/31-5}
\end{align}
which, together with \ef{CL-2.13}, \ef{r-3.28}, \ef{CL-6.17}, \ef{CL-6.18} and \ef{3/31-2}, implies
\begin{align}
& |D_t^2 \na {\rm curl} u| \le C|\mathcal{T}| (  |\na^3 p||\na u| + |\na^2 p||\na^2 u| + |\na p||\na^2 {\rm div} u|)
  +C|\na \mathcal{T}| (|\na^2 D_t p| \notag \\
&  \qquad  +|\na^2 p||\na u| + |\na p ||\na^2 u|)
 +C|\na^2 \mathcal{T}| (|\na D_t p| +|\na p||\na u|) + C |\na^2 u||\na u|^2 . \label{3/31-4}
\end{align}
Take $D_t^2$ of \ef{3/31-5} and use \ef{CL-2.13}, \ef{4/1-1}, \ef{3/31-1}, \ef{r-3.28}, \ef{r-3.29}, \ef{CL-6.17}, \ef{CL-6.18}, \ef{q-3.29}, \ef{3/31-2} and \ef{3/31-3} to obtain
\begin{align}
& |D_t^3 \na {\rm curl} u| \le  C \lt(|\mathcal{T}||\na u||\na^3 D_t p| + |\na \mathcal{T}||\na^2 D_t^2 p| +|\na^2 \mathcal{T}||\na D_t^2 p| +  \mathfrak{L}_3  \rt)  ,\label{3/31-6}
\end{align}
where $\mathfrak{L}_3$ is given by \ef{lot2}.
It follows from the definition of $\Delta$ that
\begin{align*}
\Delta u_a= g^{ce} \na_c\na_e u_a = \na_a {\rm div} u + g^{ce} \na_c (\na_e u_a - \na_a u_e)
=\na_a {\rm div} u + g^{ce} \na_c {\rm curl} u_{ea},
\end{align*}
which, together with \ef{CL-2.13}, implies
\begin{align*}
& D_t \Delta u_a
=\na_a D_t {\rm div} u + (D_t g^{ce}) \na_c {\rm curl} u_{ea} + g^{ce} D_t \na_c {\rm curl} u_{ea} \\
&|D_t^2 \Delta u| \le  | \na D_t^2  {\rm div} u | + C    ( |D_t \na u| +  |\na u|^2 )|\na^2 u| \notag\\
&  \qquad +C (|\na u||D_t \na {\rm curl} u| + |D_t^2 \na {\rm curl} u| ) , \\
& |D_t^3 \Delta u| \le   | \na D_t^3  {\rm div} u |
+ C ( |D_t^2 \na u| + |D_t \na u||\na u|  +|\na u|^3  )|\na^2 u| \\
& \qquad + C ( ( |D_t \na u| +  |\na u|^2 )|D_t \na {\rm curl} u| + |\na u||D_t^2 \na {\rm curl} u| + |D_t^3 \na {\rm curl} u| )
.
\end{align*}
With this, the lemma can be proved by noting \ef{CL-6.17}, \ef{bfordt2u}, \ef{3/31-4} and \ef{3/31-6}. \hfill $\Box$

\begin{lem} We have for $r\ge 0$,
\be\label{4/10-5}\begin{split}
  | \mathcal{T}   \Delta D_t  p +(\nabla  \mathcal{T})\cdot \nabla D_t  p |   \le & |D_t^2 {\rm div} u| + C\{ |\mathcal{T}| ( |\na^2 p| |\na u| + |\na p||\na^2 u|)  \\
&+      |\na \mathcal{T}| |\na p| |\na u|
    +   |\na u|^3 \} ,
\end{split}\ee
\be\label{nadadtp}\begin{split}
&  | \mathcal{T} \na^r  \Delta D_t p | \le    | \na^r D_t^2 {\rm div}  u |
+ C  \sum_{s_1=0}^{r+1} |\na^{s_1} \mathcal{T}| \sum_{s_2=0}^{r+1-s_1}  |\na^{s_2+1} u|
  |\na^{r+2-s_1-s_2} p|   \\
  & \qquad  + C \sum_{s=0}^r |\na^{s+1} \mathcal{T}| |\na^{r+1-s} D_t p|
+ C \sum_{s_1+s_2+s_3=r} |\na^{1+s_1} u| |\na^{1+s_2} u||\na^{1+s_3} u|,
\end{split}\ee
\be\label{dadt2p}\begin{split}
  | \mathcal{T}   \Delta D_t^2 p +(\nabla  \mathcal{T})\cdot \nabla D_t^2 p |
\le     |D_t^3{\rm div} u| + C (   |\na u| |\mathcal{T} \na^2 D_t p|   +  \mathfrak{L}_1 ),
\end{split}\ee
\be\label{174/5}\begin{split}
  | \mathcal{T} \na^r  \Delta D_t^2 p | \le    | \na^r D_t^3 {\rm div}  u |
+  C\sum_{i=1}^{r+1} |\na^i \mathcal{T}||\na^{r+2-i}D_t^2 p| + C \mathfrak{L}_4,
\end{split}\ee
\be\label{174/4}\begin{split}
| \mathcal{T}   \Delta D_t^3 p +(\nabla  \mathcal{T})\cdot \nabla D_t^3 p |
\le     |D_t^4{\rm div} u| + C   ( |\na u|  | \mathcal{T}\na^2 D_t^2 p|   + |\na u| \mathfrak{L}_1  + \mathfrak{L}_5 ),
\end{split}\ee
where $\mathfrak{L}_1$ is defined by \ef{lot},
\bee\begin{split}
   \mathfrak{L}_4 =   |\na^{r+1}(\mathcal{T}(\na u)\na D_t p)|
 +  |\na^{r+1}(\mathcal{T}(\na p) D_t {\rm div} u)| +|\na^{r+1}(\mathcal{T}(\na p)  (\na \mathcal{T}) \na p)| \\
 +|\na^{r+1}(\mathcal{T}(\na p)  (\na u) \na u)| +|\na^{r }(\mathcal{T}^2 (\na^2 p) \na^2 p  )|+ |\na^{r }( (\na u) (\na u) (\na u) \na u )|,
\end{split}\eee
and
 \bee\begin{split}
  & \mathfrak{L}_5 =   |\na D_t^2 p|(|\mathcal{T}\na^2 u|+|\na\mathcal{T}||\na u|)
 +( |\mathcal{T}\na^2 D_t p|  + |\na  \mathcal{T}| |\na D_t p| +|\na^2 u|  | \mathcal{T}\na  p|  )(|D_t{\rm div}u|  \\
 & + |\mathcal{T}\na^2 p|    + |\na  \mathcal{T}| |\na p| +|\na u|^2 )+  |\mathcal{T}\na D_t  p| (| \na D_t  {\rm div} u |   + |\na^2  \mathcal{T}| |\na p|      ) +   |\mathcal{T} \na p| |\na D_t^2 {\rm div} u|
 \\
& +  ( |\mathcal{T} \na^2 p|  + |\na \mathcal{T} ||\na p|)| D_t^2 {\rm div} u|  + \mathcal{T}| \na p|^2 |\na^2 {\rm div} u| + | \mathcal{T} \na p|| \na u|| \mathcal{T} \na^3 p|.
\end{split}\eee

\end{lem}
{\bf Proof}.  In view of \ef{CL-2.22} and  \ef{lappi}, we see that
\begin{align}
&  \Delta D_t p = D_t \Delta p + 2 (\nabla_e \nabla  p)\cdot \nabla  u^e   + (\Delta u ) \cdot \nabla p,   \label{r-3.31}\\
&  \mathcal{T} D_t \Delta p    =  - D_t^2{\rm div} u  - (D_t\mathcal{T})\Delta p    -   D_t \lt( (\nabla  \mathcal{T})\cdot \nabla  p + (\nabla_e u )\cdot\nabla u^e \rt), \notag
\end{align}
which means
 \begin{align}
 \mathcal{T} \Delta D_t p = &  - D_t^2{\rm div} u  - (D_t\mathcal{T})\Delta p
  -  \{(\nabla  D_t \mathcal{T})\cdot \nabla  p + (\nabla  \mathcal{T})\cdot \nabla D_t p   \notag \\
  &  +(D_t g^{ab}) (\na_a \mathcal{T} )\na_b p \}
    - \{g^{ab}g^{ef}( (D_t \na_e u_a)\na_b u_f +  ( \na_e u_a) D_t\na_b u_f  )
   \notag \\
 &    + ( D_t(g^{ab}g^{ef}) ) ( \na_e u_a) \na_b u_f  \}
+ 2  \mathcal{T}(\nabla_e \nabla  p)\cdot \nabla  u^e   +   \mathcal{T} (\Delta u ) \cdot \nabla p. \label{r-3.30}
\end{align}
Clearly, \ef{4/10-5} can be derived directly from \ef{r-3.30}. \ef{nadadtp} can be proved by taking  $\na^r$ of \ef{r-3.30} and using \ef{bfordtu} and \ef{r-3.28}.

It follows from \ef{CL-2.22}  that
\begin{align}
 & \mathcal{T} \Delta D_t^2 p  = \mathcal{T} D_t \Delta D_t p + 2 \mathcal{T}(\nabla_e \nabla  D_t p)\cdot \nabla  u^e   + \mathcal{T} (\Delta u ) \cdot \nabla D_t p \notag\\
& = D_t ( \mathcal{T}  \Delta D_t p ) -  (D_t \mathcal{T} ) \Delta D_t p  + 2 \mathcal{T}(\nabla_e \nabla  D_t p)\cdot \nabla  u^e   + \mathcal{T} (\Delta u ) \cdot \nabla D_t p, \label{eh3}
\end{align}
which, together with \ef{r-3.27}, implies
\begin{align}
 | \mathcal{T} \Delta D_t^2 p + (\nabla  \mathcal{T})\cdot \nabla D_t^2 p | \le & | D_t ( \mathcal{T}  \Delta D_t p ) + (\nabla  \mathcal{T})\cdot \nabla D_t^2 p|  \notag\\
 & + C  |\mathcal{T}| ( |\nabla^2  D_t p| | \nabla  u| +   |\nabla   D_t p| | \nabla^2  u| ). \label{r-3.54}
\end{align}
Take $D_t$ of \ef{r-3.30} and use \ef{CL-2.13},  \ef{CL-6.18}, \ef{r-3.28}, \ef{r-3.29}, \ef{bfordtu}, \ef{bfordt2u}, \ef{dtdau}  and  \ef{r-3.31} to get
\begin{align*}
| D_t ( \mathcal{T}  \Delta D_t p ) +(\nabla  \mathcal{T})\cdot \nabla D_t^2 p| \le |D_t^3{\rm div} u| + C (  |  \mathcal{T}||\na u| |\na^2 D_t p| +  \mathfrak{L}_1  ),
\end{align*}
where $ \mathfrak{L}_1$ is given by \ef{lot}. With \ef{r-3.54}, we then prove \ef{dadt2p}. \ef{174/5} can be proved  in  a similar way to deriving \ef{nadadtp}.

\ef{174/4} can be shown in a similar way to deriving \ef{dadt2p}.
\hfill $\Box$

\begin{lem} For $r\ge 0$ and $i=0,1,2$, we have
\begin{align}
  |\na^r D_t  {\rm div } u-\mathcal{T} \na^r \Delta {\rm div } u | \le &
   C \sum_{s=0}^{r+1}|\na^{s+1} \mathcal{T}||\na^{r-s+2}u|  + C \sum_{s=0}^r |\na^{s } {\rm div} u| |\na^{r-s } {\rm div} u|   \label{nadadivu}
\end{align}
and
\begin{align}
 |\na^i D_t^2 {\rm div} u - & \mathcal{T} \na^i \Delta D_t {\rm div} u|
\le  C \sum_{r=-1}^{i-1} |\mathcal{T} \na^{r+3} p| |\na^{i+1-r} \mathcal{T} | \notag\\ &+C\sum_{r=1}^{i+2} \lt( |\mathcal{T} \na^r {\rm div} u| +  \sum_{s=0}^{r-1}|\na^{s+1}u||\na^{r-s}\mathcal{T}|\rt) |\na^{i+3-r} u|
\notag\\
&+C\sum_{r=0}^i \lt( |\na^r D_t {\rm div} u| |\na^{i+1-r} u|+  |\na^{r+1} D_t {\rm div} u| |\na^{i+1-r} \mathcal{T} |  \rt) \notag\\
&+ C\sum_{r=-1}^i \sum_{s=0}^{r+1}
\lt(|\na^{1+s}u||\na^{r+2-s} u| +|\na^{1+s} \mathcal{T} | |\na^{r+2-s} p|  \rt)   |\na^{i-r+1}  \mathcal{T} |.
 \label{dt2divu}
\end{align}
\end{lem}
{\em Proof}. Clearly, \ef{nadadivu} follows directly from \ef{e9}. Take $D_t$ of \ef{e9} and use \ef{CL-2.22} to get
\begin{align}
&D_t^2 {\rm div} u -  \mathcal{T} \Delta D_t {\rm div} u
=  ( D_t \mathcal{T} ) \Delta {\rm div} u - \mathcal{T} \{ (\Delta u  )  \cdot \nabla  {\rm div} u + 2   (\nabla_e \nabla {\rm div} u ) \cdot \nabla  u^e \} \notag \\
& \qquad -D_t\{   (\Delta u  )  \cdot \nabla  \mathcal{T} + 2   (\nabla_e \nabla  \mathcal{T} ) \cdot \nabla  u^e  -    ({\rm div } u) ^2 - 2   (  \nabla \mathcal{T} ) \cdot \nabla  {\rm div } u\} . \label{4/2-1}
\end{align}
Taking $\na^i$ $(i=1,2)$ of \ef{4/2-1} and noticing \ef{CL-2.13}, \ef{r-3.28}, \ef{CL-6.17}, \ef{bfordtu} and \ef{dtdau}, we can obtain \ef{dt2divu}.  \hfill$\Box$

\begin{lem} We have
\begin{align}
&|D_t^3 {\rm div} u -  \mathcal{T} \Delta D_t^2 {\rm div} u|
\le  C(|\na \mathcal{T}  | |\na D_t^2 {\rm div} u| + |\na u||  D_t^2  {\rm div} u|
+ \mathfrak{L}_6), \label{dt3divu}
\\
&|D_t^4 {\rm div} u -  \mathcal{T} \Delta D_t^3 {\rm div} u | \le C(|\na \mathcal{T}  | |\na D_t^3 {\rm div} u| + |\na u||  D_t^3  {\rm div} u|
+ \mathfrak{L}_{71}+\mathfrak{L}_{72}), \label{dt4divu}
\end{align}
where
\begin{align}
 &\mathfrak{L}_6 = |\na  \mathcal{T} | \lt \{   \sum_{r=1}^3 |\na^{3-r} \mathcal{T}|  \sum_{s=1}^r |\na^s p| |\na^{r+1-s} u|   + \sum_{r=1}^2 |\na^{3-r} \mathcal{T}| |\na^{r} D_t p|\rt\} \notag \\
&\quad+  |\na^2  \mathcal{T} | \lt\{\sum_{r=1}^2 |\na^{2-r} \mathcal{T}|\lt(|\na^{r} D_t p|
+ \sum_{s=1}^r |\na^s p| |\na^{r+1-s} u| \rt)  +    |\na u|^3  \rt\}\notag\\
&\quad +|\na u| \lt\{\sum_{s=0}^2 |\na^s \mathcal{T}||\na^{2-s}D_t {\rm div} u| + |\mathcal{T}\na^3 u| |\na u|\rt\} + |D_t{\rm div} u|^2 \notag\\
&\quad +  |\na^2 u| \{|\mathcal{T}\na D_t {\rm div} u| + |\na u|
|\mathcal{T}\na  {\rm div} u|  + |\na \mathcal{T}|(|D_t {\rm div} u | + |\na u|^2) \} \notag\\
&\quad  +  |\mathcal{T} \na^2 {\rm div} u|  (|\mathcal{T}\na^2 p| + |\na \mathcal{T}||\na p|+ |\na u|^2 + |D_t{\rm div} u|) ,
\label{lot4}
\end{align}
\begin{align}
& \mathfrak{L}_{71} = |D_t^3 \na u||\na^2  \mathcal{T}| +  |D_t^3 \Delta u||\na \mathcal{T}| + |D_t^3 \na^2   \mathcal{T}| |\na u|
 +|D_t^2 \Delta u|(|\na u||\na  \mathcal{T}|
 \notag\\
& \quad+ | \mathcal{T}\na {\rm div} u|)
  +|D_t^2 \na u| (| \mathcal{T}\na^2 {\rm div} u| + |\na^2  \mathcal{T}||\na u| + |\na  \mathcal{T}||\na^2 u|)
 \notag\\
& \quad + |D_t \na^2 u||\na u|| \mathcal{T} \na {\rm div} u|
   + |D_t \Delta u| (|D_t\na u||\na  \mathcal{T}| + |\na D_t^2  \mathcal{T}|
  + | \mathcal{T} \na D_t {\rm div} u|
    \notag\\
& \quad+ |\na  \mathcal{T}||\na u|^2 + | \mathcal{T} \na {\rm div} u| |\na u|)
  + (|D_t \na u|+|\na u|^2) \{| \mathcal{T} \na^2 D_t {\rm div} u |
   \notag\\
& \quad+ |D_t^2 \na^2  \mathcal{T}| + |D_t \na^2  \mathcal{T}||\na u| + (|D_t \na u|+|\na u|^2)|\na^2  \mathcal{T}|
  + | \mathcal{T} \na^2 {\rm div} u||\na u|
    \notag\\
& \quad+ |\na^2 u| (| \mathcal{T} \na {\rm div} u|
 + |\na  \mathcal{T}||\na u|) + |\na D_t {\rm div} u||\na  \mathcal{T}|
  \}\label{lot51}
\end{align}
and
\begin{align}
& \mathfrak{L}_{72} = |\na D_t{\rm div} u|| \na D_t^2   \mathcal{T}| + |\mathcal{T} \na^2 D_t^2 {\rm div} u||\na u|
 + |\na D_t^2{\rm div} u| (|\mathcal{T}\na {\rm div} u|+|\na \mathcal{T}||\na u| \notag\\
& \quad+  |\mathcal{T}\Delta u |)
   +( |\na D_t^3 \mathcal{T}| +  |\na D_t^2  \mathcal{T}| |\na u|)|\na {\rm div} u|
 + |D_t^2 {\rm div} u| |D_t {\rm div} u|   \notag\\
&\quad + |D_t^3 \mathcal{T}||\na^2 {\rm div} u|
  +  (|\na^2 u||\na {\rm div} u|
  +|\na  u||\na^2 {\rm div} u| + |\na^2 D_t{\rm div} u|) |\mathcal{T}D_t{\rm div} u|
  \notag\\
&\quad  + |\na^2 u| (|\na u| |\na \mathcal{T}||D_t{\rm div} u| + |\na u| |\mathcal{T} D_t{\rm div} u|)  + |\Delta u || \nabla   D_t^3 \mathcal{T} |.  \label{lot52}
\end{align}

\end{lem}
{\em Proof}. Take $D_t^r$ $(r=1,2)$ of \ef{4/2-1} and use \ef{CL-2.22} to get
\begin{align}
&D_t^3 {\rm div} u -  \mathcal{T} \Delta D_t^2 {\rm div} u  \notag \\
= & ( D_t \mathcal{T} ) \Delta D_t {\rm div} u - \mathcal{T} ( (\Delta u  )  \cdot \nabla  D_t {\rm div} u + 2   (\nabla_e \nabla D_t {\rm div} u ) \cdot \nabla  u^e ) \notag\\
& + D_t \{ ( D_t \mathcal{T} ) \Delta {\rm div} u - \mathcal{T} ( (\Delta u  )  \cdot \nabla  {\rm div} u + 2   (\nabla_e \nabla {\rm div} u ) \cdot \nabla  u^e ) \} \notag \\
&   -D_t^2\{   (\Delta u  )  \cdot \nabla  \mathcal{T} + 2   (\nabla_e \nabla  \mathcal{T} ) \cdot \nabla  u^e  -    ({\rm div } u) ^2 - 2   (  \nabla \mathcal{T} ) \cdot \nabla  {\rm div } u\}\label{eh1}
\end{align}
and
\begin{align}
&D_t^4 {\rm div} u -  \mathcal{T} \Delta D_t^3 {\rm div} u  \notag \\
= &  ( D_t \mathcal{T} ) \Delta D_t^2 {\rm div} u - \mathcal{T} ( (\Delta u  )  \cdot \nabla   D_t^2 {\rm div} u + 2   (\nabla_e \nabla D_t^2 {\rm div} u ) \cdot \nabla  u^e )
\notag\\
&+D_t \{( D_t \mathcal{T} ) \Delta D_t {\rm div} u - \mathcal{T} ( (\Delta u  )  \cdot \nabla  D_t {\rm div} u + 2   (\nabla_e \nabla D_t {\rm div} u ) \cdot \nabla  u^e ) \} \notag\\
& + D_t^2 \{ ( D_t \mathcal{T} ) \Delta {\rm div} u - \mathcal{T} ( (\Delta u  )  \cdot \nabla  {\rm div} u + 2   (\nabla_e \nabla {\rm div} u ) \cdot \nabla  u^e ) \} \notag \\
&   -D_t^3\{   (\Delta u  )  \cdot \nabla  \mathcal{T} + 2   (\nabla_e \nabla  \mathcal{T} ) \cdot \nabla  u^e  -    ({\rm div } u) ^2 - 2   (  \nabla \mathcal{T} ) \cdot \nabla  {\rm div } u\}.\label{eh2}
\end{align}
With these two equations, we can obtain \ef{dt3divu} and \ef{dt4divu} by simple calculations and noticing \ef{CL-2.13}, \ef{CL-2.22}, \ef{CL-6.18}, \ef{3/31-1}, \ef{r-3.28}, \ef{r-3.29}, \ef{bfordtu}, \ef{bfordt2u}, \ef{4/1-1}, \ef{dtdau} and \ef{dt2dau}.
\hfill$\Box$

\section{Proof of Theorem \ref{mainthm}}
Let
\begin{subequations}\label{17energy}\begin{align}
& E_0 (t) =   \int_\Omega \mathcal{T}^{-1}  |u|^2 d\mu_g  \label{17e0}  \\
&  E_r (t) =   \int_\Omega  \mathcal{T}^{-1}  g^{ab} \za^{cd} Q(\na^{r-1}\na_c u_a, \na^{r-1}\na_d u_b)  d\mu_g   + \int_\Omega|\na^{r-1} {\rm curl}  u|^2    d\mu_g
 \notag\\
& +   \int_\Omega  |\na D_t^{r-1}{\rm div} u|^2  d\mu_g+   \int_{\pl\Omega}  \za^{cd} Q(\na^{r-1} \na_c p, \na^{r-1} \na_d p) (-\na_N p)^{-1} d \mu_\za , \ \ r\ge 1 , \label{17er}
\end{align}\end{subequations}
where $ Q(\al, \ba)=\varsigma^{IJ} \al_{I} \ba_{J}$.
Suppose that the following {\em a priori} assumptions are true:
\begin{subequations}\label{17/5/19}\begin{align}
& \underline{V}  \le  {\rm Vol}\Omega (t)     \le   \overline{V}  \  \  & {\rm on} \ \   &[0,T] ,  \label{vol}\\
&|\ta|+ 1/\iota_0 \le K \ \  & {\rm on} \ \  &[0,T]\times \pl \Omega ,\label{CL-7.3}\\
&-\na_N p \ge \ea_b >0   \ \  & {\rm on} \ \ &[0,T]\times \pl \Omega , \label{CL-7.4}\\
&  \sum_{i=1}^{n-1}   (|\na_N  D_t^i p| +|\na_N  D_t^i {\rm div} u|) +|\na^2 p|  \le L   \ \ &   {\rm on} \ \ &[0,T]\times \pl \Omega ,\label{CL-7.5}\\
&   |\na  p | + |\na u|  + |\na  \mathcal{T} | + |\na  {\rm div} u|   \le M      \ \  &  {\rm in} \ \ & [0,T]\times   {\Omega},\label{CL-7.6} \\
&  |D_t  p| +  | D_t  {\rm div} u |  + |\na^2 \mathcal{T}|  \le \widetilde{M}    \ \  &  {\rm in} \ \ & [0,T]\times   \Omega ,\label{CL-7.6'}
\end{align}\end{subequations}
where ${\rm Vol}\Omega (t)=\int_\Omega d\mu_g$.
Let $\iota_0$ and $\iota_1$ be as in Definitions \ref{defn1.2} and \ref{defn2.3}. Then we have, due to Lemma \ref{lemkk1} and \ef{CL-7.3}, that
$$\iota_1^{-1} \le \max\{\ea_1^{-1}|\|\ta\||_{L^\iy}, \  (2\iota_0)^{-1}\} \le \ea_1^{-1} K,$$
which means
\begin{align}
\iota_1^{-1}(t) \le K_1 \ \  {\rm on} \ \ [0,T], \ \ {\rm where} \ \ K_1=\ea_1^{-1} K. \label{w4.10}
\end{align}
Before stating the result, let us  notice  the boundary conditions and maximum principle, which are due to \ef{s110} and  \ef{e7}, as follows:
\begin{align}
& p=0, \ \  \mathcal{T}=  \mathcal{T}_b, \ \    D_t   \mathcal{T} =0    \ \ {\rm and} \ \   {\rm div} u=0 \ \ {\rm on} \ \    [0,T]  \times \pl \Omega \label{bdry}, \\
& \underline{\mathcal{T}} \le \mathcal{T} \le \overline{\mathcal{T}}
 \ \     {\rm in} \ \   [0,T]\times   \overline{\Omega}, \label{bdta}
\end{align}
where
$$\underline{\mathcal{T}}= \min\lt\{\min_{y\in \Omega} \mathcal{T}(0,y), \ \mathcal{T}_b \rt\}   \ \  {\rm and} \ \ \overline{\mathcal{T}}= \max\lt\{\max_{y\in \Omega} \mathcal{T}(0,y), \ \mathcal{T}_b \rt\}.$$

\begin{prop}\label{prop1}
Let $n=2,3$ and $1\le r\le n+2$. Then there are continuous functions
$\mathcal{F}_r$ with $\mathcal{F}_r\big|_{t=0}=1$, such that for any smooth solutions of \ef{l2}-\ef{s110} for $0\le t\le T$ satisfying \ef{17/5/19}, we have
\begin{subequations}\label{17prop1}\begin{align}
&E_0(t)\le E_0(0)+ \lt[\mathcal{F}_1\lt(t,\overline{V},K,  \ea_b^{-1},L, M, \widetilde{M}, \underline{\mathcal{T}}^{-1}, \overline{\mathcal{T}}  \rt)-1\rt] \sum_{s=1}^2 E_s(0),\label{1752-0}\\
&\sum_{s=1}^2 E_s(t)\le  \mathcal{F}_2\lt(t,\overline{V},K, \ea_b^{-1},L, M, \widetilde{M}, \underline{\mathcal{T}}^{-1}, \overline{\mathcal{T}}   \rt)  \sum_{s=1}^2 E_s(0),\label{1752-2}\\
&\sum_{s=1}^rE_s(t)\le   \mathcal{F}_r\lt(t,\overline{V},K,  \ea_b^{-1},L, M, \widetilde{M},\underline{\mathcal{T}}^{-1}, \overline{\mathcal{T}} ,  E_1(0),\cdots, E_{r-1}(0)  \rt)  \sum_{s=1}^r E_s(0), \ r\ge 3.\label{1752-r}
\end{align}\end{subequations}
\end{prop}

\begin{prop}\label{prop2} Let ${\rm Vol}\mathscr{D}_0 , K_0,  \underline{\ea}_0$ and $M_0$ be defined by \ef{17initial}, and $n=2,3$. Then there are continuous functions
$\mathscr{T}_n$   such that if
\be\label{1754}
T\le \mathscr{T}_n\lt({\rm Vol}\mathscr{D}_0 , K_0,  \underline{\ea}_0^{-1},  \underline{\mathcal{T}}^{-1}, \overline{\mathcal{T}} ,M_0,  E_0(0),\cdots,E_{n+2}(0) \rt),
\ee
then any smooth solutions of the free surface problem \ef{l2}-\ef{s110} for $0\le t\le T$ satisfies
\be\label{17/5/27-1}
 \sum_{s=0}^{n+2}E_s(t)\le   2 \sum_{s=0}^{n+2} E_s(0), \ \  0\le t\le T,
\ee
\be\label{17/5/27-2}
2^{-1} {\rm Vol}\mathscr{D}_0     \le  {\rm Vol}\Omega (t)     \le  2 {\rm Vol}\mathscr{D}_0, \ \  0\le t\le T,
\ee
\be\label{17/5/27-3}
|\|\ta(t,\cdot)\||_{L^\iy}+ \iota_0^{-1}(t)\le 18 K_0 , \ \  0\le t\le T,
\ee
\be\label{17/5/27-4}
 -\na_N p(t,  y) \ge 2^{-1}\underline{\ea}_0        \ \ {\rm for}  \ \    y \in  \pl \Omega, \ \  0\le t\le T,
\ee
\be\label{17/5/27-5}\begin{split}
    \|\na p(t,\cdot)\|_{L^\iy} + \|\na u(t,\cdot)\|_{L^\iy} + \|\na \mathcal{T}(t,\cdot)\|_{L^\iy}  \le 2 M_0, \ \  0\le t\le T,
\end{split}\ee
\be\label{17/5/27-6}\begin{split}
&\sum_{i=1}^{n-1} \lt(|\|\na  D_t^i p(t,\cdot)\||_{L^\iy} + |\|\na D_t^i {\rm div} u(t,\cdot)\||_{L^\iy}\rt)+ |\|\na^2 p  (t,\cdot)\||_{L^\iy}
 \\
& +  \|\na {\rm div} u(t,\cdot)\|_{L^\iy} + \|D_t p (t,\cdot)\|_{L^\iy}
 + \|D_t {\rm div} u (t,\cdot)\|_{L^\iy}
 + \|\na^2 \mathcal{T} (t,\cdot)\|_{L^\iy}  \\
& \le    C \lt({\rm Vol}\mathscr{D}_0 , K_0, \underline{\ea}_0^{-1}, \underline{\mathcal{T}}^{-1}, \overline{\mathcal{T}} ,  M_0, E_0(0),\cdots,E_{n+2}(0) \rt),\ \  0\le t\le T.
\end{split}\ee
\end{prop}
Clearly, Theorem \ref{mainthm} is a conclusion of this proposition. (Indeed, the compatibility condition $\mathcal{T}(0,y)=\mathcal{T}_b$ on $\pl\Oa$ implies that $\underline{\mathcal{T}}= \min_{y\in \Omega} \mathcal{T}(0,y) $  and $\overline{\mathcal{T}}=  \max_{y\in \Omega} \mathcal{T}(0,y) $.)

\subsection{Energy estimates}
 In the proof we make use of a  fact, which follows from \ef{CL-2.13}, that for a function $f=f(t,y)$,
\begin{align}
\frac{d}{dt}\int_{\Omega} f d\mu_g = \int_{\Omega} \lt(D_t f + f {\rm div} u \rt)   d\mu_g. \label{s3.7}
\end{align}

First, we deal with  the temporal derivatives of $\int_\Omega |\na  D_t^{r} {\rm div} u|^2  d\mu_g$, which can be bounded as follows:
\be\label{17/9/4}\begin{split}
&  \frac{d}{dt}\int_\Omega |\na  D_t^{r} {\rm div} u|^2  d\mu_g \le C_r(\cdots )     \lt( \sum_{i=0}^{1} \| \na^i D_t^{r } {\rm div} u\|^2
       +\sum_{i=0}^{2} \sum_{j=1}^{r-1} \| \na^i D_t^j {\rm div} u\|^2  \rt. \\
 & \lt.  +\sum_{i=1}^{2} \| \na^i D_t^{r-1} p\|^2 +  \sum_{i=1}^{3} \sum_{j=1}^{r-2}  \| \na^i D_t^j p\|^2 + \sum_{i=1}^{r+1}\lt(\|\na^i u\|^2   + \|\na^i p\|^2  \rt)+\sigma(r) \|\nabla^{r+1} {\rm div u}\|^2   \rt)
\end{split}\ee
for $r=1,2,3,4$, where  $\sigma(r)=1$ for $r=1,2$,  and $\sigma(r)=0$ for $r=3,4$. However, various quantities  that  the constants $C_r$ in \ef{17/9/4} depend are
quite different for different values of  $r$.  Identifying clearly the quantities  that   the constants $C_r$ depend will be important to closing the arguments.

\begin{lem}\label{17lemdiv} We have
\be\label{17div1}
  \frac{d}{dt}\int_\Omega |\na  {\rm div} u|^2  d\mu_g
\le C   \lt(\| \mathcal{T}^{-1}\|_{L^\iy} ,\|\na  u\|_{L^\iy},  \|\na \mathcal{T} \|_{L^\iy}\rt) \lt(  \|\na^2 \mathcal{T} \|^2 +\|\na^2 u\|^2 + \|\na u\|^2   \rt) ,
\ee
\begin{equation}\label{17div2}\begin{split}
&  \frac{d}{dt}\int_\Omega |\na D_t {\rm div} u|^2  d\mu_g \le  C  (\| \mathcal{T}^{-1}\|_{L^\iy},  \| \mathcal{T} \|_{L^\iy},\|\na  u\|_{L^\iy},  \|\na \mathcal{T} \|_{L^\iy}, \|\na p \|_{L^\iy},\|\na^2 \mathcal{T} \|_{L^\iy},\\
&\qquad  \|\na {\rm div} u \|_{L^\iy})
 \lt(\sum_{i=0}^1 \| \na^i D_t{\rm div} u\|^2   + \|\na^2 {\rm div} u \|^2+  \sum_{i=1}^{2}\lt(\|\na^i u\|^2    + \|\na^i p\|^2  \rt)\rt),
\end{split}\end{equation}
\begin{equation}\label{17div3}\begin{split}
&  \frac{d}{dt}\int_\Omega |\na  D_t^2{\rm div} u|^2  d\mu_g \le C( K_1, \| \mathcal{T}^{-1}\|_{L^\iy},\|  \mathcal{T} \|_{L^\iy}, \|\na  u\|_{L^\iy} , \|\na  \mathcal{T} \|_{L^\iy},     \|\na  p \|_{L^\iy},  \\
&\qquad   \|\na^2 \mathcal{T} \|_{L^\iy},\|\na {\rm div} u \|_{L^\iy}, \| D_t {\rm div} u \|_{L^\iy})\lt(\sum_{i=0}^1 \| \na^i D_t^2 {\rm div} u\|^2 +  \sum_{i=0}^2 \| \na^i D_t{\rm div} u\|^2  \rt.  \\
 & \lt. \qquad  + \|\na^3 {\rm div} u\|^2   + \sum_{i=1}^2 \| \na^i D_t p\|^2   + \sum_{i=1}^{3}\lt(\|\na^i u\|^2  + \|\na^i p\|^2 \rt) \rt) ,
\end{split}\end{equation}
\begin{equation}\label{17div4}\begin{split}
&  \frac{d}{dt}\int_\Omega |\na  D_t^3 {\rm div} u|^2  d\mu_g \le C( K_1,\| \mathcal{T}^{-1}\|_{L^\iy}, \|  \mathcal{T} \|_{L^\iy}, \|\na  u\|_{L^\iy} , \|\na  \mathcal{T} \|_{L^\iy},     \|\na  p \|_{L^\iy},  \\
&\qquad\|\na^2 \mathcal{T} \|_{L^\iy},   \|\na  {\rm div} u\|_{L^\iy}, \| D_t {\rm div} u \|_{L^\iy},\|\na D_t  p \|_{L^\iy}, \|\na D_t^2 {\rm div} u \|,\| D_t^2 {\rm div} u \|, \\
&\qquad
\|\na^2 D_t  {\rm div} u \|,  \|\na D_t  {\rm div} u \|,  \|\na^2    {\rm div} u \| ) \lt( \sum_{i=0}^{1} \| \na^i D_t^3 {\rm div} u\|^2
       +\sum_{i=0}^{2} \sum_{j=1}^{2}  \| \na^i D_t^j {\rm div} u\|^2   \rt. \\
 & \lt. \qquad +\sum_{i=1}^{2} \| \na^i D_t^2 p\|^2 +  \sum_{i=1}^{3}   \| \na^i D_t p\|^2 + \sum_{i=1}^{4}\lt(\|\na^i u\|^2   + \|\na^i p\|^2  \rt)   \rt) ,
\end{split}\end{equation}
\begin{equation}\label{17div5}\begin{split}
&  \frac{d}{dt}\int_\Omega |\na  D_t^4 {\rm div} u|^2  d\mu_g \le C( K_1,\| \mathcal{T}^{-1}\|_{L^\iy}, \|  \mathcal{T} \|_{L^\iy}, \|\na  u\|_{L^\iy} , \|\na  \mathcal{T} \|_{L^\iy},     \|\na  p \|_{L^\iy},  \\
&\qquad   \|\na^2 \mathcal{T} \|_{L^\iy}, \| D_t {\rm div} u \|_{L^\iy}, \|\na D_t  p \|_{L^\iy}, \|\na^2 u  \|_{L^\iy},\|\na^2 p \|_{L^\iy},\|\na^2   {\rm div} u \|_{L^\iy},  \\
&\qquad
\| D_t^2 {\rm div} u \|_{L^\iy},\|\na^2 D_t  p \|_{L^\iy} )     \lt( \sum_{i=0}^{1} \| \na^i D_t^4 {\rm div} u\|^2
       +\sum_{i=0}^{2} \sum_{j=1}^3 \| \na^i D_t^j {\rm div} u\|^2  \rt. \\
 & \lt. \qquad +\sum_{i=1}^{2} \| \na^i D_t^3 p\|^2 +  \sum_{i=1}^{3} \sum_{j=1}^{2}  \| \na^i D_t^j p\|^2 + \sum_{i=1}^{5}\lt(\|\na^i u\|^2   + \|\na^i p\|^2  \rt)   \rt) .
\end{split}\end{equation}
\end{lem}

{\em Proof}. Notice the following identity: for $r\ge 0$,
\begin{align}
&\frac{1}{2}D_t|\na D_t^r {\rm div} u|^2 + \mathcal{T}|\Delta D_t^r {\rm div} u|^2
=(D_t^{r+1}{\rm div} u- \mathcal{T}\Delta D_t^r {\rm div} u)(-\Delta D_t^r {\rm div} u)
\notag\\
&\qquad+{\rm div}\lt(  (D_t^{r+1} {\rm div} u) \na D_t^r {\rm div} u\rt) +\frac{1}{2}(D_t g^{ab})(\na_a  D_t^r {\rm div} u) \na_b D_t^r {\rm div} u.\notag
\end{align}
This, together with  \ef{s3.7}, \ef{bdry} and \ef{CL-2.13}, implies that,  for $r\ge 0$,
\begin{align}
 & \frac{d}{dt}\int_\Omega |\na D_t^r {\rm div} u|^2 d\mu_g+ \int_\Omega \mathcal{T}|\Delta D_t^r {\rm div} u|^2   d \mu_g \le  3\|\na u\|_{L^\iy} \int_\Omega |\na D_t^r {\rm div} u|^2 d\mu_g
\notag \\
& \qquad  +  \| \mathcal{T}^{-1}\|_{L^\iy} \int_\Omega |D_t^{r+1}{\rm div} u- \mathcal{T}\Delta D_t^r {\rm div} u|^2 d\mu_g   .\label{egnadivu}
\end{align}

By virtue of \ef{e9}, we have
\begin{align}
 \| D_t{\rm div} u- \mathcal{T}  \Delta {\rm div} u \|^2    \le   C   \lt(\|\na  u\|_{L^\iy},  \|\na \mathcal{T} \|_{L^\iy}\rt) \lt(  \|\na^2 \mathcal{T} \|^2 +\|\na^2 u\|^2 + \|\na u\|^2   \rt) ,  \label{4/24-1'}
 \end{align}
which, together with \ef{egnadivu}, yields \ef{17div1}.

It follows from  \eqref{dt2divu} that
\begin{align}
  \| D_t^2 {\rm div} u - \mathcal{T} \Delta  D_t {\rm div} u \|^2  \le  C  (  \| \mathcal{T} \|_{L^\iy},\|\na  u\|_{L^\iy},  \|\na \mathcal{T} \|_{L^\iy}, \|\na p \|_{L^\iy},\|\na^2 \mathcal{T} \|_{L^\iy},\notag\\
  \|\na {\rm div} u \|_{L^\iy})
 \lt(\sum_{i=0}^1 \| \na^i D_t{\rm div} u\|^2   +  \sum_{i=1}^{2}\lt(\|\na^i u\|^2    + \|\na^i p\|^2  \rt)+ \|\na^2 {\rm div} u \|^2 \rt),
    \label{new4/24-3}
\end{align}
which, together with \ef{egnadivu}, gives \ef{17div2}.

It follows from \eqref{dt3divu}, \ef{CL-A.12}, H$\ddot{o}$lder's inequality and Young's inequality that
\begin{align}
&\|D_t^3 {\rm div} u- \mathcal{T} \Delta  D_t^2 {\rm div} u \|^2 \le
C( K_1, \|  \mathcal{T} \|_{L^\iy}, \|\na  u\|_{L^\iy} , \|\na  \mathcal{T} \|_{L^\iy},     \|\na  p \|_{L^\iy},  \|\na^2 \mathcal{T} \|_{L^\iy}, \notag\\
&\qquad\|\na {\rm div} u \|_{L^\iy}, \| D_t {\rm div} u \|_{L^\iy})\lt(\sum_{i=0}^1 \| \na^i D_t^2 {\rm div} u\|^2 +  \sum_{i=0}^2 \| \na^i D_t{\rm div} u\|^2   \rt. \notag\\
 & \lt. \qquad   + \sum_{i=1}^2 \| \na^i D_t p\|^2   + \sum_{i=1}^{3}\lt(\|\na^i u\|^2  + \|\na^i p\|^2  + \|\na^i {\rm div} u\|^2\rt) \rt) , \label{new4/9-3}
\end{align}
which, together with  \ef{egnadivu}, yields \ef{17div3}.
Indeed, the following type of estimates have been used to derive \ef{new4/9-3}.
\be\label{17bala1}\begin{split}
& \|(\na^2 p) \cdot \na^2 u\|^2 \le  \|\na^2 p\|_{L^4}^2  \|\na^2 u\|_{L^4}^2
  \\
   \le &
C(K_1)\|\na  p\|_{L^\iy} \sum_{i=1}^3 \|\na^i   p\|   \|\na u\|_{L^\iy}  \sum_{i=1}^3 \|\na^i   u\|  \\
\le&  C(K_1,\|\na  p\|_{L^\iy},\|\na u\|_{L^\iy})\sum_{i=1}^3\lt(\|\na^i  p\|^2 +\|\na^i  u\|^2\rt),
\end{split}\ee
and
\be\label{17bala2}\begin{split}
& \| (\na D_t{\rm div} u ) \cdot \na^2 u    \|^2 \le   \|\na D_t{\rm div} u \|_{L^4}^2  \|\na^2 u\|_{L^4}^2
  \\ \le &
C(K_1)\|D_t{\rm div} u \|_{L^\iy} \sum_{i=0}^2 \|\na^i D_t{\rm div} u \|   \|\na u\|_{L^\iy}  \sum_{i=1}^3 \|\na^i   u\| \\
\le&  C(K_1,\|D_t{\rm div} u\|_{L^\iy},\|\na u\|_{L^\iy})\sum_{i=1}^3\lt(\|\na^{i-1}  D_t{\rm div} u\|^2 +\|\na^i  u\|^2\rt).
\end{split}\ee

By virtue of \ef{dt3divu}, \ef{CL-A.15}, \ef{CL-A.12}, H$\ddot{o}$lder's inequality and Young's inequality, we have
\begin{align}
&\|D_t^4 {\rm div} u- \mathcal{T} \Delta  D_t^3 {\rm div} u \|^2 \le
C( K_1, \|  \mathcal{T} \|_{L^\iy}, \|\na  u\|_{L^\iy} , \|\na  \mathcal{T} \|_{L^\iy},     \|\na  p \|_{L^\iy}, \|\na^2 \mathcal{T} \|_{L^\iy}, \notag\\
&\quad  \|\na  {\rm div} u\|_{L^\iy}, \| D_t {\rm div} u \|_{L^\iy},\|\na D_t  p \|_{L^\iy}, \|\na D_t^2 {\rm div} u \|,\| D_t^2 {\rm div} u \|, \|\na^2 D_t  {\rm div} u \|,\notag\\
&\quad \|\na D_t  {\rm div} u \|,  \|\na^2    {\rm div} u \| ) \lt( \sum_{i=0}^{1} \| \na^i D_t^3 {\rm div} u\|^2
       +\sum_{i=0}^{2} \lt( \| \na^i D_t^2 {\rm div} u\|^2 +  \| \na^i D_t  {\rm div} u\|^2 \rt) \rt.\notag\\
 & \lt. \quad +\sum_{i=1}^{2} \| \na^i D_t^2 p\|^2 +  \sum_{i=1}^{3}   \| \na^i D_t p\|^2 + \sum_{i=1}^{4}\lt(\|\na^i u\|^2   + \|\na^i p\|^2  \rt)   \rt),\label{nt2}
\end{align}
which, together with   \ef{egnadivu}, gives \ef{17div4}.
In addition to \ef{17bala1} and \ef{17bala2}, the following type of estimates have been used to derive \ef{nt2}.
\be\label{17bala4}\begin{split}
 \| (\na D_t^2 {\rm div} u ) \cdot \na^2 u    \|  \le   \|\na D_t^2{\rm div} u \|   \|\na^2 u\|_{L^\iy}  \le
C(K_1)  \|\na D_t^2{\rm div} u \|   \sum_{i=2}^4 \|\na^i u\|,
\end{split}\ee
\be\label{17bala5}\begin{split}
 \|    ( \na^2  {\rm div} u  ) D_t^2 {\rm div} u  \|  \le  \| \na^2  {\rm div} u\|   \|  D_t^2{\rm div} u \|_{L^\iy}     \le
C(K_1)  \| \na^2  {\rm div} u\|    \sum_{i=0}^2 \|\na^i  D_t^2{\rm div} u \|,
\end{split}\ee
and
\be\label{17bala3}\begin{split}
& \|(\na^2 p) \cdot \na^3 u\|  \le  \|\na^2 p\|_{L^6}   \|\na^3 u\|_{L^3}
  \\
   \le &
C(K_1)\|\na  p\|_{L^\iy}^{ {2}/{3}} \lt(\sum_{i=1}^4 \|\na^i   p\| \rt)^{ {1}/{3}}  \|\na u\|_{L^\iy}^{ {1}/{3}}  \lt(\sum_{i=1}^4 \|\na^i   u\| \rt)^{ {2}/{3}}   \\
\le&  C(K_1,\|\na  p\|_{L^\iy},\|\na u\|_{L^\iy})\sum_{i=1}^4\lt(\|\na^i  p\|  +\|\na^i  u\| \rt).
\end{split}\ee

Similarly, we can obtain \ef{17div5}.
\hfill $\Box$

\begin{lem}\label{lemloworder} It holds that
\be\frac{d}{dt} \int_\Omega \mathcal{T}^{-1}  |u|^2 d\mu_g \le
     \|p\|^2 +   \|{\rm div} u\|^2    , \label{4/12-9'}\ee
\be\begin{split}
  & \frac{d}{dt}\int_\Omega \mathcal{T}^{-1} g^{ac}\zeta^{bd} (\nabla_b u_a) \nabla_d u_c  d\mu_g
+ \frac{d}{dt}\int_\Omega |{\rm curl} u|^2 d\mu_g   \\
& \quad \le     C( K, \|\mathcal{T}^{-1}\|_{L^\iy}, \|\nabla u\|_{L^\iy},\|\nabla \mathcal{T}\|_{L^\iy})
      \lt(   \|\nabla u\|^2 +    \|\na {\rm div} u \|^2 + \|\na p\|^2 \rt) . \label{4/12-9}
\end{split}\ee
\end{lem}

{\em Proof}. It follows from \ef{e6}, \ef{e7} and \ef{e25} that
\begin{align*}
& \frac{1}{2}D_t\lt(\mathcal{T}^{-1}  |u|^2  \rt)+ \frac{1}{2} \mathcal{T}^{-1}  |u|^2 {\rm div} u =-  {\rm div}(pu)+ p {\rm div} u
\end{align*}
and
\begin{align*}
& \frac{1}{2}D_t\lt( \mathcal{T}^{-1} g^{ac}\zeta^{bd}  (\nabla_b u_a) \nabla_d u_c \rt) + \frac{1}{2}  \mathcal{T}^{-1} g^{ac}\zeta^{bd}  (\nabla_b u_a) (\nabla_d u_c ){\rm div} u \notag\\
  = & - {\rm div} \lt(  \zeta^{bd} (\na_b p) \na_d u   \rt)  +  \frac{1}{2} \mathcal{T}^{-1} \lt(D_t\lt(g^{ac}\zeta^{bd}  \rt) \rt) (\nabla_b u_a) \nabla_d u_c +  (\na_d u) \cdot (\na \zeta^{bd}) \na_b p   \notag \\
& \ +  \zeta^{bd}    (\na_b p) \na_d {\rm div} u  + \zeta^{bd} \mathcal{T}^{-1} \lt( (\na_d u)\cdot (\na u_e) \na_b u^e - (\na_b \mathcal{T} ) (\na_d u)\cdot \na p    \rt) ,\notag
\end{align*}
which, together with \ef{s3.7}, \ef{bdry},   \ef{CL-2.13} and \ef{CL-3.28}, imply that
\begin{align}
\frac{d}{dt} \int_\Omega \mathcal{T}^{-1}  |u|^2 d\mu_g
= 2\int_\Omega p {\rm div} u d\mu_g \le \int_\Omega (p^2 + |{\rm div} u|^2 ) d\mu_g  \label{170214-1}
\end{align}
and
\begin{align}
 &\frac{d}{dt}\int_\Omega  \mathcal{T}^{-1} g^{ac}\zeta^{bd} (\nabla_b u_a) \nabla_d u_c  d\mu_g \notag\\
\le & C( K, \|\mathcal{T}^{-1}\|_{L^\iy}, \|\nabla u\|_{L^\iy},\|\nabla \mathcal{T}\|_{L^\iy})
    \int_\Omega \lt(   |\nabla u|^2 +  (|\na u|  + |\na {\rm div} u |) |\na p| \rt) d\mu_g \notag \\
    \le &  C( K, \|\mathcal{T}^{-1}\|_{L^\iy}, \|\nabla u\|_{L^\iy},\|\nabla \mathcal{T}\|_{L^\iy})
    \int_\Omega \lt(   |\nabla u|^2 +    |\na {\rm div} u |^2 + |\na p|^2 \rt) d\mu_g \label{170214-2}.
\end{align}
It follows from \ef{s3.7}, \ef{CL-2.13} and  \ef{e27} that
\begin{align}
 \frac{d}{dt}\int_\Omega |{\rm curl} u|^2 d\mu_g  \le & C(\|\nabla u\|_{L^\iy})   \int_\Omega \lt( |{\rm curl} u|  | D_t {\rm curl} u | +    |{\rm curl} u|^2 \rt)d\mu_g  \notag\\
\le   & C(\|\nabla u\|_{L^\iy},\|\nabla \mathcal{T}\|_{L^\iy})  \int_\Omega \lt( |\na u|^2  + |\na p |^2 \rt)d\mu_g .  \notag
\end{align}
This, together with \ef{170214-1} and \ef{170214-2}, gives \ef{4/12-9'} and \ef{4/12-9}.
$\hfill \Box$

\begin{lem} For $r\ge 2$, we have
\begin{align}
&  \frac{d}{dt}\int_\Omega  \mathcal{T}^{-1} g^{ab}\zeta^{cd} \zeta^{IJ} (\na^{r-1}_I\nabla_c u_a) \na^{r-1}_J\nabla_d u_b  d\mu_g + \frac{d}{dt}\int_\Omega  |\na^{r-1} {\rm curl} u|^2 d\mu_g  \notag\\
& +  \frac{d}{dt} \int_{\pl\Omega}      \zeta^{cd} \zeta^{IJ} (\na^{r -1 }_I \na_c p)( \na^{r -1}_J \na_d  p )(-\na_N p)^{-1} d \mu_\za \notag \\
\le & C( K, \|\mathcal{T}^{-1}\|_{L^\iy}, \|\nabla u\|_{L^\iy}, |\| (\na_N p)^{-1} \||_{L^\iy}, |\|\nabla_N p\||_{L^\iy}, |\|\nabla_N D_t p \||_{L^\iy}, |\|\nabla u\||_{L^\iy})  \notag\\
&  \times \lt(\| \na^r u\|^2 +  \| \na^r {\rm div} u\|^2 + \|D_t \na^{r }u + \mathcal{T}\na^{r+1}p\|^2    + \|D_t \na^{r-1}  {\rm curl} u\|^2 + \| \na^r p\|^2 \rt.\notag\\
& \lt. + |\| \Pi  \na^r p\||^2 +
|\|\Pi (D_t \na^{r }p + (\na^r u) \cdot \na p)\||^2    \rt). \label{4/13-4}
\end{align}

\end{lem}
{\em Proof}. When $r\ge 1$, simple calculations give that
\begin{align}
& \frac{1}{2}D_t\lt( \mathcal{T}^{-1} g^{ab}\zeta^{cd} \zeta^{IJ} (\na^r_I\nabla_c u_a) \na^r_J\nabla_d u_b \rt) + \frac{1}{2} \mathcal{T}^{-1} g^{ab}\zeta^{cd} \zeta^{IJ} (\na^r_I\nabla_c u_a) (\na^r_J\nabla_d u_b ){\rm div} u \notag\\
&  =  - {\rm div} \lt(  \zeta^{cd} \zeta^{IJ} (\na^r_I \na_c p) \na^r_J \na_d u   \rt)  +  \frac{1}{2} \mathcal{T}^{-1} \lt(D_t\lt(g^{ab}\zeta^{cd} \za^{IJ}  \rt) \rt) (\na^r_I\nabla_c u_a) \na^r_J \nabla_d u_b  \notag \\
& \quad + \mathcal{T}^{-1}  g^{ab}\zeta^{cd} \za^{IJ}  (D_t \na^r_I\nabla_c u_a + \mathcal{T}\na^r_I\na_c\na_a p) \na^r_J \nabla_d u_b
\notag\\
& \quad +  ( \na^r_J \nabla_d u )\cdot \lt( \na (\za^{cd}\za^{IJ}) \rt) (\na^r_I \na_c p)
+  \za^{cd}\za^{IJ} (\na^r_I \na_c p)  \na^r_J \nabla_d {\rm div} u . \label{4/12-11}
\end{align}
Due to \ef{bdry}, we have $\na_a p =\overline{\na}_a p + N_a \na_N p = N_a \na_N p$ on $\pl \Omega$ and
\begin{align}
&\int_\Omega  {\rm div} \lt(  \zeta^{cd} \zeta^{IJ} (\na^r_I \na_c p) \na^r_J \na_d u   \rt) d \mu_g  = \int_{\pl\Omega}    N_a  \zeta^{cd} \zeta^{IJ} (\na^r_I \na_c p)( \na^r_J \na_d u^a ) d \mu_\za  \notag \\
& = - \int_{\pl\Omega}    N_a  (\na_N p) \zeta^{cd} \zeta^{IJ} (\na^r_I \na_c p)( \na^r_J \na_d u^a )(-\na_N p)^{-1} d \mu_\za \notag\\
& =  -  \int_{\pl\Omega}     \zeta^{cd} \zeta^{IJ} (\na^r_I \na_c p)(  D_t \na^r_J \na_d  p +  (\na^r_J \na_d u^a)\na_a p )(-\na_N p)^{-1} d \mu_\za \notag\\
&\quad
+\int_{\pl\Omega}     \zeta^{cd} \zeta^{IJ} (\na^r_I \na_c p)( D_t \na^r_J \na_d  p )(-\na_N p)^{-1} d \mu_\za.
\label{4/12-12}
\end{align}
Moreover, it follows from \ef{CL-3.23} that
\begin{align}
&  \int_{\pl\Omega}     \zeta^{cd} \zeta^{IJ} (\na^r_I \na_c p)( D_t \na^r_J \na_d  p )(-\na_N p)^{-1} d \mu_\za \notag\\
= & \frac{1}{2} \frac{d}{dt} \int_{\pl\Omega}     \zeta^{cd} \zeta^{IJ} (\na^r_I \na_c p)( \na^r_J \na_d  p )(-\na_N p)^{-1} d \mu_\za
\notag\\
& - \frac{1}{2}\int_{\pl\Omega}  \lt( D_t(  \zeta^{cd} \zeta^{IJ} (-\na_N p)^{-1} ) \rt) (\na^r_I \na_c p)( \na^r_J \na_d  p ) d \mu_\za \notag\\
& - \frac{1}{2} \int_{\pl\Omega}     \zeta^{cd} \zeta^{IJ} (\na^r_I \na_c p)( \na^r_J \na_d  p )(-\na_N p)^{-1} ({\rm div} u -h_{NN}) d \mu_\za.
\label{4/12-13}
\end{align}
Notice that on $\pl\Oa$,
\begin{align}
|D_t (-\na_N p)^{-1}|\le &  | (-\na_N p)^{-1}|^2 |D_t \na_N p|
 \le   | (-\na_N p)^{-1}|^2(|(D_t N^c)\na_c p|+|\na_N D_t p|) \notag\\
 = & | (-\na_N p)^{-1}|^2(|\lt(-2h^c_d N^d + h_{NN}N^c\rt)\na_c p|+|\na_N D_t p|)\notag\\
 = & | (-\na_N p)^{-1}|^2(| h_{NN}\na_N p|+|\na_N D_t p|),\label{dtnanp}
\end{align}
because of \ef{CL-3.21} and  \ef{bdry}.
We then obtain, with the help of \ef{4/12-11}-\ef{4/12-13}, \ef{s3.7},    \ef{CL-2.13},  \ef{CL-3.22} and \ef{CL-3.28},  that for $r\ge 2$,
 \begin{align}
&  \frac{d}{dt}\int_\Omega \mathcal{T}^{-1} g^{ab}\zeta^{cd} \zeta^{IJ} (\na^{r-1}_I\nabla_c u_a) \na^{r-1}_J\nabla_d u_b  d\mu_g  \notag\\
&+  \frac{d}{dt} \int_{\pl\Omega}     \zeta^{cd} \zeta^{IJ} (\na^{r-1}_I \na_c p)( \na^{r-1}_J \na_d  p )(-\na_N p)^{-1} d \mu_\za
\notag \\
\le & C( K, \|\mathcal{T}^{-1}\|_{L^\iy}, \|\nabla u\|_{L^\iy})  \int_\Omega \{|\na^{r } u|( |\na^{r } u| + |D_t \na^{r }u + \mathcal{T}\na^{r+1}p|+  |\na^{r }p |) \notag\\
&+ |\na^{r }p ||\na^{r }{\rm div} u| \} d\mu_g + C( |\| (\na_N p)^{-1} \||_{L^\iy}, |\|\nabla_N p\||_{L^\iy}, |\|\nabla_N D_t p \||_{L^\iy}, |\|\nabla u\||_{L^\iy})\notag\\
 & \times\int_{\pl\Omega} (  |\Pi \na^r p|  |\Pi (D_t \na^{r }p + (\na^r u) \cdot \na p) | + | \Pi\na^r p|^2  ) d\mu_\za.
 \label{4/13-1}
\end{align}
It follows from \ef{s3.7} and \ef{CL-2.13} that for $r\ge 2$
\begin{align}
& \frac{d}{dt}\int_\Omega |\na^{r-1} {\rm curl} u|^2 d\mu_g \notag\\
 \le & C(\|\nabla u\|_{L^\iy})   \int_\Omega \lt( |\na^{r-1} {\rm curl} u|^2  + |D_t \na^{r-1}  {\rm curl} u ||\na^{r-1} {\rm curl} u|  \rt)d\mu_g ,  \notag
\end{align}
which, together with \ef{4/13-1}, gives \ef{4/13-4}.
\hfill $\Box$

\begin{lem}\label{lemhighorder} For $r=2,3,4,5$, we have
 \begin{align}
&  \frac{d}{dt}\int_\Omega  \mathcal{T}^{-1} g^{ab}\zeta^{cd} \zeta^{IJ} (\na^{r-1}_I\nabla_c u_a) \na^{r-1}_J\nabla_d u_b  d\mu_g + \frac{d}{dt}\int_\Omega  |\na^{r-1} {\rm curl} u|^2 d\mu_g  \notag\\
& +  \frac{d}{dt} \int_{\pl\Omega}      \zeta^{cd} \zeta^{IJ} (\na^{r -1 }_I \na_c p)( \na^{r -1}_J \na_d  p )(-\na_N p)^{-1} d \mu_\za \notag \\
&\le  C( K, K_1,\|\mathcal{T}^{-1}\|_{L^\iy}, \|\nabla u\|_{L^\iy}, \|\nabla p\|_{L^\iy},\|\nabla \mathcal{T}\|_{L^\iy},|\| (\na_N p)^{-1} \||_{L^\iy}, |\|\nabla_N p\||_{L^\iy},\notag\\
&|\|\nabla_N D_t p \||_{L^\iy}, |\|\nabla u\||_{L^\iy}, |\|\nabla^2 p\||_{L^\iy})   \lt(\sum_{i=1}^r \lt(\|\na^i u\|^2  + \|\na^i p\|^2 + \|\na^i \mathcal{T}\|^2 \rt)   +   \| \na^r {\rm div} u\|^2\rt.\notag\\
& \lt.+ |\| \Pi  \na^r  p\||^2 + |\| \Pi  \na^r D_t p\||^2+\sum_{i=2}^{r-1}|\|{\na}^i u\||^2 +\sum_{i=3}^r|\|{\na}^i p\||^2 + \sum_{i=2}^{r-3}|\|\overline{\na}^i \ta \||^2
   \rt).  \label{hfin4.58}
\end{align}
\end{lem}

{\em Proof}. It follows from \ef{CL-6.17} and H$\ddot{\rm o}$lder's inequality   that
\begin{equation*}\begin{split}
&\|D_t \na^{r }u + \mathcal{T}\na^{r+1}p\|+ \|D_t \na^{r-1}  {\rm curl} u\| \le C( \|\nabla u\|_{L^\iy}, \|\nabla p\|_{L^\iy},  \|\na \mathcal{T}\|_{L^\iy}) \\
 & \times \lt\{ \begin{split}
&\|\na^2 u\|  + \|\na^2 p\|  + \|\na^2 \mathcal{T}\|,  & r=2, \\
&  \|\na^3 u\|  + \|\na^3 p\|  + \|\na^3 \mathcal{T}\| +  \|\na^2 u\|_{L^4}^2 +  \|\na^2  p\|_{L^4}   \|\na^2 \mathcal{T}  \|_{L^4}  ,  & r=3, \\
&  \|\na^4 u\|  + \|\na^4 p\|  + \|\na^4 \mathcal{T}\| +  \|\na^2 u\|_{L^6}  \|\na^3 u\|_{L^3}  & \\
& \quad +  \|\na^2 p\|_{L^6}  \|\na^3  \mathcal{T}\|_{L^3} + \|\na^3 p\|_{L^3}  \|\na^2 \mathcal{T}\|_{L^6}  ,  & r=4 ,\\
&  \|\na^5 u\|  + \|\na^5 p\|  + \|\na^5 \mathcal{T}\| +  \|\na^2 u\|_{L^8}  \|\na^4 u\|_{L^{\frac{8}{3}}} +  \|\na^3 u\|_{L^4}^2 & \\
& \quad +  \|\na^2 p\|_{L^8}  \|\na^4  \mathcal{T}\|_{L^{\frac{8}{3}}} +  \|\na^3 p\|_{L^4} \|\na^3 \mathcal{T}\|_{L^{4}} +   \|\na^4  p\|_{L^{\frac{8}{3}}} \|\na^2 \mathcal{T}\|_{L^8}   ,  & r=5 ,
\end{split}\rt.
\end{split}\end{equation*}
which, together with \ef{CL-A.12} and Young's inequality, implies for $r=2,3,4,5$,
\be\label{4/13-5}\begin{split}
& \|D_t \na^{r }u + \mathcal{T}\na^{r+1}p\|+ \|D_t \na^{r-1}  {\rm curl} u\| \\
\le & C(K_1, \|\nabla u\|_{L^\iy}, \|\nabla p\|_{L^\iy},  \|\na \mathcal{T}\|_{L^\iy}) \sum_{i=1}^r \lt(\|\na^i u\|  + \|\na^i p\|  + \|\na^i \mathcal{T}\| \rt)  .
\end{split}\ee
It follows from \ef{CL-6.18'} and \ef{bdry} that
\begin{equation}\label{ww4.64}\begin{split}
&  |\|\Pi (D_t \na^{r}p + (\na^r u) \cdot \na p)\|| \le  |\|\Pi  \na^{r }D_t p  \||  \\
& \quad + \lt\{ \begin{split}
& 0   ,  & r=2, \\
&      C    |\|\na^2 u\|| |\|\na^2  p\||_{L^\iy} ,  & r=3, \\
&     C(   |\|\na^3 u\|| |\|\na^2  p\||_{L^\iy}  +
 |\|\Pi((\na^2 u)\cdot \na^3  p)\|| ) ,  & r=4 ,\\
&   C(   |\|\na^4 u\|| |\|\na^2  p\||_{L^\iy}  +   |\|\Pi((\na^3 u)\cdot \na^3  p)\|| +
 |\|\Pi((\na^2 u)\cdot \na^4  p)\|| ) ,  & r=5
 .\\
\end{split}\rt.
\end{split}\end{equation}
Notice that on $\pl\Oa$,
\begin{align}
& |\Pi((\na^2 u)\cdot \na^3  p)|\le |\overline{\na}\na u||\overline{\na}\na^2 p|,\notag\\
& |\Pi((\na^3 u)\cdot \na^3  p)| \le |\overline{\na}^2\na u||\overline{\na}\na^2 p|+ CK |\overline{\na}\na u||\overline{\na}\na^2 p| , \notag\\
& | \Pi((\na^2 u)\cdot \na^4  p) | \le |\overline{\na}\na u||\overline{\na}^2 \na^2 p|+ CK |\overline{\na}\na u||\overline{\na}\na^2 p| . \notag
 \end{align}
We then obtain, with the help of H$\ddot{\rm o}$lder's inequality, \ef{CL-A.4} and Young's inequality, that
\begin{align}
&|\|\Pi((\na^2 u)\cdot \na^3  p)\||\le     |\|\overline{\na}\na u\||_{L^4}|\|\overline{\na}\na^2 p\||_{L^4} \notag\\
\le & C
|\| \na u\||_{L^\iy}^{\frac{1}{2}} |\|\overline{\na}^2 \na u\||^{\frac{1}{2}} |\|\na^2 p\||_{L^\iy}^{\frac{1}{2}}  |\|\overline{\na}^2 \na^2 p\||^{\frac{1}{2}}\notag \\
\le & C\lt( |\| \na u\||_{L^\iy},  |\|\na^2 p\||_{L^\iy}\rt) \lt( |\|\overline{\na}^2 \na u\|| +|\|\overline{\na}^2 \na^2 p\|| \rt), \label{ww4.65}
\end{align}
\begin{align}
& |\|\Pi((\na^3 u)\cdot \na^3  p) \||  + |\|\Pi((\na^2 u)\cdot \na^4  p) \|| \notag\\
 \le  &   |\|\overline{\na}^2\na u\||_{L^3}|\|\overline{\na}\na^2 p\||_{L^6}  + |\|\overline{\na}\na u\||_{L^6}|\|\overline{\na}^2 \na^2 p\||_{L^3}  + CK  |\|\overline{\na}\na u\||_{L^4}|\|\overline{\na}\na^2 p\||_{L^4} \notag\\
\le & C
|\| \na u\||_{L^\iy}^{\frac{1}{3}} |\|\overline{\na}^3 \na u\||^{\frac{2}{3}} |\|\na^2 p\||_{L^\iy}^{\frac{2}{3}}  |\|\overline{\na}^3 \na^2 p\||^{\frac{1}{3}}+ CK  |\|\overline{\na}\na u\||_{L^4}|\|\overline{\na}\na^2 p\||_{L^4} \notag\\
& + C
|\| \na u\||_{L^\iy}^{\frac{2}{3}} |\|\overline{\na}^3 \na u\||^{\frac{1}{3}} |\|\na^2 p\||_{L^\iy}^{\frac{1}{3}}  |\|\overline{\na}^3 \na^2 p\||^{\frac{2}{3}}\notag \\
\le & C( K, |\| \na u\||_{L^\iy},  |\|\na^2 p\||_{L^\iy}) \lt(\sum_{i=2}^3|\|\overline{\na}^i \na u\|| +\sum_{i=2}^3||\|\overline{\na}^i \na^2 p\||   \rt). \label{170215-1}
\end{align}
For a $(0,2)$ tensor $\al$, we have on $\pl\Oa$,
$$|\overline{\na}^2 \al| \le |\na^2 \al| + C K |\na \al| \ \ {\rm and} \ \ |\overline{\na}^3 \al| \le |\na^3 \al| + C \lt(K |\na^2  \al| + K^2 |\na \al | + |\overline{\na} \ta||\overline{\na} \al| \rt).$$
This, together with \ef{ww4.64}, \ef{ww4.65} and \ef{170215-1}, implies that
\begin{equation*}\label{}\begin{split}
&  |\|\Pi (D_t \na^{r}p + (\na^r u) \cdot \na p)\|| \le  |\|\Pi  \na^{r }D_t p  \||  \\
&  + \lt\{ \begin{split}
&   C\lt( |\| \na u\||_{L^\iy},  |\|\na^2 p\||_{L^\iy}\rt) \lt(\sum_{i=2}^3|\|{\na}^i u\|| +\sum_{i=3}^4|\|{\na}^i p\||   \rt)  ,  & r=4 ,\\
&    C( K, |\| \na u\||_{L^\iy},  |\|\na^2 p\||_{L^\iy}) \lt(\sum_{i=2}^4|\|{\na}^i u\|| +\sum_{i=3}^5|\|{\na}^i p\|| + |\|\overline{\na}^2 \ta \||  \rt)  ,  & r=5,\\
\end{split}\rt.
\end{split}\end{equation*}
which implies, using \ef{ww4.64}, that for $r=2,3,4,5$,
\begin{align}
&|\|\Pi (D_t \na^{r}p + (\na^r u) \cdot \na p)\|| \le  |\|\Pi  \na^{r }D_t p  \||
\notag\\
& \; + C( K, |\| \na u\||_{L^\iy},  |\|\na^2 p\||_{L^\iy}) \lt(\sum_{i=2}^{r-1}|\|{\na}^i u\|| +\sum_{i=3}^r|\|{\na}^i p\|| + \sum_{i=2}^{r-3}|\|\overline{\na}^i \ta \||  \rt). \label{170215-2}
\end{align}
So, \ef{hfin4.58} follows from \ef{4/13-4}, \ef{4/13-5} and \ef{170215-2}.
 \hfill $\Box$

\subsection{Elliptic estimates}\label{17elliptic}
Before starting with the proof of \ef{17prop1}, let us first see what a bound for the energy \ef{17energy}  implies.

\begin{lem}\label{17lem01} We have
\be\label{4/10-7}
\|D_t^{r-1}  {\rm div} u\|^2 \le C({\rm Vol} \Omega) \|\na D_t^{r-1}  {\rm div} u\|^2  \le  C({\rm Vol} \Omega) E_r(t), \ \ \ r\ge 1,
\ee
\be\label{17t2}
 \| \mathcal{T}-\mathcal{T}_b\|^2+ \|\na  \mathcal{T}\|^2 + \|\na^2  \mathcal{T}\|^2  \le C(K,{\rm Vol} \Oa) E_1,\ee
\be\label{17u01}\|  u\|^2 \le  \|  \mathcal{T} \|_{L^\iy}  E_0, \ \  \| \na u\|^2 \le C(  {\rm Vol} \Omega ,  \|  \mathcal{T} \|_{L^\iy} )  E_1, \ee
\be\label{17but}|\| u\||^2 \le C(K,  {\rm Vol} \Omega ,  \|  \mathcal{T} \|_{L^\iy} ) (E_0+ E_1), \  \
|\| \na  \mathcal{T}\||^2 \le C(K,  {\rm Vol} \Omega  ) E_1.
\ee
\end{lem}

{\em Proof}. \ef{4/10-7} follows from \ef{bdry} and \ef{CL-A.17}. \ef{17t2} follows from
 \ef{bdry}, \ef{310-1},  \ef{e7} and \ef{4/10-7}. \ef{17u01} follows from \ef{CL-5.16} and \ef{4/10-7}. \ef{17but} follows from \ef{CL-5.19}, \ef{17u01} and \ef{17t2}. \hfill $\Box$

\begin{lem}\label{17lem2} We have
\be
 \| \na^2 u\|^2 +|\| \na u\||^2 \le C(  {\rm Vol} \Omega ,   \|  \mathcal{T} \|_{L^\iy} ) \lt( E_1+ E_2\rt),  \label{r-4.14}
 \ee
 \be
  \|  p \|^2 +\| \na p \|^2    \le  C({\rm Vol} \Oa,  \|  \mathcal{T}^{-1} \|_{L^\iy}, \|  \mathcal{T} \|_{L^\iy},\|\na  u\|_{L^\iy}   ) \lt(E_1 + E_2\rt), \label{17p1}
  \ee
\be \label{ediv2} \begin{split}
& \|\na^2 {\rm div} u \|^2 + |\| \na {\rm div} u  \||^2 + \| \na^2 p \|^2 +|\|\na_N p \||^2 \\
&\quad \le  C( K,{\rm Vol} \Oa,  \|  \mathcal{T}^{-1} \|_{L^\iy}, \|  \mathcal{T}\|_{L^\iy}, \|\na  u\|_{L^\iy} , \|\na  \mathcal{T} \|_{L^\iy}  )\lt( E_1+E_2\rt) ,
\end{split}\ee
\be
|\|\Pi \na^r p\||^2 \le  |\| \na_N p\||_{L^\iy} E_r, \ \  r\ge 2, \label{rt2}
\ee
\be\label{17bp2}\begin{split}
 |\|\na^2 p \||^2  \le & C( K,K_1, {\rm Vol} \Oa,  \|  \mathcal{T}^{-1} \|_{L^\iy}, \|  \mathcal{T} \|_{L^\iy},\|\na  u\|_{L^\iy} , \|\na \mathcal{T} \|_{L^\iy} , \|\na p \|_{L^\iy}, \\ &  |\|\na_N p\||_{L^\iy} ) \lt(E_1 + E_2\rt) ,
\end{split}\ee
\be
|\|  \ta \||^2
\le    |\| (\na_N p)^{-1}\||_{L^\iy}^2 |\| \na_N p\||_{L^\iy} E_2 ,\label{17ta}
\ee
\be|\|\na^2 \mathcal{T} \||^2  \le  C( K,{\rm Vol} \Oa,   |\|\na_N p\||_{L^\iy} , |\|(\na_N p)^{-1}\||_{L^\iy},  |\|\na_N  \mathcal{T}\||_{L^\iy}) \lt(E_1 + E_2\rt) , \label{17bt2}
\ee
\be
|\|\Pi \na^2 D_t p\||^2 \le C(  |\| \na_N D_t p\||_{L^\iy},
|\|\na_N p\||_{L^\iy} , |\|(\na_N p)^{-1}\||_{L^\iy}  )E_2 . \label{17dtp}
\ee
\end{lem}

{\em Proof}. The bound for $\|\na^2 u\|$ in \ef{r-4.14} follows from \ef{CL-5.16}. The bound for $|\|\na u\||$ in \ef{r-4.14} follows from \ef{CL-5.19}, \ef{17u01} and the  bound just obtained for $\|\na^2 u\|$.
The bound for $\|\na^2 {\rm div} u\|$ in \ef{ediv2} follows from \ef{4/24-1'}, \ef{310-1},  \ef{r-4.14} and Lemma \ref{17lem01}.
The bound for $|\|\na {\rm div} u\||$ in \ef{ediv2} follows from \ef{CL-5.19} and the bound just obtained for $\|\na^2 {\rm div} u\|$. Let $q$ be a function satisfying $q=0$ on $\pl\Omega$, we have for any $\da>0$,
\begin{align}
\int_\Oa  \mathcal{T} | \na q |^2 d\mu_g= & \int_\Omega  {\rm div}(\mathcal{T} q \na q ) d\mu_g -\int_\Omega (\mathcal{T}   \Delta q +(\nabla  \mathcal{T})\cdot \nabla q) q d\mu_g  \notag\\
\le & \frac{1}{2 \da }\int_\Omega |\mathcal{T}   \Delta q +(\nabla  \mathcal{T})\cdot \nabla q|^2 d\mu_g + \frac{\da}{2}\int_\Omega q^2 d\mu_g.  \notag
\end{align}
Due to \ef{CL-A.17}, we can then choose a  suitably small $\da$    to obtain
\begin{align}
\| \na q\|^2\le C({\rm Vol}\Oa,  \|\mathcal{T}^{-1}\|_{L^\iy})  \|\mathcal{T}   \Delta q +(\nabla  \mathcal{T})\cdot \nabla q\|^2  . \label{4/10-6}
\end{align}
It follows from \ef{lappi} that
\be
\|\mathcal{T} \Delta p + (\na \mathcal{T})\cdot \na p\| \le \| D_t {\rm div} u \|+ \|\na u\|_{L^\iy} \|\na u\|, \label{17nap}
\ee
which, together with  \ef{4/10-6}, \ef{4/10-7} and \ef{r-4.14}, gives the bound for $\|\na p\|$ in \ef{17p1}.
The bound for $\| p\|$ in \ef{17p1} follows from \ef{CL-A.17} and the bound just obtained for $\|\na p\|$. The bound for $\|\na^2 p\|$ in \ef{ediv2} follows from  \ef{310-1}, \ef{17nap},  \ef{17p1} and Lemma \ref{17lem01}.
The bound for $|\|\na_N p\||$ in \ef{ediv2} follows from \ef{CL-5.19}, \ef{17p1} and the bound  just obtained for $\|\na^2 p\|$. Clearly, \ef{rt2} holds.
Due to \ef{CL-6.17}, \ef{CL-A.12}, H$\ddot{\rm o}$lder's inequality and Young's inequality, we have for $r=1,2,3,4$,
\begin{align}
     \|  \mathcal{T} \na^r \Delta p \| ^2  \le    C     \| \na^r D_t{\rm div} u \| ^2 +  C(K_1, \|\na  u\|_{L^\iy}, \|\na \mathcal{T} \|_{L^\iy}, \|\na  p\|_{L^\iy}) \notag \\
       \times\sum_{i=1}^{r+1}\lt(\|\na^i u\|^2
     +  \|\na^i \mathcal{T} \|^2 + \|\na^i p\|^2 \rt)  . \label{4/24-2}
\end{align}
\ef{17bp2} follows from \ef{CL-5.28}, \ef{rt2}, \ef{4/24-2}, \ef{ediv2}, \ef{17t2} and \ef{r-4.14}. \ef{17ta} follows from \ef{4.25.1} and \ef{rt2}.
\ef{17bt2} follows from \ef{CL-5.28}, \ef{hb3}, \ef{e7}, \ef{17ta} and \ef{4/10-7}. \ef{17dtp} follows from \ef{hb3} and \ef{17ta}.
\hfill $\Box$

\begin{lem}\label{17lem3} We have
\be\begin{split}
 \|\na^3  u \|^2 + |\| \na^2    u\||^2 \le    C(   K,{\rm Vol} \Oa,  \|  \mathcal{T}^{-1} \|_{L^\iy}, \|  \mathcal{T}\|_{L^\iy}, \|\na  u\|_{L^\iy} , \|\na  \mathcal{T} \|_{L^\iy}  )\sum_{i=1}^3 E_i, \label{t-1}
\end{split}\ee
\be\begin{split}
&\|   D_t  p\|^2 + \|  \na D_t  p\|^2 + \| \na^2  D_t  p\|^2 +  |\|  \na D_t  p\||^2  \\
& \qquad\le   C(K, {\rm Vol}\Oa,    \|  \mathcal{T}^{-1} \|_{L^\iy}, \|  \mathcal{T} \|_{L^\iy},  \|\na  u\|_{L^\iy} , \|\na  \mathcal{T} \|_{L^\iy},     \|\na  p \|_{L^\iy}) \sum_{i=1}^3 E_i, \label{4/25-1}
\end{split}\ee
\be\begin{split}
 \|\na^3 p\|^2
\le  &  C(   K,  K_1, {\rm Vol}\Oa,   \|  \mathcal{T}^{-1} \|_{L^\iy}, \|  \mathcal{T}  \|_{L^\iy},     \|\na  u\|_{L^\iy} , \|\na  \mathcal{T} \|_{L^\iy},  \\
   & \|\na  p \|_{L^\iy},|\| \na_N p\||_{L^\iy})\sum_{i=1}^3 E_i
  , \label{t-2}
\end{split}\ee
\be\begin{split}
  |\|\overline{\na} \ta \||^2
\le &   C(K,  K_1, {\rm Vol}\Oa,      \|  \mathcal{T}^{-1} \|_{L^\iy}, \|  \mathcal{T}  \|_{L^\iy},   \|\na  u\|_{L^\iy} , \|\na  \mathcal{T} \|_{L^\iy},   \\
&    \|\na  p \|_{L^\iy},|\| (\na_N p)^{-1}\||_{L^\iy},|\| \na_N p\||_{L^\iy} )
 \sum_{i=1}^3 E_i , \label{r-4.21}
\end{split}\ee
\be\begin{split}
&\| \na^3   \mathcal{T} \|^2 +  |\| \na^3   \mathcal{T} \||^2   +\| \na^2 D_t {\rm div} u \|^2 + |\|\na  D_t {\rm div} u\||^2+  |\| \na^3 p  \||^2
\\
   & \qquad \le      C(K,  K_1, {\rm Vol}\Oa,     \|  \mathcal{T}^{-1} \|_{L^\iy}, \|  \mathcal{T}  \|_{L^\iy},    \|\na  u\|_{L^\iy} ,  \|\na  \mathcal{T} \|_{L^\iy},  \|\na  p \|_{L^\iy}, \\
   & \qquad  |\| (\na_N p)^{-1}\||_{L^\iy},   |\| \na_N p \||_{L^\iy}, |\| \na_N \mathcal{T} \||_{L^\iy})
\sum_{i=1}^3 E_i   ,\label{step3}
\end{split}\ee
\be\begin{split}
&  |\| \na^2   {\rm div} u\||^2 + |\| \Pi \na^3   {\rm div} u\||^2 +  \|\na^3 {\rm div} u  \|^2     \\
&\qquad \le      C(K,  K_1, {\rm Vol}\Oa,      \|  \mathcal{T}^{-1} \|_{L^\iy}, \|  \mathcal{T}  \|_{L^\iy},   \|\na  u\|_{L^\iy} , \|\na  \mathcal{T} \|_{L^\iy},    \|\na  p \|_{L^\iy},   \\
&\qquad|\| (\na_N p)^{-1}\||_{L^\iy},   |\| \na_N p \||_{L^\iy}, |\| \na_N \mathcal{T} \||_{L^\iy}, |\| \na_N  {\rm div} u \||_{L^\iy})
 \sum_{i=1}^3 E_i ,  \label{step4}
\end{split}\ee
\be\begin{split}
 &  |\| \na^2 D_t p \||^2  +   |\|\Pi \na^3 D_t p \||^2 + \|\na^3 D_t p \|^2     \\
&\qquad \le    C(K,  K_1, {\rm Vol}\Oa,      \|  \mathcal{T}^{-1} \|_{L^\iy}, \|  \mathcal{T}  \|_{L^\iy},  \|\na  u\|_{L^\iy} ,   \|\na  \mathcal{T} \|_{L^\iy},   \|\na  p \|_{L^\iy},  \| D_t  p \|_{L^\iy},  \\
&\qquad |\| (\na_N p)^{-1}\||_{L^\iy},  |\| \na_N p \||_{L^\iy}, |\| \na_N \mathcal{T} \||_{L^\iy},  |\| \na_N D_t p \||_{L^\iy})
\sum_{i=1}^3 E_i.  \label{step4''}
\end{split}\ee

\end{lem}

{\em Proof}. The bound for $\|\na^3 u\|$ in \ef{t-1}  follows from \ef{CL-5.16} and \ef{ediv2}. The bound for $|\|\na^2 u\||$ in \ef{t-1}  follows from
 \ef{CL-5.19}, \ef{r-4.14} and the bound just obtained for $\|\na^3 u\|$.
The bound for $\|\na D_tp \|$ in \ef{4/25-1} follows from \ef{4/10-5}, \ef{4/10-6}, and Lemmas \ref{17lem01} and \ref{17lem2}. The bound for $\|  D_tp \|$ in \ef{4/25-1} follows from \ef{CL-A.17} and the bound just obtained for $\|\na D_tp \|$.
The bound for $\|\na^2 D_tp \|$ in \ef{4/25-1} follows from \ef{310-1}, \ef{4/10-5}, and the bound just obtained for $\|\na D_tp \|$. The bound for $|\|\na D_tp \||$ in \ef{4/25-1} follows from \ef{CL-5.19}, and the bounds just obtained for
$\|\na D_tp \|$ and $\|\na^2 D_tp \|$.
\ef{t-2} follows from \ef{CL-5.29},  \ef{4/24-2} and Lemmas \ref{17lem01} and \ref{17lem2}. \ef{r-4.21} follows from \ef{4.25.2} and Lemma \ref{17lem2}. The bounds for $\|\na^3 \mathcal{T}\|$ and $|\|\na^3 \mathcal{T}\||$ in \ef{step3} follow from \ef{CL-5.28}, \ef{hb4}, \ef{e7}, \ef{r-4.21} and Lemmas \ref{17lem01} and \ref{17lem2}.
It follows   from \ef{dt2divu},   \ef{CL-A.12},  H$\ddot{o}$lder's inequality and Young's inequality that
\begin{align}
  \| D_t^2 {\rm div} u - \mathcal{T} \Delta  D_t {\rm div} u \|^2 \le  C  (K_1, \| \mathcal{T} \|_{L^\iy},\|\na  u\|_{L^\iy},  \|\na \mathcal{T} \|_{L^\iy}, \|\na p \|_{L^\iy})\notag\\
\times \lt(\sum_{i=0}^1 \| \na^i D_t{\rm div} u\|^2   +  \sum_{i=1}^{3}\lt(\|\na^i u\|^2 +  \|\na^i \mathcal{T} \|^2 + \|\na^i p\|^2 \rt)\rt).
    \label{4/24-3}
\end{align}
The bound for $\|\na^2 D_t {\rm div} u\|$ in \ef{step3} follows from \ef{310-1}, \ef{4/24-3}, Lemmas \ref{17lem01} and \ref{17lem2}, and the bounds just obtained for $\|\na^3 u\|$, $\|\na^3 \mathcal{T}\|$ and $\|\na^3 p\|$. The bound for $|\|\na D_t {\rm div} u\||$ in \ef{step3} follows from \ef{CL-5.19} and the bound just obtained for $\|\na^2 D_t {\rm div} u\|$.
The bound for $|\|\na^3 p  \||$ in \ef{step3} follows from \ef{CL-5.28},  \ef{4/24-2}, Lemmas \ref{17lem01} and \ref{17lem2}, and the bounds just obtained for $\|\na^2 D_t {\rm div} u\|$, $\|\na^3 u\|$, $\|\na^3 \mathcal{T}\|$ and $\|\na^3 p\|$.
 It follows from \ef{nadadivu},   \ef{CL-A.12},  H$\ddot{o}$lder's inequality and Young's inequality that for $r= 1,2, 3$,
\be\label{4/24-1}
 \|  \mathcal{T} \na^r \Delta {\rm div} u \|^2    \le  C     \| \na^r D_t{\rm div} u \| ^2 + C(K_1,\|\na  u\|_{L^\iy},
  \|\na \mathcal{T} \|_{L^\iy}) \sum_{i=1}^{r+2}\lt(\|\na^i u\|^2 +  \|\na^i \mathcal{T} \|^2   \rt)    .
\ee
The bound for $|\|\na^2 {\rm div} u\||$ in \ef{step4} follows from \ef{CL-5.28}, \ef{hb3}, \ef{4/24-1}, Lemmas \ref{17lem01} and \ref{17lem2}, and the bounds just obtained for $\|\na^3 u\|$ and $\|\na^3 \mathcal{T}\|$.
The bound for $|\|\Pi \na^3 {\rm div} u\||$ in \ef{step4} follows from \ef{hb4}, \ef{r-4.21}, \ef{ediv2}, and the bound just obtained for $|\|\na^2 {\rm div} u\||$.
The bound for $\|\na^3 {\rm div} u\|$ in \ef{step4} follows from \ef{CL-5.29},  \ef{4/24-1}, Lemmas \ref{17lem01} and \ref{17lem2}, and the bounds just obtained for $|\|\Pi\na^3 {\rm div} u\||$, $\|\na^3 u\|$ and $\|\na^3 \mathcal{T}\|$.
It follows from  \ef{nadadtp},   \ef{CL-A.12},  H$\ddot{o}$lder's inequality and Young's inequality  that for $r=1,2, 3$,
\begin{align}
&  \| \mathcal{T} \na^r \Delta  D_t p \|^2 \le C\| \na^r  D_t^2 {\rm div} u\|^2 +  C(    K_1,   \|  \mathcal{T} \|_{L^\iy}, \|\na  u\|_{L^\iy} , \|\na  \mathcal{T} \|_{L^\iy},     \|\na  p \|_{L^\iy}, \notag\\
& \quad \| D_t  p \|_{L^\iy})\lt(     \sum_{i=1}^{r+2}\lt(\|\na^i u\|^2 +  \|\na^i \mathcal{T} \|^2 + \|\na^i p\|^2 \rt)  + \sum_{i=0}^{r+1} \| \na^i D_t p\|^2 \rt)  \label{4/25-2}.
\end{align}
With \ef{4/25-2}, we can obtain \ef{step4''} by use of a similar way to the derivation of \ef{step4}.
\hfill $\Box$

\begin{lem}\label{17lem4} We have
\be\begin{split}
 & \sum_{i=0}^2 \|\na^i D_t^2 p \|^2   +   |\| \na  D_t^2 p \||^2
  \le      C(K,  K_1, {\rm Vol}\Oa,      \|  \mathcal{T}^{-1} \|_{L^\iy},   \|  \mathcal{T}  \|_{L^\iy},  \|\na  u\|_{L^\iy} ,   \\
   &\quad  \|\na  \mathcal{T} \|_{L^\iy},   \|\na  p \|_{L^\iy},  \| D_t  p \|_{L^\iy}, \| D_t  {\rm div} u \|_{L^\iy},  |\| \na_N p \||_{L^\iy}  )
\sum_{i=1}^4 E_i , \label{step3'}
\end{split}\ee
\be\begin{split}
&\|\na^4 p\|^2 + |\|\overline{\na}^2 \ta \||^2  + \|\na^4  \mathcal{T}\|^2
\le     C(K,  K_1, {\rm Vol}\Oa,      \|  \mathcal{T}^{-1} \|_{L^\iy}, \|  \mathcal{T}  \|_{L^\iy},   \|\na  u\|_{L^\iy} ,  \\
   & \quad \|\na  \mathcal{T} \|_{L^\iy},  \|\na  p \|_{L^\iy},   |\| (\na_N p)^{-1}\||_{L^\iy},  |\| \na_N p \||_{L^\iy}, |\| \na_N \mathcal{T} \||_{L^\iy})\sum_{i=1}^4 E_i,
  \label{17p4}
\end{split}\ee
\be\begin{split}
&\|\na^4  u \|^2 + |\| \na^3    u\||^2 +  |\| \na^4 \mathcal{T}  \||^2 +|\| \na^3   {\rm div} u\||^2 + |\| \Pi \na^4   {\rm div} u\||^2 +  \|\na^4 {\rm div} u  \|^2  \\
&\quad   \le     C(K,  K_1, {\rm Vol}\Oa,      \|  \mathcal{T}^{-1} \|_{L^\iy}, \|  \mathcal{T}  \|_{L^\iy},   \|\na  u\|_{L^\iy} , \|\na  \mathcal{T} \|_{L^\iy},    \|\na  p \|_{L^\iy},   \\
&\quad  |\| (\na_N p)^{-1}\||_{L^\iy},  |\| \na_N p \||_{L^\iy}, |\| \na_N \mathcal{T} \||_{L^\iy}, |\| \na_N  {\rm div} u \||_{L^\iy})
\sum_{i=1}^4 E_i  , \label{17u4}
\end{split}\ee
\be\begin{split}
 &  |\| \na^2 D_t {\rm div} u \||^2  +   |\|\Pi \na^3 D_t {\rm div} u \||^2 + \|\na^3 D_t {\rm div} u \|^2  +  |\|  \na^4 p \||^2  \\
&\quad  \le   C(K,  K_1, {\rm Vol}\Oa,      \|  \mathcal{T}^{-1} \|_{L^\iy}, \|  \mathcal{T}  \|_{L^\iy},  \|\na  u\|_{L^\iy} ,   \|\na  \mathcal{T} \|_{L^\iy},   \|\na  p \|_{L^\iy}, \|D_t  {\rm div} u \|_{L^\iy},   \\
   & \quad  |\| (\na_N p)^{-1}\||_{L^\iy},  |\| \na_N p \||_{L^\iy}, |\| \na_N \mathcal{T} \||_{L^\iy},   |\| \na_N  {\rm div} u \||_{L^\iy}, |\| \na_N D_t {\rm div} u \||_{L^\iy})
\sum_{i=1}^4 E_i ,  \label{17na3dt}
\end{split}\ee
\be\begin{split}
 &  \| \na^2 D_t^2 {\rm div} u \|^2  +   |\|\na  D_t^2 {\rm div} u \||^2 +  |\| \na^3 D_t p \||^2  +   |\|\Pi \na^4 D_t p \||^2 + \|\na^4 D_t p \|^2    \le    C(K, \\
  &\quad K_1,    {\rm Vol}\Oa,      \|  \mathcal{T}^{-1} \|_{L^\iy}, \|  \mathcal{T}  \|_{L^\iy},  \|\na  u\|_{L^\iy} ,   \|\na  \mathcal{T} \|_{L^\iy},   \|\na  p \|_{L^\iy},  \| D_t  p \|_{L^\iy},   \| D_t {\rm div} u \|_{L^\iy},  \\
   &\quad |\| (\na_N p)^{-1}\||_{L^\iy},  |\| \na_N p \||_{L^\iy}, |\| \na_N \mathcal{T} \||_{L^\iy},|\| \na_N  {\rm div} u \||_{L^\iy},    |\| \na_N D_t p \||_{L^\iy})
\sum_{i=1}^4 E_i,  \label{17na4pt}
\end{split}\ee
\be\begin{split}
& |\| \na^2 D_t^2 p    \||^2 +  |\| \Pi \na^3 D_t^2 p    \||^2 +\| \na^3 D_t^2 p \|^2
 \le   C(K,  K_1, {\rm Vol}\Oa,      \|  \mathcal{T}^{-1} \|_{L^\iy}, \|  \mathcal{T}  \|_{L^\iy},  \\
  & \quad \|\na  u\|_{L^\iy} , \|\na  \mathcal{T} \|_{L^\iy},    \|\na  p \|_{L^\iy}, \| D_t  p \|_{L^\iy}, \| D_t^2  p \|_{L^\iy}, \| D_t {\rm div} u \|_{L^\iy},    |\| (\na_N p)^{-1}\||_{L^\iy}, \\
&\quad |\| \na_N p \||_{L^\iy},  |\| \na_N \mathcal{T} \||_{L^\iy}, |\| \na_N  {\rm div} u \||_{L^\iy},|\| \na_N  D_t  p\||_{L^\iy},|\| \na_N  D_t^2 p\||_{L^\iy})
\sum_{i=1}^4 E_i  , \label{17d2tp}
\end{split}\ee
\be\begin{split}
& |\| \na^2 D_t^2 {\rm div} u    \||^2 +  |\| \Pi \na^3 D_t^2 {\rm div} u   \||^2 +\| \na^3 D_t^2 {\rm div} u  \|^2
 \le     C(K,  K_1, {\rm Vol}\Oa,      \|  \mathcal{T}^{-1} \|_{L^\iy}, \\
& \quad \|  \mathcal{T}  \|_{L^\iy},  \|\na  u\|_{L^\iy} ,   \|\na  \mathcal{T} \|_{L^\iy},   \|\na  p \|_{L^\iy},  \| D_t  p \|_{L^\iy}, \| D_t {\rm div} u \|_{L^\iy},    |\| (\na_N p)^{-1}\||_{L^\iy},\\
   & \quad  |\| \na_N p \||_{L^\iy}, |\| \na_N \mathcal{T} \||_{L^\iy}, |\| \na_N  {\rm div} u \||_{L^\iy},|\| \na_N D_t {\rm div} u \||_{L^\iy},   \\
   & \quad |\| \na_N D_t p \||_{L^\iy},|\| \na_N D_t^2 {\rm div} u \||_{L^\iy},\sum_{i=1}^4 E_i)
\sum_{i=1}^4 E_i  . \label{17na3d2tdivu}
\end{split}\ee
\end{lem}

{\em Proof}. It follows from  \ef{dadt2p}, \ef{CL-A.12},  H$\ddot{o}$lder's inequality and Young's inequality that
  \begin{align}
& \|\mathcal{T} \Delta D_t^2 p + (\na \mathcal{T})\cdot \na D_t^2 p\|^2 \notag\\
 \le & C \|   D_t^3 {\rm div} u\|^2  + C(K_1,  \|  \mathcal{T} \|_{L^\iy},    \|\na  u\|_{L^\iy} , \|\na  \mathcal{T} \|_{L^\iy},   \|\na  p \|_{L^\iy},  \| D_t  p \|_{L^\iy}, \| D_t {\rm div} u \|_{L^\iy})  \notag \\
 & \times\lt(     \| \na  D_t{\rm div} u\|^2 + \|\na^2 \mathcal{T} \|^2      + \sum_{i=0}^2 \| \na^i D_t p\|^2  + \sum_{i=1}^{3}\lt(\|\na^i u\|^2 +   \|\na^i p\|^2 \rt)     \rt).\label{4/25-3}
\end{align}
With \ef{4/25-3}, we can obtain  \ef{step3'} by use of a similar way to the derivation of \ef{4/25-1}.
In a similar way to deriving \ef{t-1} and \ef{t-2}, we can obtain, respectively, the bounds for $\|\na^4 u\|$ and $|\|\na^3 u\||$ in \ef{17u4}, and the bound for $\|\na^4 p \| $ in \ef{17p4}. The bound for $|\|\overline{\na}^2 \ta \|| $ in \ef{17p4} follows from \ef{hb2}, \ef{rt1} and Lemmas \ref{17lem2} and \ref{17lem3}.
The bound for $\|\na^4  \mathcal{T} \| $ in \ef{17p4} follows from
\ef{CL-5.29}, \ef{hb1},  the bound just obtained for $|\|\overline{\na}^2 \ta \|| $, and Lemmas \ref{17lem01}-\ref{17lem3}.
The bound for $|\|\na^4  \mathcal{T} \|| $ in \ef{17u4} follows from
\ef{CL-5.28}, \ef{hb1},  the bound just obtained for $|\|\overline{\na}^2 \ta \|| $, and Lemmas \ref{17lem01}-\ref{17lem3}.
The bound for $|\|\na^3 {\rm div} u\||$ in \ef{17u4} follows from \ef{CL-5.28},   \ef{4/24-1}, Lemmas \ref{17lem01}-\ref{17lem3}, and the bounds just obtained for $\|\na^4 u\|$ and $\|\na^4 \mathcal{T}\|$.
The bound for $|\|\Pi \na^4 {\rm div} u\||$ in \ef{17u4} follows from \ef{hb1}, \ef{r-4.21}, \ef{ediv2}, \ef{step4}, and the bounds just obtained for $|\| \overline{\na}^2 \ta \||$ and $|\|\na^3 {\rm div} u\||$.
The bound for $\|\na^4 {\rm div} u\|$ in \ef{17u4} follows from \ef{CL-5.29},  \ef{4/24-1}, Lemmas \ref{17lem01}-\ref{17lem3}, and the bounds just obtained for $|\|\Pi\na^4 {\rm div} u\||$, $\|\na^4 u\|$ and $\|\na^4 \mathcal{T}\|$.
It follows   from \ef{dt2divu},   \ef{CL-A.12},  H$\ddot{o}$lder's inequality and Young's inequality that for $r= 1,2$,
 \begin{align}
&  \| \mathcal{T}\na^r  \Da D_t{\rm div} u\|^2 \le C \|\na^r  D_t^2 {\rm div} u\|^2
 + C( K_1 , \|   \mathcal{T} \|_{L^\iy},   \|\na  u\|_{L^\iy}, \|\na  \mathcal{T} \|_{L^\iy},     \|\na  p \|_{L^\iy},  \notag\\
  & \quad \| D_t {\rm div} u \|_{L^\iy} )   \lt(  \sum_{i=1}^{r+3}\lt(\|\na^i u\|^2 +  \|\na^i \mathcal{T} \|^2 + \|\na^i p\|^2  \rt)   + \sum_{i=0}^{r+1} \| \na^i D_t{\rm div} u\|^2  \rt) . \label{4/24-4}
\end{align}
With \ef{4/24-4}, we can obtain the bounds for $|\| \na^2 D_t {\rm div} u \||$, $|\|\Pi \na^3 D_t {\rm div} u \||^2$, and $\|\na^3 D_t {\rm div} u \|$ in \ef{17na3dt} by use of a similar way to the derivation of \ef{step4}. The bound for $|\|\na^4 p\||$ in \ef{17na3dt} follows from \ef{CL-5.28}, \ef{4/24-2}, Lemmas \ref{17lem01}-\ref{17lem3}, and the bounds just obtained for $\|\na^3 D_t {\rm div} u \|$, $\|\na^4 u\|$, $\|\na^4 \mathcal{T}\|$ and $\|\na^4 p\|$. It follows from \ef{dt3divu}, \ef{CL-A.12},  H$\ddot{o}$lder's inequality and Young's inequality that
\begin{align}
&\|D_t^3 {\rm div} u- \mathcal{T} \Delta  D_t^2 {\rm div} u \|^2 \le
C( K_1, \|  \mathcal{T} \|_{L^\iy}, \|\na  u\|_{L^\iy} , \|\na  \mathcal{T} \|_{L^\iy},     \|\na  p \|_{L^\iy}, \| D_t  p \|_{L^\iy}, \notag\\
&\qquad\| D_t {\rm div} u \|_{L^\iy})\lt(\sum_{i=0}^1 \| \na^i D_t^2 {\rm div} u\|^2  + \sum_{i=1}^{4}\lt(\|\na^i u\|^2 +  \|\na^i \mathcal{T} \|^2 + \|\na^i p\|^2 \rt) \rt. \notag\\
 & \lt. \qquad +  \sum_{i=0}^2 \| \na^i D_t{\rm div} u\|^2 + \sum_{i=0}^3 \| \na^i D_t p\|^2     \rt) . \label{4/9-3}
\end{align}
The bound for $\| \na^2 D_t^2 {\rm div} u \| $ in \ef{17na4pt} follows from \ef{310-1}, \ef{4/9-3}, Lemmas \ref{17lem01}-\ref{17lem3}, and the bounds just obtained for $\|\na^4 u\|$, $\|\na^4 p\|$ and $\|\na^4 \mathcal{T}\|$.
The bound for $|\|\na  D_t^2 {\rm div} u \||$ in \ef{17na4pt}  follows from \ef{CL-5.19}, and the bound just obtained for $\| \na^2 D_t^2 {\rm div} u \| $.
In a similar way to deriving the bounds for $|\| \na^3   {\rm div} u\||$, $|\| \Pi \na^4   {\rm div} u\||$ and $\|\na^4 {\rm div} u  \|$ in \ef{17u4}, we can obtain the bounds  for $|\| \na^3 D_t p \||$,   $|\|\Pi \na^4 D_t p \||$ and $\|\na^4 D_t p \|$ in \ef{17na4pt} by use of \ef{4/25-2} and the bound just obtained for
$\| \na^2 D_t^2 {\rm div} u \| $.
It follows  from \ef{174/5},   \ef{CL-A.12},  H$\ddot{o}$lder's inequality and Young's inequality that for $r= 1,2$,
\begin{align}
&  \| \mathcal{T} \na^r \Delta  D_t^2 p \|^2 \le C\| \na^r  D_t^3 {\rm div} u\|^2 +  C(    K_1,   \|  \mathcal{T} \|_{L^\iy}, \|\na  u\|_{L^\iy} , \|\na  \mathcal{T} \|_{L^\iy},     \|\na  p \|_{L^\iy}, \notag\\
& \quad \| D_t  {\rm div}u \|_{L^\iy},\| D_t  p \|_{L^\iy},\| D_t^2  p \|_{L^\iy})\lt(     \sum_{i=1}^{r+2}\lt(\|\na^i u\|^2 +  \|\na^i \mathcal{T} \|^2 + \|\na^i p\|^2 \rt)  \rt.\notag\\
&\quad  \lt. + \sum_{i=0}^{r+2} \| \na^i D_t p\|^2 + \sum_{i=0}^{r+1} \lt(\| \na^i D_t^2 p\|^2 + \| \na^i D_t{\rm div} u\|^2\rt) + \|\na^{r+3} u\|^2   + \|\na^{r+3}  p\|^2  \rt)  \label{174/10}.
\end{align}
With \ef{174/10}, we can obtain   \ef{17d2tp} in a similar way to deriving \ef{step4}.
It follows from \ef{eh1},   \ef{CL-A.12},  H$\ddot{o}$lder's inequality and Young's inequality that for $r= 1,2$,
\begin{align}
&  \| \mathcal{T}\na^r  \Da D_t^2 {\rm div} u\|^2 \le C \|\na^r  D_t^3 {\rm div} u\|^2
 + C( K_1 , \|   \mathcal{T} \|_{L^\iy},   \|\na  u\|_{L^\iy}, \|\na  \mathcal{T} \|_{L^\iy},     \|\na  p \|_{L^\iy},  \|\na^2  u\|_{L^\iy},  \notag\\
  &   \|\na^2  \mathcal{T} \|_{L^\iy},    \|\na^2  p \|_{L^\iy}, \| D_t {\rm div} u \|_{L^\iy}, \| D_t^2 {\rm div} u \|_{L^\iy}, \|\na D_t p \|_{L^\iy}) \lt( \sum_{i=0}^{r+1} \| \na^i D_t^2 {\rm div} u\|^2 \rt. \notag\\
    &   \lt.  + \sum_{i=0}^{r+2} \| \na^i D_t{\rm div} u\|^2 + \sum_{i=1}^{r+2} \| \na^i D_t p \|^2+ \sum_{i=1}^{r+3}\lt(\|\na^i u\|^2 +  \|\na^i \mathcal{T} \|^2 + \|\na^i p\|^2  \rt)   \rt) . \label{174/11}
\end{align}
Due to \ef{CL-A.15}, we have
\begin{align}
&     \|\na^2  u\|_{L^\iy}^2 + \|\na^2  \mathcal{T} \|_{L^\iy}^2+ \|\na^2  p \|_{L^\iy}^2 + \| D_t^2 {\rm div} u \|_{L^\iy}^2 + \|\na D_t p \|_{L^\iy}^2  \notag\\
 \le & C(K_1) \sum_{i=0}^2  \lt( \|\na^{2+i}  u\|^2 + \|\na^{2+i}  \mathcal{T} \|^2+ \|\na^{2+i}  p \|^2 + \|\na^i D_t^2 {\rm div} u \|^2 + \|\na^{i+1} D_t p \|^2 \rt) \notag \\
 \le & C(K,  K_1, {\rm Vol}\Oa,      \|  \mathcal{T}^{-1} \|_{L^\iy}, \|  \mathcal{T}  \|_{L^\iy},  \|\na  u\|_{L^\iy} ,   \|\na  \mathcal{T} \|_{L^\iy},   \|\na  p \|_{L^\iy},  \| D_t  p \|_{L^\iy}, \| D_t {\rm div} u \|_{L^\iy}, \notag\\
   &  |\| (\na_N p)^{-1}\||_{L^\iy},  |\| \na_N p \||_{L^\iy}, |\| \na_N \mathcal{T} \||_{L^\iy}, |\| \na_N  {\rm div} u \||_{L^\iy}, |\| \na_N D_t p \||_{L^\iy})
\sum_{i=1}^4 E_i. \label{174/11-1}
\end{align}
With \ef{174/11} and \ef{174/11-1}, we can obtain \ef{17na3d2tdivu} in a similar way to deriving \ef{step4}.
\hfill $\Box$

\begin{lem}\label{17lem5} We have
\be\begin{split}
&\|\na^5  u \|^2 + |\| \na^4   u\||^2
 \le     C(K,  K_1, {\rm Vol}\Oa,      \|  \mathcal{T}^{-1} \|_{L^\iy}, \|  \mathcal{T}  \|_{L^\iy},   \|\na  u\|_{L^\iy} , \|\na  \mathcal{T} \|_{L^\iy},    \\
&\quad \|\na  p \|_{L^\iy},   |\| (\na_N p)^{-1}\||_{L^\iy},  |\| \na_N p \||_{L^\iy}, |\| \na_N \mathcal{T} \||_{L^\iy}, |\| \na_N  {\rm div} u \||_{L^\iy})
\sum_{i=1}^5 E_i  , \label{17u5}
\end{split}\ee
\be\begin{split}
 &   \|\na^5 p \|^2
  \le     C(K,  K_1, {\rm Vol}\Oa,      \|  \mathcal{T}^{-1} \|_{L^\iy}, \|  \mathcal{T}  \|_{L^\iy},  \|\na  u\|_{L^\iy} ,   \|\na  \mathcal{T} \|_{L^\iy},   \|\na  p \|_{L^\iy}, \|D_t  {\rm div} u \|_{L^\iy},  \\
   &\quad  |\| (\na_N p)^{-1}\||_{L^\iy},  |\| \na_N p \||_{L^\iy}, |\| \na_N \mathcal{T} \||_{L^\iy},   |\| \na_N  {\rm div} u \||_{L^\iy}, |\| \na_N D_t {\rm div} u \||_{L^\iy})
\sum_{i=1}^5 E_i ,  \label{17p5}
\end{split}\ee
\be\begin{split}
 & \sum_{i=0}^2 \|\na^i D_t^3 p \|^2   +   |\| \na  D_t^3 p \||^2
  \le      C(K,  K_1, {\rm Vol}\Oa,      \|  \mathcal{T}^{-1} \|_{L^\iy},   \|  \mathcal{T}  \|_{L^\iy},  \|\na  u\|_{L^\iy} ,   \\
   & \quad \|\na  \mathcal{T} \|_{L^\iy},  \|\na  p \|_{L^\iy},  \| D_t  p \|_{L^\iy}, \| D_t  {\rm div} u \|_{L^\iy}, \| D_t^2  p \|_{L^\iy}, \| D_t^2  {\rm div} u \|_{L^\iy},   |\| (\na_N p)^{-1}\||_{L^\iy},  \\
   & \quad  |\| \na_N p \||_{L^\iy},  |\| \na_N \mathcal{T} \||_{L^\iy}, |\| \na_N D_t p \||_{L^\iy},   |\| \na_N  {\rm div} u \||_{L^\iy}  )
\sum_{i=1}^5 E_i , \label{174/14}
\end{split}\ee
\be\begin{split}
&  |\|\overline{\na}^3 \ta \||^2 +    \|\na^5 \mathcal{T}\| ^2 +   |\| \na^4   {\rm div} u\||^2 + |\| \Pi \na^5   {\rm div} u\||^2 +  \|\na^5 {\rm div} u  \|^2
\le    C(K,  K_1, {\rm Vol}\Oa,  \\
     & \quad \|  \mathcal{T}^{-1} \|_{L^\iy}, \|  \mathcal{T}  \|_{L^\iy},  \|\na  u\|_{L^\iy} ,   \|\na  \mathcal{T} \|_{L^\iy},   \|\na  p \|_{L^\iy}, \|D_t  {\rm div} u \|_{L^\iy},   |\| (\na_N p)^{-1}\||_{L^\iy},  \\
   &\quad  |\| \na_N p \||_{L^\iy}, |\| \na_N \mathcal{T} \||_{L^\iy},   |\| \na_N  {\rm div} u \||_{L^\iy}, |\| \na_N D_t {\rm div} u \||_{L^\iy},\sum_{i=1}^4 E_i)
\sum_{i=1}^5 E_i ,
  \label{17ta3}
\end{split}\ee
\be\begin{split}
 &  |\| \na^3 D_t {\rm div} u \||^2  +   |\|\Pi \na^4 D_t {\rm div} u \||^2 + \|\na^4 D_t {\rm div} u \|^2  +  |\|  \na^5 p \||^2  \\
&\quad  \le    C(K,  K_1, {\rm Vol}\Oa,      \|  \mathcal{T}^{-1} \|_{L^\iy}, \|  \mathcal{T}  \|_{L^\iy},  \|\na  u\|_{L^\iy} ,   \|\na  \mathcal{T} \|_{L^\iy},   \|\na  p \|_{L^\iy},  \| D_t  p \|_{L^\iy},  \\
   &\quad \| D_t {\rm div} u \|_{L^\iy},  |\| (\na_N p)^{-1}\||_{L^\iy},  |\| \na_N p \||_{L^\iy}, |\| \na_N \mathcal{T} \||_{L^\iy}, |\| \na_N  {\rm div} u \||_{L^\iy}, \\
    & \quad|\| \na_N D_t {\rm div} u \||_{L^\iy},  |\| \na_N D_t p \||_{L^\iy},\sum_{i=1}^4 E_i)
\sum_{i=1}^5 E_i ,  \label{174/3}
\end{split}\ee
\be\begin{split}
 &   \|\na^2 D_t^3 {\rm div} u \|^2
  \le     C(K,  K_1, {\rm Vol}\Oa,      \|  \mathcal{T}^{-1} \|_{L^\iy}, \|  \mathcal{T}  \|_{L^\iy},  \|\na  u\|_{L^\iy} ,   \|\na  \mathcal{T} \|_{L^\iy},   \|\na  p \|_{L^\iy},  \\
   & \quad \| D_t  p \|_{L^\iy}, \| D_t {\rm div} u \|_{L^\iy}, |\| (\na_N p)^{-1}\||_{L^\iy},    |\| \na_N p \||_{L^\iy}, |\| \na_N \mathcal{T} \||_{L^\iy}, |\| \na_N  {\rm div} u \||_{L^\iy}, \\
   &\quad   |\| \na_N D_t p \||_{L^\iy},\sum_{i=1}^4 E_i)
\sum_{i=1}^5 E_i ,  \label{174/3'}
\end{split}\ee
\be\begin{split}
 &|\| \na^4 D_t p \||^2  +   |\|\Pi \na^5 D_t p \||^2 + \|\na^5 D_t p \|^2  \\
& \quad \le   C(K,  K_1, {\rm Vol}\Oa,      \|  \mathcal{T}^{-1} \|_{L^\iy}, \|  \mathcal{T}  \|_{L^\iy},  \|\na  u\|_{L^\iy} ,   \|\na  \mathcal{T} \|_{L^\iy},   \|\na  p \|_{L^\iy},
   \| D_t  p \|_{L^\iy}, \\
   &\quad \| D_t {\rm div} u \|_{L^\iy},   |\| (\na_N p)^{-1}\||_{L^\iy},  |\| \na_N p \||_{L^\iy}, |\| \na_N \mathcal{T} \||_{L^\iy}, |\| \na_N  {\rm div} u \||_{L^\iy}, \\
   &\quad |\| \na_N D_t {\rm div} u \||_{L^\iy},  |\| \na_N D_t p \||_{L^\iy},|\| \na_N D_t^2 {\rm div} u \||_{L^\iy},\sum_{i=1}^4 E_i)
\sum_{i=1}^5 E_i, \label{174/11'}
\end{split}\ee
\be\begin{split}
& |\| \na^3 D_t^2 p    \||^2 +  |\| \Pi \na^4 D_t^2 p    \||^2 +\| \na^4 D_t^2 p \|^2  \\
&\quad  \le    C(K,  K_1, {\rm Vol}\Oa,      \|  \mathcal{T}^{-1} \|_{L^\iy}, \|  \mathcal{T}  \|_{L^\iy},   \|\na  u\|_{L^\iy} , \|\na  \mathcal{T} \|_{L^\iy},    \|\na  p \|_{L^\iy}, \| D_t  p \|_{L^\iy},  \\
&\quad \| D_t^2  p \|_{L^\iy},  \| D_t {\rm div} u \|_{L^\iy},   |\| (\na_N p)^{-1}\||_{L^\iy},  |\| \na_N p \||_{L^\iy}, |\| \na_N \mathcal{T} \||_{L^\iy}, |\| \na_N  {\rm div} u \||_{L^\iy},
 \\
&\quad |\| \na_N D_t {\rm div} u \||_{L^\iy}, |\| \na_N  D_t  p\||_{L^\iy},|\| \na_N  D_t^2 p\||_{L^\iy})
\sum_{i=1}^5 E_i , \label{17na3d2tp}
\end{split}\ee
\be\begin{split}
 & |\| \na^3 D_t^2 {\rm div} u    \||^2 +  |\| \Pi \na^4 D_t^2 {\rm div} u   \||^2 +\| \na^4 D_t^2 {\rm div} u  \|^2  \\
&\quad \le   C(K,  K_1, {\rm Vol}\Oa,      \|  \mathcal{T}^{-1} \|_{L^\iy}, \|  \mathcal{T}  \|_{L^\iy},  \|\na  u\|_{L^\iy} ,   \|\na  \mathcal{T} \|_{L^\iy},   \|\na  p \|_{L^\iy},  \| D_t  p \|_{L^\iy}, \\
   &\quad  \| D_t {\rm div} u \|_{L^\iy},   |\| (\na_N p)^{-1}\||_{L^\iy},  |\| \na_N p \||_{L^\iy}, |\| \na_N \mathcal{T} \||_{L^\iy}, |\| \na_N  {\rm div} u \||_{L^\iy}, \\
   &\quad|\| \na_N D_t {\rm div} u \||_{L^\iy},  |\| \na_N D_t p \||_{L^\iy},|\| \na_N D_t^2 {\rm div} u \||_{L^\iy},\sum_{i=1}^4 E_i)
\sum_{i=1}^5 E_i .\label{174/12}
\end{split}\ee
\end{lem}

{\em Proof}.
In a similar way to deriving \ef{t-1} and \ef{t-2}, we can obtain, respectively, the bounds for $\|\na^5 u\|$ and $|\|\na^4 u\||$ in \ef{17u5}, and the bound for $\|\na^5 p \| $ in \ef{17p5}.
It follows from  \ef{174/4}, \ef{CL-A.12},  H$\ddot{o}$lder's inequality and Young's inequality that
  \begin{align}
& \|\mathcal{T} \Delta D_t^3 p + (\na \mathcal{T})\cdot \na D_t^3 p\|^2
 \le   C \|   D_t^4 {\rm div} u\|^2  + C(K_1,  \|  \mathcal{T} \|_{L^\iy},    \|\na  u\|_{L^\iy} , \|\na  \mathcal{T} \|_{L^\iy},   \|\na  p \|_{L^\iy}, \notag\\
  &\quad \| D_t  p \|_{L^\iy}, \| D_t {\rm div} u \|_{L^\iy},  \| D_t^2  p \|_{L^\iy}, \| D_t^2 {\rm div} u \|_{L^\iy})   (     \| \na  D_t^2 {\rm div} u\|^2 + \sum_{i=0}^2 (\| \na^i  D_t  {\rm div} u\|^2 \notag\\
  &\quad     + \| \na^i D_t^2 p\|^2 )  + \sum_{i=0}^3 \| \na^i D_t  p\|^2  + \sum_{i=1}^3 \|\na^i \mathcal{T} \|^2   + \sum_{i=1}^{4} (\|\na^i u\|^2 +   \|\na^i p\|^2)      ).\label{17414}
\end{align}
With \ef{17414}, we can obtain  \ef{174/14} by use of a similar way to the derivation of \ef{4/25-1}.
The bound for $|\|\overline{\na}^3 \ta \|| $ in \ef{17ta3} follows from \ef{17hb2}, \ef{rt1} and Lemmas \ref{17lem2}-\ref{17lem4}.
The bound for $|\|\na^4 {\rm div} u\||$ in \ef{17ta3} follows from \ef{CL-5.28},   \ef{4/24-1}, Lemmas \ref{17lem01}-\ref{17lem4}, and the bounds just obtained for $\|\na^5 u\|$ and $\|\na^5\mathcal{T}\|$.
The bound for $|\|\Pi \na^5 {\rm div} u\||$ in \ef{17ta3} follows from \ef{17hb1}, Lemmas \ref{17lem2}-\ref{17lem4}, and the bounds just obtained for $|\| \overline{\na}^3 \ta \||$ and $|\|\na^4 {\rm div} u\||$.
The bound for $\|\na^5 {\rm div} u\|$ in \ef{17ta3} follows from \ef{CL-5.29},  \ef{4/24-1}, Lemmas \ref{17lem01}-\ref{17lem4}, and the bounds just obtained for $|\|\Pi\na^5 {\rm div} u\||$, $\|\na^5 u\|$ and $\|\na^5 \mathcal{T}\|$.
In a similar way to deriving the bounds for $|\| \na^3   {\rm div} u\||$, $|\| \Pi \na^4   {\rm div} u\||$ and $\|\na^4 {\rm div} u  \|$ in \ef{17u4}, we can obtain the bounds  for $|\| \na^3 D_t {\rm div} u  \||$,   $|\|\Pi \na^4 D_t {\rm div} u  \||$ and $\|\na^4 D_t {\rm div} u  \|$ in \ef{174/3} by use of \ef{4/24-4}.
The bound for $|\|\na^5 p\||$ in \ef{174/3} follows from \ef{CL-5.28}, \ef{4/24-2}, Lemmas \ref{17lem01}-\ref{17lem4}, and the bounds just obtained for $\|\na^4 D_t {\rm div} u \|$, $\|\na^5 u\|$, $\|\na^5 \mathcal{T}\|$ and $\|\na^5 p\|$.
It follows from \ef{nt2}, \ef{CL-A.15}, and Lemmas \ref{17lem01}-\ref{17lem4} that
\begin{align}
&\|D_t^4 {\rm div} u- \mathcal{T} \Delta  D_t^3 {\rm div} u \|^2
\le     C(K,  K_1, {\rm Vol}\Oa,      \|  \mathcal{T}^{-1} \|_{L^\iy}, \|  \mathcal{T}  \|_{L^\iy},  \|\na  u\|_{L^\iy} ,   \|\na  \mathcal{T} \|_{L^\iy},  \notag \\
&\quad    \|\na  p \|_{L^\iy},  \| D_t  p \|_{L^\iy}, \| D_t {\rm div} u \|_{L^\iy}, |\| (\na_N p)^{-1}\||_{L^\iy},  |\| \na_N p \||_{L^\iy},   |\| \na_N \mathcal{T} \||_{L^\iy}, \notag\\
   &\quad |\| \na_N  {\rm div} u \||_{L^\iy},    |\| \na_N D_t p \||_{L^\iy},\sum_{i=1}^4 E_i)
\sum_{i=1}^4 E_i,\label{17nt2}
\end{align}
which, together with \ef{310-1} and \ef{4/10-7}, gives  \ef{174/3'}.
In a similar way to deriving the bounds for $|\| \na^4   {\rm div} u\||$, $|\| \Pi \na^5   {\rm div} u\||$ and $\|\na^5 {\rm div} u  \|$ in \ef{17ta3}, we can obtain \ef{174/11'}   by use of  \ef{4/25-2}.
In a similar way to deriving the bounds for $|\| \na^3   {\rm div} u\||$, $|\| \Pi \na^4   {\rm div} u\||$ and $\|\na^4 {\rm div} u  \|$ in \ef{17u4}, we can obtain, respectively,  \ef{17na3d2tp}  and \ef{174/12}   by use of \ef{174/10} and \ef{174/11}.
\hfill$\Box$

\subsection{Proof of Proposition \ref{prop1}}

After having seen in Section \ref{17elliptic} what a bound for the energy implies, we can now prove \ef{17prop1}.

It follows from \ef{4/12-9'}, and Lemmas \ref{17lem01} and \ref{17lem2} that
\be\label{17et0}
   \frac{d}{dt} E_0
 \le  C(  {\rm Vol}\Oa, \|\mathcal{T}^{-1}\|_{L^\iy}, \|\mathcal{T} \|_{L^\iy},\|\nabla u\|_{L^\iy} ) (E_1+ E_2) .
\ee
It follows from \ef{17div1}, \ef{4/12-9},  and Lemmas \ref{17lem01} and \ref{17lem2} that
\be\label{17et1}
   \frac{d}{dt} E_1
 \le  C( K,  {\rm Vol}\Oa, \|\mathcal{T}^{-1}\|_{L^\iy}, \|\mathcal{T} \|_{L^\iy},\|\nabla u\|_{L^\iy},  \|\nabla \mathcal{T}\|_{L^\iy} ) (E_1+ E_2).
\ee
It follows from \ef{17div2}, \ef{hfin4.58},  and Lemmas \ref{17lem01} and \ref{17lem2} that
\be\label{17et2}\begin{split}
&  \frac{d}{dt} E_2
 \le  C( K, K_1,{\rm Vol}\Oa, \|\mathcal{T}^{-1}\|_{L^\iy}, \|\mathcal{T} \|_{L^\iy},\|\nabla u\|_{L^\iy}, \|\nabla p\|_{L^\iy},\|\nabla \mathcal{T}\|_{L^\iy},\|\nabla {\rm div} u\|_{L^\iy}, \\
&\|\na^2 \mathcal{T} \|_{L^\iy},|\| (\na_N p)^{-1} \||_{L^\iy}, |\|\nabla_N p\||_{L^\iy},|\|\nabla_N D_t p \||_{L^\iy}, |\|\nabla u\||_{L^\iy}, |\|\nabla^2 p\||_{L^\iy}) (E_1+ E_2) .
\end{split}\ee
It follows from \ef{17div3}, \ef{hfin4.58},  and Lemmas \ref{17lem01}-\ref{17lem3} that
\be\label{17et3}\begin{split}
&   \frac{d}{dt} E_3
 \le   C( K, K_1,{\rm Vol}\Oa, \|\mathcal{T}^{-1}\|_{L^\iy}, \|\mathcal{T} \|_{L^\iy},\|\nabla u\|_{L^\iy}, \|\nabla p\|_{L^\iy},\|\nabla \mathcal{T}\|_{L^\iy},\\
&\|\nabla {\rm div} u\|_{L^\iy}, \|\na^2 \mathcal{T} \|_{L^\iy},\|D_t{\rm div} u\|_{L^\iy},\|D_t p\|_{L^\iy},|\| (\na_N p)^{-1} \||_{L^\iy}, |\|\nabla_N p\||_{L^\iy},\\
& |\|\nabla_N  \mathcal{T}\||_{L^\iy},  |\|\nabla_N  {\rm div} u \||_{L^\iy}, |\|\nabla_N D_t p \||_{L^\iy}, |\|\nabla u\||_{L^\iy}, |\|\nabla^2 p\||_{L^\iy}) \sum_{i=1}^3 E_i.
\end{split}\ee
It follows from \ef{CL-A.15} that
$$\|\na D_t p \|_{L^\iy}\le C(K_1)\sum_{i=1}^3 \|\na^i D_t p\|,$$
which, together with \ef{17div4}, \ef{hfin4.58},  and Lemmas \ref{17lem01}-\ref{17lem4}, implies that
\be\label{17et4}\begin{split}
&  \frac{d}{dt} E_4
 \le  C( K, K_1,{\rm Vol}\Oa, \|\mathcal{T}^{-1}\|_{L^\iy}, \|\mathcal{T} \|_{L^\iy},\|\nabla u\|_{L^\iy}, \|\nabla p\|_{L^\iy},\|\nabla \mathcal{T}\|_{L^\iy},\|\nabla {\rm div} u\|_{L^\iy}, \\
&\|\na^2 \mathcal{T} \|_{L^\iy},\|D_t{\rm div} u\|_{L^\iy},\|D_t p\|_{L^\iy},|\| (\na_N p)^{-1} \||_{L^\iy}, |\|\nabla_N p\||_{L^\iy},|\|\nabla_N  \mathcal{T}\||_{L^\iy}, \\
&|\|\nabla_N  {\rm div} u \||_{L^\iy},  |\|\nabla_N D_t p \||_{L^\iy}, |\|\nabla_N D_t {\rm div} u \||_{L^\iy}, |\|\nabla u\||_{L^\iy}, |\|\nabla^2 p\||_{L^\iy}, \sum_{i=1}^3 E_i) \sum_{i=1}^4 E_i.
\end{split}\ee
It follows from \ef{CL-A.15} that
\bee\begin{split}
 \|D_t^2 p\|_{L^\iy}+\|D_t^2 {\rm div} u\|_{L^\iy}+ \|\na  D_t p\|_{L^\iy} + \|\na^2 D_t p\|_{L^\iy}+\|\na^2 u\|_{L^\iy}+ \|\na^2 p\|_{L^\iy} \\
+\|\na^2 {\rm div} u\|_{L^\iy}
\le   C(K_1)\sum_{j=0}^2 ( \|\na^j D_t^2 p\| +\|\na^j D_t^2 {\rm div} u\| + \|\na^{1+ j} D_t p\|\\
 + \|\na^{2+ j} D_t p\| +\|\na^{2+j} u\|  + \|\na^{2+j} p\|  +\|\na^{2+j} {\rm div} u\|  ),
\end{split}\eee
which, together with \ef{17div5}, \ef{hfin4.58},  and Lemmas \ref{17lem01}-\ref{17lem5}, implies that
\be\label{17et5}\begin{split}
&  \frac{d}{dt} E_5
 \le  C( K, K_1,{\rm Vol}\Oa, \|\mathcal{T}^{-1}\|_{L^\iy}, \|\mathcal{T} \|_{L^\iy},\|\nabla u\|_{L^\iy}, \|\nabla p\|_{L^\iy},\|\nabla \mathcal{T}\|_{L^\iy},\|\nabla {\rm div} u\|_{L^\iy}, \\
&\|\na^2 \mathcal{T} \|_{L^\iy},\|D_t{\rm div} u\|_{L^\iy},\|D_t p\|_{L^\iy},|\| (\na_N p)^{-1} \||_{L^\iy}, |\|\nabla_N p\||_{L^\iy},|\|\nabla_N  \mathcal{T}\||_{L^\iy},  \\
&|\|\nabla_N  {\rm div} u \||_{L^\iy},  |\|\nabla_N D_t p \||_{L^\iy}, |\|\nabla_N D_t {\rm div} u \||_{L^\iy},|\|\nabla_N D_t^2 p \||_{L^\iy}, |\|\nabla_N D_t^2 {\rm div} u \||_{L^\iy}, \\
&  |\|\nabla u\||_{L^\iy}, |\|\nabla^2 p\||_{L^\iy}, \sum_{i=1}^4 E_i) \sum_{i=1}^5 E_i .
\end{split}\ee

It is produced by substituting \ef{17/5/19}, \ef{w4.10}, \ef{bdta} into \ef{17et0}-\ef{17et5} that  there are continuous functions $C_r$  $(0\le r\le n+2)$ such that
\bee
\frac{d}{dt} E_0(t) \le C_0\lt( \overline{V},    M, \underline{\mathcal{T}}^{-1}, \overline{\mathcal{T}}  \rt) \lt(E_1(t)+E_2(t)\rt) ,
\eee
\bee
\frac{d}{dt}  E_1(t) \le C_1 \lt( \overline{V},K,   M, \underline{\mathcal{T}}^{-1}, \overline{\mathcal{T}}   \rt)\lt(E_1(t)+E_2(t)\rt),
\eee
\bee
 \frac{d}{dt}  E_2(t) \le C_2 \lt( \overline{V},K,  \ea_b^{-1},L, M, \widetilde{M}, \underline{\mathcal{T}}^{-1}, \overline{\mathcal{T}}    \rt)\lt(E_1(t)+E_2(t)\rt),
 \eee
 \bee
 \frac{d}{dt}  E_r(t) \le C_r \lt( \overline{V},K,  \ea_b^{-1},L, M, \widetilde{M}, \underline{\mathcal{T}}^{-1}, \overline{\mathcal{T}} , \sum_{s=1}^{r-1} E_s(t)  \rt)\sum_{s=1}^r E_s(t), \ \   3 \le r\le n+2.
\eee
This concludes the proof of Proposition \ref{prop1}.

\subsection{Proof of Proposition \ref{prop2}}
Let us now show how Proposition \ref{prop2} follows.

\begin{lem}\label{17lem6}
 Let $n=2,3$. Then there are continuous functions $T_{n1}>0$ such that
\be\label{17/5/20-4}\begin{split}
&\sum_{s=0}^{r} E_s(t)\le 2 \sum_{s=0}^{r}   E_s(0),      \ \ 2 \le r \le n+2,  \ \  \ \  2^{-1} {\rm Vol}\mathscr{D}_0     \le  {\rm Vol}\Omega (t)     \le  2 {\rm Vol}\mathscr{D}_0, \\   & |\|\ta(t,\cdot)\||_{L^\iy}+ \iota_0^{-1}(t)\le 18 K_0 , \ \  \ \
 -\na_N p(t,  y) \ge 2^{-1}\underline{\ea}_0        \ \ {\rm for}  \ \    y \in  \pl \Omega, \\
&     \|\na p(t,\cdot)\|_{L^\iy} + \|\na u(t,\cdot)\|_{L^\iy} + \|\na \mathcal{T}(t,\cdot)\|_{L^\iy} \le 2 M_0,
\end{split}\ee
for $ t\le T_{n1}({\rm Vol}\mathscr{D}_0 , K_0,  \underline{\ea}_0^{-1},  L, M, \widetilde{M}, \underline{\mathcal{T}}^{-1}, \overline{\mathcal{T}} , M_0, E_0(0),\cdots,E_{n+2}(0) ) $.
\end{lem}

{\em Proof}. The cases $n=2$ and $n=3$ can be shown in the same manner, so we present here only the proof of the case of $n=2$.
Let $n=2$ in the rest of this proof.

It follows from \ef{s3.7} that
$$
\lt|\frac{d}{dt}{\rm Vol}\Omega (t) \rt| = \lt|\frac{d}{dt} \int_{\Omega} d\mu_g \rt|=   \lt|\int_{\Omega} {\rm div} u d\mu_g \rt| \le M {\rm Vol}\Omega (t),
$$
which implies
$
 {\rm Vol}\Omega(0) \exp\{-Mt\}  \le  {\rm Vol}\Omega (t)     \le  {\rm Vol}\Omega(0) \exp\{Mt\}.
$
Thus we  have,  due to the fact ${\rm Vol} \Omega(0)= {\rm Vol}\mathscr{D}_0$, that for  $t \le M^{-1}\ln 2$,
\be\label{bforvol}
2^{-1} {\rm Vol}\mathscr{D}_0     \le  {\rm Vol}\Omega (t)     \le  2 {\rm Vol}\mathscr{D}_0 .
\ee

It follows from \ef{dtnanp}, \ef{CL-3.21}, \ef{CL-7.4}, \ef{CL-7.5} and \ef{CL-7.6} that
\bee
|D_t (-\na_N p)^{-1}|\le   | (-\na_N p)^{-1}|^2(| h_{NN}\na_N p|+|\na_N D_t p|) \le \ea_b^{-2}(M^2+L) \ \ {\rm on} \ \ \pl\Omega,
\eee
which, together with $\pl_\mathcal{N}p=\na_N p$ and  \ef{s118}, implies that for $t \le  \ea_b^2(M^2+ L)^{-1} \underline{\ea}_0^{-1}$,
\be\label{bfornanp}
 -\na_N p(t,  y) \ge 2^{-1}\underline{\ea}_0        \ \ {\rm for}  \ \    y \in  \pl \Omega .
\ee

Let $\ea_1\in (0,1/2]$ be a fixed constant (for example, $\ea_1=1/4$), $\iota_1(0)$ the largest number such that
\be\label{CL-7.55}\begin{split}
&|\mathcal{N}( x(0,y_1))-\mathcal{N}( x(0,y_2))|\le \ea_1/2 ,
\\
& \ \  \ \  \ \  {\rm whenever} \ \
|x(0,y_1) -x(0,y_2)| \le  2\iota_1(0),  \ \
 y_1, y_2 \in \pl\Oa.
\end{split}\ee
Thus  we have from  Lemma \ref{lemkk1}  and \ef{s118} that
\be\label{iota10}
\iota_1(0) \ge 2^{-1}K_0^{-1}  \ea_1 .
\ee
Due to $\pl_t x (t,y)=v(t, x(t,y))$ in $\bar\Omega$,
$|D_t \mathcal{N}|\le 2 |\na u|\le 2 M$ on $\pl\Oa$,  and
\be\label{CL-7.47}
\|v(t,x(t,\cdot))\|_{L^\iy(\bar{\Omega})} \le 2 \|v(0,x(0,\cdot))\|_{L^\iy(\bar{\Omega})}  \ \ {\rm for} \ \  t \le  ({\overline{\mathcal{T}}  M })^{-1}  \|v(0,x(0,\cdot))\|_{L^\iy(\bar{\Omega})} ,
\ee
we  have
\begin{align}
&|x(t,y)-x(0,y)|\le 2^{-1}\iota_1(0) \ \  {\rm for} \ \    y\in \bar\Oa \ \  {\rm and} \ \ t\le T_1   , \label{CL-7.40} \\
&|\mathcal{N}( x(t,   y))-\mathcal{N}( x(0,   y))| \le 4^{-1}\ea_1  \ \ {\rm for} \ \   y\in \pl\Oa  \ \  {\rm and} \ \ t\le 8^{-1}M^{-1} \ea_1, \label{CL-7.39}
\end{align}
where
 $ T_1   = \min \{({\overline{\mathcal{T}}  M })^{-1}  \|v(0,x(0,\cdot))\|_{L^\iy(\bar{\Omega})} ,  \ 4^{-1}\|v(0,x(0,\cdot))\|_{L^\iy(\bar{\Omega})}^{-1} \iota_1(0) \}  $.
Indeed, the bound $|D_t v|= |\mathcal{T} \pl p|= |\mathcal{T} \na p| \le \overline{\mathcal{T}} M $ in $\bar{\Oa}$, which follows from \ef{l2-a}, \ef{bdta} and \ef{CL-7.6},
 has been used to derive \ef{CL-7.47}.
It follows from \ef{CL-7.55}, \ef{CL-7.40} and \ef{CL-7.39} that
for $t \le \min\{T_1, \ 8^{-1}M^{-1} \ea_1\}$,
\be\label{CL-7.56}\begin{split}
&|\mathcal{N}( x(t,y_1))-\mathcal{N}( x(t,y_2))|\le \ea_1  ,
\\
& \ \  \ \  \ \  {\rm whenever} \ \
|x(t,y_1) -x(t,y_2)| \le   \iota_1(0),  \ \
 y_1, y_2 \in \pl\Oa   .
\end{split}\ee
This, together with \ef{iota10},  implies that for  $t \le \min\{T_1, \ 8^{-1}M^{-1} \ea_1\}$,
\be\label{17514-3}
\iota_1(t) \ge \iota_1(0) \ge 2^{-1}K_0^{-1}  \ea_1 .
\ee

It follows from Proposition \ref{prop1} that there are  continuous functions $T_r>0$ such that
\be\label{bforet}
\sum_{s=0}^{r} E_s(t)\le 2 \sum_{s=0}^{r}   E_s(0)  , \ \  2\le r \le 4  ,
\ee
for $  t \le T_r(\overline{V}, K,  \ea_b^{-1},L, M, \widetilde{M}, \underline{\mathcal{T}}^{-1}, \overline{\mathcal{T}} , E_1(0),\cdots,E_{r-1}(0))$.
This, together with \ef{CL-A.15},  \ef{CL-A.8}, Lemmas \ref{17lem2}-\ref{17lem4}, \ef{17/5/19}, \ef{w4.10} and \ef{bdta}, gives that
\be\label{17/5/20-1}\begin{split}
 \|\na D_t p(t,\cdot)\|_{L^\iy}^2
\le & C(K_1)\sum_{i=1}^3\|\na^i D_t p(t,\cdot)\|^2  \\
  \le & C(\overline{V}, K,  \ea_b^{-1},L, M, \widetilde{M}, \underline{\mathcal{T}}^{-1}, \overline{\mathcal{T}})\sum_{s=0}^3E_s(0), \ \  t\le T_3 ,
\end{split}\ee
\be\label{17/5/20-2}\begin{split}
 \|\na^2 p(t,\cdot)\|_{L^\iy}^2 +   |\|\na^2 u(t,\cdot)\||_{L^\iy}^2 \le & C(K_1)\lt(\sum_{i=2}^4 \|\na^i p(t,\cdot)\|^2 + \sum_{i=2}^3 |\|\na^i u(t,\cdot)\||^2\rt)
    \\
  \le & C(\overline{V}, K,  \ea_b^{-1},  M,  \underline{\mathcal{T}}^{-1}, \overline{\mathcal{T}}) \sum_{s=0}^4 E_s(0), \ \  t\le T_4.
\end{split}\ee
Notice that
\bee
 |D_t \na p| = |\na D_t p| \le \|\na D_t  p \|_{L^\iy}  \ \ {\rm in} \ \  \Oa ,  \  \  |D_t \na p| = |\na D_t p| \le L \ \ {\rm on} \ \ \pl\Oa,
 \eee
 \be\label{17514-1}
 |D_t \na u| \le \mathcal{T}|\na^2  p| +|\na \mathcal{T}||\na u| +  |\na u|^2 \le  \overline{\mathcal{T}}\|\na^2  p\|_{L^\iy} + 2M^2 \ \ {\rm in} \  \ \Oa,
 \ee
 \be\label{17514-2}
  | D_t \na \mathcal{T}| \le \mathcal{T}|\na {\rm div} u| + |\na \mathcal{T}||\na u| \le    \overline{\mathcal{T}} M + M^2  \ \ {\rm in} \ \  \Oa,
  \ee
\bee  |D_t \ta| \le |\na^2 u| + C |\ta||\na u| \le  |\|\na^2 u\||_{L^\iy} + C KM  \ \ {\rm on} \ \ \pl\Oa. \eee
Here \ef{17514-1} (or respectively, \ef{17514-2}) follows from \ef{e25} (or respectively, \ef{r-3.28}), \ef{bdta}  and \ef{CL-7.6}. Then we have from \ef{s118}, \ef{s113}, the fact that
$|\pl p|=|\na p|$, $|\pl v|=|\na u|$ and $|\pl \mathcal{T}|=|\na \mathcal{T} |$, \ef{17/5/20-1} and \ef{17/5/20-2} that
there is a continuous function $T_5>0$ such that
\be\label{17/5/20-3}
 \|\na p(t,\cdot)\|_{L^\iy} + \|\na u(t,\cdot)\|_{L^\iy} + \|\na \mathcal{T}(t,\cdot)\|_{L^\iy} + |\| \na p(t,\cdot)\||_{L^\iy} \le 2 M_0, \ee
\be\label{bforta0}
|\|\ta(t,\cdot)\||_{L^\iy}\le 2 K_0 ,
\ee
for  $t \le   T_5(\overline{V}, K,  \ea_b^{-1},L, M, \widetilde{M}, \underline{\mathcal{T}}^{-1}, \overline{\mathcal{T}} , E_0(0),\cdots,E_{4}(0), M_0, K_0)$.
Moreover, we can derive from Lemma \ref{lemkk1}, \ef{bforta0} and \ef{17514-3} that   for  $t \le \min\{T_1, \ 8^{-1}M^{-1} \ea_1, \ T_5\}$,
\be\label{4/30-1}
 \iota_0^{-1}(t) \le \max\{2  \iota_1^{-1}(t), \ 2 K_0\} \le 4 \ea_1^{-1}K_0 .
\ee

Clearly, there is a continuous function  $ T_{6}>0$ such that
\ef{bforvol}, \ef{bfornanp}, \ef{bforet}, \ef{17/5/20-3}, \ef{bforta0} and \ef{4/30-1} hold
for
 $t \le  T_{6}(\overline{V}, K,  \ea_b^{-1},L, M, \widetilde{M}, \underline{\mathcal{T}}^{-1}, \overline{\mathcal{T}} , E_0(0),\cdots,E_{4}(0), \underline{\ea}_0, M_0, K_0,{\rm Vol}\mathscr{D}_0) $, since
$\|v(0,x(0,\cdot))\|_{L^\iy(\bar{\Omega})}^2 \le C({\rm Vol}\mathscr{D}_0, \overline{\mathcal{T}})\sum_{s=0}^2 E_s(0)$ which follows from $|v|=|u|$, \ef{CL-A.15},  \ef{CL-A.8},  Lemmas \ref{17lem01} and \ref{17lem2}, and \ef{bdta}.
So,   \ef{17/5/20-4} holds for $t\le T_{21}$ for some continuous function $$T_{21}({\rm Vol}\mathscr{D}_0 , K_0,  \underline{\ea}_0^{-1},L, M, \widetilde{M}, \underline{\mathcal{T}}^{-1}, \overline{\mathcal{T}} , E_0(0),\cdots,E_{4}(0),  M_0 )>0,$$ by choosing $\overline{V}=4 {\rm Vol}\mathscr{D}_0$, $\ea_b=4^{-1}\underline{\ea}_0$, and $K=2(2+4 \ea_1^{-1})K_0 = 36 K_0$ with $\ea_1$ being $4^{-1}$ in $T_{6}$.
\hfill $\Box$

\begin{lem}
Let $n=2,3$. Then there are continuous functions $T_{n2}>0$ such that
\be\label{17/5/22}\begin{split}
&\sum_{i=1}^{n-1} \lt(|\|\na  D_t^i p(t,\cdot)\||_{L^\iy} + |\|\na D_t^i {\rm div} u(t,\cdot)\||_{L^\iy}\rt)+ |\|\na^2 p  (t,\cdot)\||_{L^\iy}
 \\
& +  \|\na {\rm div} u(t,\cdot)\|_{L^\iy} + \|D_t p (t,\cdot)\|_{L^\iy}
 + \|D_t {\rm div} u (t,\cdot)\|_{L^\iy}
+ \|\na^2 \mathcal{T} (t,\cdot)\|_{L^\iy}     \\
& \le    C \lt({\rm Vol}\mathscr{D}_0 , K_0, \underline{\ea}_0^{-1}, \underline{\mathcal{T}}^{-1}, \overline{\mathcal{T}} ,  M_0, E_0(0),\cdots,E_{n+2}(0) \rt)
\end{split}\ee
for $ t\le T_{n2}({\rm Vol}\mathscr{D}_0 , K_0,  \underline{\ea}_0^{-1},  \underline{\mathcal{T}}^{-1}, \overline{\mathcal{T}} ,M_0,  E_0(0),\cdots,E_{n+2}(0) ) $.

Moreover, \ef{17/5/20-4} also holds for $t\le T_{n2}$.

\end{lem}

{\em Proof}. The proof consists of two cases of $n=2$ and $n=3$.

{\em Case 1}. Let n=2.
It follows from \ef{CL-5.34} and  \ef{17ta} that
\bee\label{}\begin{split}
 |\|\na \mathcal{T} \||_{L^\iy} \le & \|\na \Da \mathcal{T}\| + C\lt(K, K_1, |\|\ta\||, {\rm Vol} \Oa\rt) \|\Da \mathcal{T}\| \\
 \le & \| \na \Da \mathcal{T} \| + C\lt(K,   K_1, |\| (\na_N p)^{-1}\||_{L^\iy},  |\| \na_N p \||_{L^\iy}, E_2, {\rm Vol} \Oa \rt) \| \Da \mathcal{T} \|,
\end{split}\eee
which, together with \ef{e7},  Lemma \ref{17lem01}, \ef{w4.10}, \ef{17/5/20-3}  and \ef{17/5/20-4}, gives that for $ t\le T_{21}$,
\be\label{17/5/21-1}
|\|\na \mathcal{T}(t,\cdot)\||_{L^\iy}\le C_1\lt({\rm Vol}\mathscr{D}_0 , K_0, \underline{\ea}_0^{-1}, M_0 ,  E_0(0),\cdots,E_{2}(0) \rt) .
\ee
It follows from \ef{CL-A.15} and  \ef{CL-A.8}  that
\bee\label{ }\begin{split}
 \|D_t p\|_{L^\iy}
\le   C(K_1)\sum_{i=0}^2 \|\na^i D_t p \|   \ \  {\rm and } \ \   |\|\na  u \||_{L^\iy}  \le & C(K_1)  \sum_{i=1}^2 |\|\na^i u \||,
\end{split}\eee
which, together with  Lemmas \ref{17lem2} and \ref{17lem3}, \ef{w4.10}, \ef{17/5/20-3}  and \ef{17/5/20-4}, gives  that for $ t\le T_{21}$,
\be\label{17/5/21-2}
\|D_t p (t,\cdot)\|_{L^\iy} + |\|\na u(t,\cdot)\||_{L^\iy}\le C_2\lt({\rm Vol}\mathscr{D}_0 , K_0, \underline{\mathcal{T}}^{-1}, \overline{\mathcal{T}} ,  M_0, E_0(0),\cdots,E_{3}(0) \rt).
\ee
Similarly, we have that for $ t\le T_{21}$,
\be\label{17/5/21-3}\begin{split}
&  |\|\na {\rm div} u(t,\cdot)\||_{L^\iy} + |\|\na D_t p(t,\cdot)\||_{L^\iy} +  \|\na {\rm div} u(t,\cdot)\|_{L^\iy} + \|D_t {\rm div} u (t,\cdot)\|_{L^\iy}   \\
& \quad + |\|\na^2 p  (t,\cdot)\||_{L^\iy}  \le C_3\lt({\rm Vol}\mathscr{D}_0 , K_0, \underline{\ea}_0^{-1}, \underline{\mathcal{T}}^{-1}, \overline{\mathcal{T}} ,  M_0, E_0(0),\cdots,E_{3}(0) \rt),
\end{split}\ee
\be\label{17/5/21-4}\begin{split}
&\|\na^2 \mathcal{T} (t,\cdot)\|_{L^\iy} + |\|\na D_t {\rm div} u(t,\cdot)\||_{L^\iy}  \\
& \quad \le C_4\lt({\rm Vol}\mathscr{D}_0 , K_0, \underline{\ea}_0^{-1}, \underline{\mathcal{T}}^{-1}, \overline{\mathcal{T}} ,  M_0, E_0(0),\cdots,E_{4}(0) \rt) .
\end{split}\ee
Indeed, the bounds for  $|\|\na {\rm div} u \||_{L^\iy}$  and  $|\|\na D_t p \||_{L^\iy}$ in \ef{17/5/21-3} follows from  \ef{CL-5.34},  \ef{17ta}, \ef{4/24-1}, \ef{4/25-2}, Lemmas \ref{17lem01}-\ref{17lem3},  the bounds just obtained for $|\|\na \mathcal{T} \||_{L^\iy}$ and  $\|D_t p\|_{L^\iy}$, \ef{w4.10}, \ef{17/5/20-3}   and \ef{17/5/20-4}. The bounds for  $\|\na {\rm div} u \|_{L^\iy}$,   $\|D_t {\rm div} u  \|_{L^\iy}$   and $|\|\na^2 p  \||_{L^\iy}$  in \ef{17/5/21-3} follows from \ef{CL-A.15},  \ef{CL-A.8}, Lemmas \ref{17lem01}-\ref{17lem3},  the bounds just obtained for $|\|\na \mathcal{T} \||_{L^\iy}$ and  $|\|\na {\rm div} u \||_{L^\iy}$, \ef{w4.10}, \ef{17/5/20-3}   and \ef{17/5/20-4}.
The bound for $ \|\na^2 \mathcal{T} \|_{L^\iy}$ in \ef{17/5/21-4} follows from \ef{CL-A.15}, Lemmas \ref{17lem01}-\ref{17lem4}, \ef{w4.10}, \ef{17/5/20-3}   and \ef{17/5/20-4}.
The bound for $ |\|\na D_t {\rm div} u \||_{L^\iy} $
in \ef{17/5/21-4} follows from  \ef{CL-5.34},  \ef{17ta},  \ef{4/24-4},  Lemmas \ref{17lem01}-\ref{17lem4},  the bounds  just obtained for   $\|D_t {\rm div} u  \|_{L^\iy}$, $|\|\na \mathcal{T} \||_{L^\iy}$ and  $|\|\na {\rm div} u \||_{L^\iy}$, \ef{w4.10}, \ef{17/5/20-3}   and \ef{17/5/20-4}.

So, it follows from \ef{17/5/21-1}-\ef{17/5/21-4} that \ef{17/5/22}  holds for $t\le T_{22}$ for some continuous function $ T_{22}({\rm Vol}\mathscr{D}_0 , K_0,  \underline{\ea}_0^{-1},  \underline{\mathcal{T}}^{-1}, \overline{\mathcal{T}} ,M_0,  E_0(0),\cdots,E_{4}(0) ) >0$, by choosing  $L=2\sum_{i=3}^4 C_i$, $M=4M_0+2\sum_{i=1}^3 C_i$, and $\widetilde{M}=2\sum_{i=2}^4 C_i$ in the continuous function $T_{21}$ given by Lemma \ref{17lem6}.
Clearly, \ef{17/5/20-4}  holds for $t\le T_{22}$.

{\em Case 2}. Let $n=3$. In this case, we can use the arguments  similar to the way which we dealt with case 1 to obtain that  for $ t\le T_{31}$,
\be\label{17/5/23}\begin{split}
&\|D_t p (t,\cdot)\|_{L^\iy} \le C_5\lt({\rm Vol}\mathscr{D}_0 , K_0, \underline{\mathcal{T}}^{-1}, \overline{\mathcal{T}} ,  M_0, E_0(0),\cdots,E_{3}(0) \rt),\\
& |\|\na \mathcal{T}(t,\cdot)\||_{L^\iy} +   \|D_t {\rm div} u (t,\cdot)\|_{L^\iy}   \le C_6\lt({\rm Vol}\mathscr{D}_0 , K_0, \underline{\ea}_0^{-1}, \underline{\mathcal{T}}^{-1}, \overline{\mathcal{T}} ,  M_0, E_0(0),\cdots,E_{3}(0) \rt),\\
& \|\na^2 \mathcal{T} (t,\cdot)\|_{L^\iy}    \le C_7\lt({\rm Vol}\mathscr{D}_0 , K_0, \underline{\ea}_0^{-1}, \underline{\mathcal{T}}^{-1}, \overline{\mathcal{T}} ,  M_0, E_0(0),\cdots,E_{4}(0) \rt) .
\end{split}\ee
It follows from \ef{CL-5.29}, \ef{hb4}  and \ef{bdry}  that for any $\da\in (0,1]$,
\begin{align*}
 \|\na^3  {\rm div} u\|  +  |\| \na^2   {\rm div} u \||
 \le & \da |\|\Pi \na^3  {\rm div} u\|| + C(\da^{-1}, K, {\rm Vol}\Oa)\sum_{s=0}^1\|\na^s\Da {\rm div} u \| \notag\\
\le  & \da \lt(     |\| \na_N {\rm div} u \||_{L^\iy}  |\|\overline{\na} \ta \||  + 3K |\|  \na^2 {\rm div} u\|| + 2K^2 |\|  \na {\rm div} u \||    \rt) \notag\\
& + C(\da^{-1}, K, {\rm Vol}\Oa)\sum_{s=0}^1\|\na^s\Da {\rm div} u \|
      ,     \notag
\end{align*}
which implies, by choosing  $\da = \min\{(6K)^{-1},1\}$, that
\begin{align}
  \|\na^3  {\rm div} u\|  + 2^{-1} |\| \na^2   {\rm div} u \||
\le  &   |\| \na_N {\rm div} u  \||_{L^\iy} |\|\overline{\na} \ta \||  + C(K)   |\|  \na  {\rm div} u \||    \notag\\
&+  C(  K, {\rm Vol}\Oa)\sum_{s=0}^1\|\na^s\Da {\rm div} u \|.  \notag
\end{align}
This, together with \ef{CL-5.16}, gives
\be \label{noting}
  \|\na^4   u\|^2 \le    C\|\na^3  {\rm div} u\|^2  + C (\|\mathcal{T}\|_{L^\iy}+ 1)E_4
   \le    C  |\| \na_N {\rm div} u  \||_{L^\iy}^2 |\|\overline{\na} \ta \||^2 + L_1 ,
\ee
where
$$L_1=C(K)   |\|  \na  {\rm div} u \||^2
 +  C(  K, {\rm Vol}\Oa)\sum_{s=0}^1\|\na^s\Da {\rm div} u \|^2  + C (\|\mathcal{T}\|_{L^\iy}+ 1)E_4.$$
It follows from \ef{CL-5.34},  \ef{4/24-1} and \ef{noting} that for any $\da\in (0,1]$,
\begin{align}
|\| \na  {\rm div} u  \||_{L^\iy}^2 \le &\da \|\na^2  \Da {\rm div} u  \|^2 +  C  (\da^{-1},K,K_1,|\|\overline{\na}\ta\||, {\rm Vol}\Oa )  \sum_{s=0}^1\|\na^s\Da {\rm div} u \|^2 \notag \\
\le & \da   \|  \mathcal{T}^{-1}  \|_{L^\iy}^2 C(K_1, \|\na u\|_{L^\iy}, \|\na \mathcal{T}\|_{L^\iy} ) \|\na^4 u \|^2 + L_2 \notag\\
\le & \da   \|  \mathcal{T}^{-1}  \|_{L^\iy}^2 C(K_1, \|\na u\|_{L^\iy}, \|\na \mathcal{T}\|_{L^\iy} )|\| \na_N {\rm div} u  \||_{L^\iy}^2 |\|\overline{\na} \ta \||^2  + L_3, \label{17/5/23-1}
\end{align}
where
\begin{align*}
L_2= &C  (\da^{-1},K,K_1,|\|\overline{\na}\ta\||, {\rm Vol}\Oa )  \sum_{s=0}^1\|\na^s\Da {\rm div} u \|^2 + \da   \|  \mathcal{T}^{-1}  \|_{L^\iy}^2   C     \| \na^2 D_t{\rm div} u \| ^2 \\
&+ \da   \|  \mathcal{T}^{-1}  \|_{L^\iy}^2 C(K_1,\|\na  u\|_{L^\iy},
  \|\na \mathcal{T} \|_{L^\iy})\lt(\sum_{i=1}^{3} \|\na^i u\|^2 +\sum_{i=1}^{4}   \|\na^i \mathcal{T} \|^2   \rt) ,
\end{align*}
and
$$L_3=  \da   \|  \mathcal{T}^{-1}  \|_{L^\iy}^2 C(K_1, \|\na u\|_{L^\iy}, \|\na \mathcal{T}\|_{L^\iy} ) L_1 + L_2.$$
By choosing suitable small $\da $ in  \ef{17/5/23-1}, and using
Lemmas \ref{17lem01}-\ref{17lem4}, the bound just obtained for $|\|\na\mathcal{T}\||_{L^\iy}$  in \ef{17/5/23}, \ef{w4.10} and \ef{17/5/20-4}, we have  that for $ t\le T_{31}$,
\be\label{17/5/23-2}
|\| \na  {\rm div} u (t,\cdot)  \||_{L^\iy}   \le C_8\lt({\rm Vol}\mathscr{D}_0 , K_0, \underline{\ea}_0^{-1}, \underline{\mathcal{T}}^{-1}, \overline{\mathcal{T}} ,  M_0, E_0(0),\cdots,E_{4}(0) \rt).
\ee
In a similar way to deriving \ef{17/5/23-2}, we have, using \ef{4/25-2}, \ef{4/24-4}, \ef{174/10} and \ef{174/11}, that for $ t\le T_{31}$,
\begin{align}
&|\| \na  D_t p (t,\cdot)  \||_{L^\iy}   \le C_9\lt({\rm Vol}\mathscr{D}_0 , K_0, \underline{\ea}_0^{-1}, \underline{\mathcal{T}}^{-1}, \overline{\mathcal{T}} ,  M_0, E_0(0),\cdots,E_{4}(0) \rt),\label{17/5/23-3}\\
& |\| \na  D_t  {\rm div} u (t,\cdot)  \||_{L^\iy} + |\| \na  D_t^2 p  (t,\cdot)  \||_{L^\iy}   + |\| \na  D_t^2    {\rm div} u (t,\cdot)  \||_{L^\iy} \notag\\
 & \qquad \le C_{10}\lt({\rm Vol}\mathscr{D}_0 , K_0, \underline{\ea}_0^{-1}, \underline{\mathcal{T}}^{-1}, \overline{\mathcal{T}} ,  M_0, E_0(0),\cdots,E_{5}(0) \rt) . \label{17/5/23-4}
\end{align}
With these bounds, we can obtain, in the same manner as the case of $n=2$, that for $ t\le T_{31}$,
\be\label{17/5/23-5}
 \| \na  {\rm div} u (t,\cdot)  \|_{L^\iy} + |\|\na u(t,\cdot)\||_{L^\iy}   \le C_{11}\lt({\rm Vol}\mathscr{D}_0 , K_0, \underline{\ea}_0^{-1}, \underline{\mathcal{T}}^{-1}, \overline{\mathcal{T}} ,  M_0, E_0(0),\cdots,E_{4}(0) \rt),
 \ee
 \be\label{17/5/23-6}
  |\| \na^2 p    (t,\cdot)  \||_{L^\iy}     \le C_{12}\lt({\rm Vol}\mathscr{D}_0 , K_0, \underline{\ea}_0^{-1}, \underline{\mathcal{T}}^{-1}, \overline{\mathcal{T}} ,  M_0, E_0(0),\cdots,E_{5}(0) \rt).
  \ee

It is produced  from \ef{17/5/23}, \ef{17/5/23-1}-\ef{17/5/23-6} that \ef{17/5/22}  holds for $t\le T_{32}$ for some continuous function $ T_{32}({\rm Vol}\mathscr{D}_0 , K_0,  \underline{\ea}_0^{-1},  \underline{\mathcal{T}}^{-1}, \overline{\mathcal{T}} ,M_0,  E_0(0),\cdots,E_{5}(0) ) >0$, by choosing  $L= 2(C_9+C_{10}+C_{12})$, $M=4 M_0+ 2(C_6+C_8+C_{11}) $, and $\widetilde{M}= 2\sum_{i=5}^7 C_i $ in the continuous function $T_{31}$ given by Lemma \ref{17lem6}. Clearly, \ef{17/5/20-4}  holds for $t\le T_{32}$.
\hfill $\Box$

\vskip 0.5cm

\noindent{\bf Acknowledgements} Luo's research was  supported in part by a GRF grant CityU 11303616 of RGC (Hong Kong). Zeng's research was supported in part by NSFC  Grant 11671225, and the Center of Mathematical Sciences and Applications, Harvard University.

\bibliographystyle{plain}

\noindent {Tao Luo}\\
Department of Mathematics\\
City University of Hong Kong,\\
Tat Chee Ave. Hong Kong\\
Email: taoluo@cityu.edu.hk\\

\noindent Huihui Zeng\\
Yau Mathematical Sciences Center\\
Tsinghua University\\
Beijing, 100084, China;\\
Center of Mathematical Sciences and Applications\\
 Harvard University\\
  Cambridge, MA 02318, USA.\\
Email: hhzeng@mail.tsinghua.edu.cn

\end{document}